\documentclass[12pt,letterpaper,table]{article}
\usepackage{appendix}
\usepackage{mathtools}
\usepackage{graphicx}
\usepackage{setspace}
\usepackage{listings}
\usepackage{xcolor}

\usepackage[margin=0.65in]{geometry}
\usepackage[version=3]{mhchem} 
\usepackage[english]{babel}
\usepackage{chemfig}
\usepackage{comment}
\usepackage{hyperref}

\usepackage{subcaption}

\usepackage[utf8]{inputenc}
\usepackage{amsmath}
\usepackage{amsfonts}
\usepackage{amssymb}
\usepackage{mathtools}
\usepackage{graphicx}
\usepackage{algpseudocode}
\usepackage{algorithm}
\usepackage{tikzsymbols}
\usepackage{subcaption}
\usepackage{cancel}
\usepackage{mathrsfs}
\usepackage{cases}
\usepackage{pdflscape}
\usepackage{authblk}

\usepackage{float}
\usepackage[font=small,labelfont=bf]{caption}
\emergencystretch=\maxdimen
\usepackage{lipsum}
\usepackage{listings}
\usepackage{float}
\usepackage{longtable}
\date{\vspace{-5ex}}

\newenvironment{proof}{\paragraph{Proof:}}{\hfill$\square$}

\usepackage{hyperref}

\usepackage{amsmath}

\newtheorem{theorem}{Theorem}[section]
\newtheorem{lemma}{Lemma}[section]
\newtheorem{remark}{Remark}[section]

\newtheorem{proposition}{Proposition}[section]

\definecolor{codegreen}{rgb}{0,0.6,0}
\definecolor{codegray}{rgb}{0.5,0.5,0.5}
\definecolor{codepurple}{rgb}{0.58,0,0.82}
\definecolor{backcolour}{rgb}{0.95,0.95,0.92}

\lstdefinestyle{mystyle}{
    backgroundcolor=\color{backcolour},   
    commentstyle=\color{codegreen},
    keywordstyle=\color{magenta},
    numberstyle=\tiny\color{codegray},
    stringstyle=\color{codepurple},
    basicstyle=\ttfamily\footnotesize,
    breakatwhitespace=false,         
    breaklines=true,                 
    captionpos=b,                    
    keepspaces=true,                 
    numbers=left,                    
    numbersep=5pt,                  
    showspaces=false,                
    showstringspaces=false,
    showtabs=false,                  
    tabsize=2
}

 \usepackage{tikz} 
 \usetikzlibrary{arrows}
 \usepackage{float}

 \usepackage{subcaption}

\definecolor{boxcolor_tumour}{HTML}{B9DCFF}
\definecolor{boxcolor_lymph}{HTML}{90EE90}
\usepackage[utf8]{inputenc}

\title{\centering On the Treatment of Melanoma: \\ A Mathematical Model of Oncolytic Virotherapy}

\author[]{Tedi Ramaj}
\author[]{Xingfu Zou}

\affil[]{{\small Department of Mathematics, Western University, London, On. Canada}}

\begin{document}

\maketitle

\begin{abstract}
    We develop and analyze a mathematical model of oncolytic virotherapy in the treatment of melanoma. We begin with a special, local case of the model, in which we consider the dynamics of the tumour cells in the presence of an oncolytic virus at the primary tumour site. We then consider the more general regional model, in which we incorporate a linear network of lymph nodes through which the tumour cells and the oncolytic virus may spread. The modelling also considers the impact of hypoxia on the disease dynamics. The modelling takes into account both the effects of hypoxia on tumour growth and spreading, as well as the impact of hypoxia on oncolytic virotherapy as a treatment modality. We find that oxygen-rich environments are favourable for the use of adenoviruses as oncolytic agents, potentially suggesting the use of complementary external oxygenation as a key aspect of treatment. Furthermore, the delicate balance between a virus' infection capabilities and its oncolytic capabilities should be considered when engineering an oncolytic virus. If the virus is too potent at killing tumour cells while not being sufficiently effective at infecting them, the infected tumour cells are destroyed faster than they are able to infect additional tumour cells, leading less favourable clinical results. Numerical simulations are performed in order to support the analytic results and to further investigate the impact of various parameters on the outcomes of treatment. Our modelling provides further evidence indicating the importance of three key factors in treatment outcomes: tumour microenvironment oxygen concentration, viral infection rates, and viral oncolysis rates. The numerical results also provide some estimates on these key model parameters which may be useful in the engineering of oncolytic adenoviruses. 
\end{abstract}

\vspace{5mm}

\textbf{Keywords}: mathematical modelling, oncolytic virotherapy, disease dynamics, differential equation modelling

\section{Introduction}

Melanoma is considered the most deadly type of skin cancer. Melanoma begins in melanocytes - the cells responsible for producing melanin - and can develop in various parts of the body \cite{melanocyte_paper, melanocyte_paper2}.  Melanoma is the fifth most common cancer in adults in the United States \cite{melanoma_stats1}.  While melanoma rates have been steadily rising, mortality has not followed this same trend. This decreased mortality is attributed to various factors such as early detection, increased protection against UV radiation, and improvements in treatment \cite{melanoma_stats2}.  Metastatic melanoma continues to be a major issue contributing the cancer mortality, due to the increased difficulty of treating the disease once it has spread beyond its original site \cite{melanoma_stats3}. Various forms of therapy, including chemotherapy, immunotherapy, and radiotherapy are used in the treatment of advanced melanoma. Developing new forms of therapy and enhancing existing therapy is always desirable in increasing survival rates of the disease. 

Oncolytic virotherapy is a method of cancer treatment in which viruses are used to selectively infect and destroy cancer cells via a variety of direct and indirect mechanisms, while leaving surrounding healthy cells unharmed \cite{ OV_paper2, OV_paper1}.  These viruses are called oncolytic viruses (OVs). {These therapeutics include both genetically modified viruses and non-modified viruses, such as live attenuated viruses (i.e., the measles virus \cite{measles_virus_paper}).  } The genetically modified herpes simplex virus Talimogene laherparepvec (T-VEC) has been used in clinical trials to treat inoperable melanoma \cite{TVECpaper1, TVECpaper2}.  The treatment is often performed in combination with other therapies, such as being followed up with the use of adjuvant radiotherapy. The oncolytic virus is typically administered via direct subcutaneous injection into the lesion \cite{TVECpaper3}. The idea is for the virus to selectively infect cancer cells and use them to replicate and perform oncolysis to destroy the neoplasm. The viral infection may also destroy the cancer cells through indirect mechanisms such as activating the immune system and aiding the immune response against the cancer cells \cite{ OV_paper3, OV_paper2, OV_paper1}. 

Other OVs which have been studied (not necessarily in melanoma trials, but in the context of other cancers, such as colorectal cancer) include the adenoviruses ONYX-015 and ZD55-IL-24 \cite{hu_adenovirus}. {ZD55-IL-24's primary mechanism of action is through inducing a systematic anti-tumour cell immune response. There is also evidence that this virus may inhibit tumour cell growth by inhibiting angiogenesis, as was previously observed in an immuno-competent mouse model \cite{hu_adenovirus}. Such immune-mediated effects of viruses like ZD55-IL-24 are more established in the existing literature. A mechanism of action of ONYX-015 involves replication and lysis of tumour cells that are p-53 deficient \cite{ ONYX015_carlapaper, ONYX015_RiesPaper}. The importance of this direct mechanism of action of is currently under debate \cite{davola} and is considered in this present paper via mechanistic modelling.  }

Mathematical modelling of cancer treatment has seen widespread use in the last few decades. These models frequently take the form of ODE, PDE, and delay models in the continuous setting. By studying the effect of disease treatment from a quantitative perspective, based on biological and physical mechanistic modelling, new insights may be obtained to guide future treatment direction. This type of modelling has also been used to study the treatment of cancer via oncolytic virotherapy \cite{mathmodel_onc1, mathmodel_onc2}.  The recent work of Wang et al \cite{mathmodel_onc3}. in mathematical modelling of virotherapy as a treatment modality for melanoma, the models were able to provide insights concerning virus treatment thresholds as well as how immunosupperssive drugs may work in tandem with OVs.  The work by Urenda-Cazares et al. examined the use of OVs in combination with chemotherapy to treat glioma. As a result of these types of models, some results were obtained on how to optimize treatment in a clinical setting \cite{mathmodel_onc4}. 

In this paper, we model the effect of hypoxic environments on oncolytic virotherapy treatment through the use of ordinary differential equation (ODE) modelling. Hypoxic refers to oxygen-poor environments. Typically, viruses which are more efficient at infecting cells in oxygen-rich environments tend to lose their infectivity under hypoxic conditions \cite{viral_friend_foe}. This is particularly true of adenoviruses such as ONYX-015 \cite{onyx_adeno1}. Hypoxia has a negative effect on the efficacy of OVs as well as any adjuvant radiotherapy which may be administered \cite{hypoxiapaper1, hypoxiapaper2}.  In the context of melanoma treatment, hypoxic environments can inhibit the action of OVs, such as their ability to infect cancer cells and their ability to induce the death of cancer cells. {Due to the lack of dynamical modelling of this phenomenon, we explore the relationship between tumour microenvironment oxygen concentration and the efficacy of the OV with the objective of contributing to the existing oncology literature from a quantitative perspective. The application of mathematical modelling can capture some elements of the complex interplay between oxygen concentration conditions and OV efficacy. } In our model, we study the effect of oxygen concentration when the OV is applied directly to the primary lesion. More specifically, we study the impact which parameters such as the infectiousness of the OV on the efficacy of the treatment under different oxygen conditions.

The structure of this paper is organized as follows. In Section 2, we formulate an ODE model and give the assumptions on our functions. We refer to this as our local model, since we are studying the effect of OV directly on the primary tumour. We also perform non-dimensionalization of the model for the purposes of mathematical analysis. We explain the meaning of our model in terms of the biological context. In Section 3, we perform an analysis of the local model. This includes proofs on  the well-posedness results. In Subsection 3.1, we first look at the case where we do not take into account the oxygen concentration dependence. In Subsection 3.2, we look at the case of oxygen concentration dependence. We perform an analysis of the stability of the relevant steady states of our system. In Section 4, we perform numerical simulations and give biological interpretations of these results. In Section 5, we extend our model to a regional model, where we take into account the movement of tumour cells into the surrounding lymph nodes. In Section 6, we perform numerical simulations on the regional model. We complete this paper with some conclusions and discuss possible directions for future work in Section 7.

\section{Local oncolytic virotherapy model}

We begin by considering a melanoma tumour, initially consisting of some initial quantity of proliferating  tumour cells. At this initial point in time, a localized treatment of oncolytic virotherapy begins at the  tumour site, by introducing the OV via direct injection into the lesion. We consider the {use of a virus with oncolytic and replication rates down-regulated by hypoxia. Such a virus shares these features with adenoviruses. The rationale for considering adenoviruses (or OVs with similar hypoxia down-regulating properties as adenoviruses), comes from this consideration of the effects of hypoxia on the action of the OV. Namely, while hypoxic tumour microenvironments reduce the efficacy of adenoviruses, they also promote melanoma tumour progression \cite{melanoma_progression_paper}. One of the goals of this present work is to mathematically capture and model the dynamics of OVs under the same unfavourable hypoxic conditions which typically have an inverse (favourable) impact on tumour progression.  Indeed, the  modelling presented in this paper is not only limited to adenoviruses, but to any virus which experiences similar down-regulation in hypoxic tumour microenvironments. } The OV then proceeds to infect the  tumour cells. The model consists of three variables, the density of uninfected  tumour cells, the density of infected  tumour cells, and oxygen concentration, at time $t$, respectively represented by $u(t), n(t),$ and $c(t)$. Then, we have the following model:
\begingroup
\addtolength{\jot}{7pt}
\begin{align}
\dfrac{\mathrm{d}u}{\mathrm{d}t} &= r_1 u \left( 1- \dfrac{u+n}{K}  \right) - \dfrac{\theta (c) nu}{\alpha + n}, \label{3eq1}\\
\dfrac{\mathrm{d}n}{\mathrm{d}t} &= r_2 n \left( 1 - \dfrac{u+n}{K}  \right) + \dfrac{\theta (c) nu}{\alpha + n} - \gamma (c) n, \label{3eq2}  \\
\dfrac{\mathrm{d}c}{\mathrm{d}t} &= \phi - \beta c - q_1uc - q_2nc.  \label{3eq3}
\end{align}
\endgroup
Note that we are considering cell-to-cell infections, which have been observed as a mode of infection used by oncolytic viruses \cite{cell_cell_virus_evidence}. {Oncolytic adenoviruses which exhibit cell-to-cell spreading, such as VRX-009, have also previously been constructed \cite{adenovirus_cell_to_cell}. Importantly, while VRX-009 was not tested as a treatment modality for melanoma, its production provides a proof of concept of the idea of an oncolytic adenovirus with a cell-to-cell spreading mechanism. Our work hence provides a theoretical modelling framework for cell-to-cell spreading of adenoviruses (or adenovirus-like OVs) which may infect melanoma cancer cells.}  Previous mathematical models of cell-to-cell viral infection made use of a mass-action-like terms to represent infection \cite{ Webb_cell_to_cell, zou_cell_to_cell} and we adopt a similar approach in our model, but  with the infection mechanism of a Holling type II functional response function. {We make this consideration to model the saturating effect of melanoma cells that have already been infected by the OV.} 

We prescribe the initial conditions $u(0) = u_0, n(0) = n_0$ and $c(0) = c_0$ to be non-negative quantities. We assume that both classes of  tumour cells exhibit logistic growth. The carrying capacity of the  tumour cells is given by $K$ and the growth rates of the uninfected  tumour cells and the infected  tumour cells are given by $r_1$ and $r_2$, respectively. We further assume that $r_1 > r_2$ to reflect  that the infected  tumour cells are less effective at proliferating due their cell machinery being hijacked by the OV. Following the approach of \cite{oxygen_con_paper}, we use mass-action terms to express the oxygen consumption by the tumour cells. To that end, the parameters $q_1$ and $q_2$ give the oxygen consumption rate by the  uninfected tumour cells and the infected tumour cells, respectively. The rate of oxygenation, assumed constant (due to having some control over this parameter, i.e., through certain therapies \cite{tibbles}), is given by $\phi$ and the rate of oxygen consumption by surrounding \textit{non-cancerous cells (or healthy cells)} is given by $\beta$.  

We use a Hill function to represent the transition of a  tumour cell from uninfected by an OV to infected by an OV. The parameter $\theta  \in \mathrm{C}^1 (\mathbb{R_+})$ represents the virus infection rate, which is dependent on available oxygen concentration. The other oxygen dependent parameter $\gamma \in \mathrm{C}^1 (\mathbb{R_+}) $ is the virus-induced death rate of the infected  tumour cells. Note that the terms \textit{virus-induced death rate} and \textit{oncolysis rate} are used interchangeably. The adenovirus is inhibited by a hypoxic environment and hence we assume that as oxygen concentration is locally decreased, the OV will become less effective, both in infecting the  tumour cells and inducing  tumour cell death \cite{onyx_adeno1}. Hence, we have the following conditions on $\theta(c)$ and $\gamma(c)$:
\begin{equation}
\begin{dcases}
\theta'(c) \geq 0, \quad \gamma'(c) \geq 0, \quad \text{for} \quad c \in (0, \infty), \\[5pt]
\theta(0) = \theta_0 \geq 0,  \quad \gamma(0) = \gamma_0 \geq 0, \\[5pt]
\lim_{c \rightarrow \infty} \theta(c) = \theta_\infty > \theta_0, \quad \lim_{c \rightarrow \infty} \gamma (c) = \gamma_\infty > \gamma_0,
\end{dcases} \label{oxygen conditions}
\end{equation}
where $\theta_\infty$ and $\gamma_\infty$ give the OV efficacy in response to high oxygen environments. Note that in hypoxic environments, oncolytic virotherapy will not be as efficient as an adenovirus is being used. 

We non-dimensionalize the model by making the following substitutions:
\begin{equation*}
x:= \dfrac{u}{K}, \quad y:= \dfrac{n}{K}, \quad z := \dfrac{\beta c}{\phi}, \quad \tau := r_1 t.
\end{equation*}
Then we obtain the system,
\begingroup
\allowdisplaybreaks
\addtolength{\jot}{7pt}
\begin{align}
\dfrac{\mathrm{d}x}{\mathrm{d}\tau} &= x(1-x-y) - \dfrac{\hat{\theta} (z) xy}{\hat{\alpha} + y}, \\
\dfrac{\mathrm{d}y}{\mathrm{d}\tau} &= ry(1-x-y) + \dfrac{\hat{\theta} (z) xy}{\hat{\alpha} + y} - \hat{\gamma} (z) y,  \\
\dfrac{\mathrm{d}z}{\mathrm{d}\tau} &= \hat{\beta} (1-z) - \hat{q}_1 xz - \hat{q}_2 yz,
\end{align}
\endgroup
where we define
\begin{equation*}
\hat{\theta} (z) := \dfrac{1}{r_1} \theta \left(\dfrac{\phi z}{\beta} \right), \quad \hat{\gamma} (z) := \dfrac{1}{r_1} \gamma\left( \dfrac{\phi z}{\beta} \right),
\end{equation*}
\begin{equation*}
\hat{r} := \dfrac{r_2}{r_1}, \quad \hat{\alpha} := \dfrac{\alpha}{K}, \quad \hat{\beta} := \dfrac{\beta}{r_1}, \quad \hat{q}_1 :=  \dfrac{q_1 K}{r_1}, \quad \hat{q}_2 := \dfrac{q_2 K}{r_1}.
\end{equation*}
Note that the properties of $\theta(c)$ and $\gamma(c)$ given in (\ref{oxygen conditions}) are preserved by $\hat{\theta}(z)$ and $\hat{\gamma}(z)$, respectively, up to some scaling. The most notable change is in the long-term behavior: $\hat{\theta}$ will approach $\hat{\theta}_\infty := \theta_\infty /r_1$ and $\hat{\gamma}$ will approach $\hat{\gamma}_\infty := \gamma_\infty / r_1$ as $z \rightarrow \infty$. We now drop the tilde and replace $\tau$ with $t$ for notational convenience and hence, for the subsequent analysis, we consider the following model:
\begingroup
\addtolength{\jot}{7pt}
\begin{align}
\dfrac{\mathrm{d}x}{\mathrm{d}t} &= x(1-x-y) - \dfrac{{\theta} (z) xy}{{\alpha} + y}, \label{e1} \\
\dfrac{\mathrm{d}y}{\mathrm{d}t} &= ry(1-x-y) + \dfrac{{\theta} (z) xy}{{\alpha} + y} - {\gamma} (z) y, \label{e2}  \\
\dfrac{\mathrm{d}z}{\mathrm{d}t} &= \beta (1-z) - q_1 xz - q_2 yz, \label{e3}
\end{align}
\endgroup
with non-negative initial conditions:
\begin{equation}
x(0) = x_0 \geq 0, \quad y(0) = y_0 \geq 0, \quad z(0) = z_0 \geq 0. \label{IC}
\end{equation}
The functions $\theta$ and $\gamma$ once again have the properties given in (\ref{oxygen conditions}). 
In Section 3, we perform a mathematical analysis of the rescaled model (\ref{e1}) - (\ref{e3}) to explore some predictions related to the effect of available oxygen concentration on the OV treatment. 

\section{Analysis of the local model}

We begin by considering the well-posedness of the model. Existence and uniqueness of the solution of (\ref{e1}) - (\ref{e3}), subject to initial conditions (\ref{IC}), follow from the elementary theory of ODEs. We consider the solution of this initial value problem, $(x(t), y(t), z(t)) \in \mathbb{R}^3$. Since the variables represent densities and concentration of physical quantities, the system must remain non-negative for all $t \geq 0$. We begin with equation (\ref{e1}). From this equation, it follows that
\begin{equation*}
x(t) = x_0 \cdot \exp \left[ \int_0^t  \left( 1 - x(s) - y(s) - \dfrac{\theta(z(s))y(s)}{\alpha + y(s)} \right)  \mathrm{d}s   \right],
\end{equation*}
and so, $x(t) \geq 0$ for all $t \geq 0$. Similarly, it follows from equation (\ref{e2}) that
\begin{equation*}
y(t) = y_0 \cdot \exp \left[ \int_0^t  \left( r(1-x-y) + \dfrac{\theta(z(s)) x(s)}{\alpha + y(s)} - \gamma(z(s))y(s)   \right)  \mathrm{d}s  \right].
\end{equation*}
Therefore, $y(t) \geq 0$ for all $t \geq 0$. Finally, equation (\ref{e3}) gives
\begin{equation*}
z(t) = z_0 \cdot\exp \left({-\int_0^t (\beta + q_1 x(s) + q_2 y(s)) \mathrm{d}s} \right) + \beta \int_0^t \exp \left({-\int_s^t (\beta + q_1 x(\xi) + q_2 y(\xi)) \mathrm{d} \xi}\right) \mathrm{d}s.
\end{equation*}
This shows that $z(t) \geq 0$ for all $t \geq 0$. In fact, if $t > 0$, then $z$ is strictly positive. 

Next, we address the boundedness of the solution. To this end, we apply a comparison argument. From equations (\ref{e1}) and (\ref{e3}), a solution of the system satisfies the inequalities
\begin{equation*}
\dfrac{\mathrm{d}x}{\mathrm{d}t} \leq x(1-x),  \quad \dfrac{\mathrm{d}z}{\mathrm{d}t} \leq \beta (1-z).
\end{equation*}
Then, it follows that
\begin{equation*}
\limsup_{t \rightarrow \infty} x(t) \leq 1, \quad \limsup_{t \rightarrow \infty} z(t) \leq 1.
\end{equation*}
Hence, $x(t)$ and $z(t)$ are bounded functions. Let $\bar{x}$ be an upper bound for $x(t)$, i.e., $x(t) \leq \bar{x}$ for all $t \geq 0.$ It then follows from equation $(\ref{e2})$ that
\begin{equation*}
\dfrac{\mathrm{d}y}{\mathrm{d}t} \leq ry(1-y) + \dfrac{\theta_\infty \bar{x} y}{\alpha + y} \implies \dfrac{\mathrm{d}y}{\mathrm{d}t} \leq ry(1-y) + \theta_\infty \bar{x}.
\end{equation*}
Therefore, by a comparison argument,
\begin{equation*}
\limsup_{t \rightarrow \infty} y(t) \leq \dfrac{r + \sqrt{r^2 + 4r \theta_\infty \bar{x}}}{2r},
\end{equation*}
which shows that $y(t)$ is a bounded function. 

Summarizing these results, we have the following theorem:
\begin{theorem}
The solution of the initial value problem (\ref{e1}) - (\ref{e3}), satisfying initial conditions (\ref{IC}), is non-negative and bounded. \label{theorem_pos_bounded}
\end{theorem}

\subsection{Dynamics of the local model -- case I: no oxygen dependence}

We consider first the case with no oxygen dependence. That is, we set $\theta(z) = \theta$ and $\gamma(z) = \gamma$, where $\theta$ and $\gamma$ are positive constants. In this case, system (\ref{e1}) - (\ref{e3}) reduces to the following two-variable system:
\begingroup
\addtolength{\jot}{7pt}
\begin{align}
\dfrac{\mathrm{d}x}{\mathrm{d}t} &= x(1-x-y) - \dfrac{{\theta} xy}{{\alpha} + y}, \label{e11} \\
\dfrac{\mathrm{d}y}{\mathrm{d}t} &= ry(1-x-y) + \dfrac{{\theta} xy}{{\alpha} + y} - {\gamma} y. \label{e22} 
\end{align}
\endgroup
If we consider system (\ref{e11}) - (\ref{e22}) over the region $(x,y) \in \mathbb{R}_+^2$, we can rule out the existence of non-constant periodic orbits.
\begin{proposition}
Consider system (\ref{e11}) - (\ref{e22}) over the region $\mathbb{R}_+^2$. There are no closed orbits contained entirely $\mathbb{R}_+^2$. \label{prop_dulac}
\end{proposition}

\begin{proof}
Let $S(x,y) = 1/(xy)$ for $x, y > 0$. Then,
\begin{equation*}
\dfrac{\partial}{\partial x} \left[ S(x,y) \left(  x(1-x-y) - \dfrac{{\theta} xy}{{\alpha} + y}  \right)   \right] + \dfrac{\partial}{\partial y} \left[ S(x,y) \left(  ry(1-x-y) + \dfrac{{\theta} xy}{{\alpha} + y} - {\gamma} y \right) \right]
\end{equation*}
may be computed to give
\begin{equation}
-\dfrac{1}{y} - \dfrac{r}{x} - \dfrac{\theta}{(\alpha + y)^2} < 0.
\end{equation}
Since this function does not change sign on $\mathbb{R}_+^2$, we conclude by the Dulac-Bendixson Theorem that there are no closed orbits contained entirely in $\mathbb{R}_+^2$. 
\end{proof}

Next, we determine the steady states of system (\ref{e11}) - (\ref{e22}) by solving the algebraic system
\begingroup
\addtolength{\jot}{7pt}
\begin{align}
x(1-x-y) - \dfrac{{\theta} xy}{{\alpha} + y} = 0, \label{e11_alg} \\
ry(1-x-y) + \dfrac{{\theta} xy}{{\alpha} + y} - {\gamma} y = 0. \label{e22_alg} 
\end{align}
\endgroup
It can be readily seen that $(x,y) = (0,0)$ and $(x,y) = (1,0)$ are solutions of this system for all parameter values. Another solution which may be easily seen is $(x,y) = (0, (r - \gamma)/r)$, which only exists if $r > \gamma$. It can be shown that the remaining steady states (if any exist) are determined by solving the system 
\begin{equation}
(\theta - r \theta + \gamma)y^2 + (\theta^2 + \alpha \theta + 2\alpha \gamma - r \alpha \theta - \theta)y + \alpha(\alpha \gamma - \theta) = 0, \label{alg1} 
\end{equation}
\begin{equation}
x = 1 - y - \dfrac{\theta y}{\alpha + y}.\label{alg2}
\end{equation}
We linearize the system at its steady states by first computing the Jacobian matrix  
\begin{equation*}
J(x,y) = \begin{pmatrix} 1 - 2x - y - \dfrac{\theta y}{\alpha + y} & -x - \dfrac{\theta \alpha x}{(\alpha + y)^2}    \\[18pt] -ry + \dfrac{\theta y}{\alpha + y} & r - rx - 2ry + \dfrac{\alpha \theta x}{(\alpha + y)^2} - \gamma \end{pmatrix}.
\end{equation*}
We begin with the assumption $r < \gamma$ in order to discount the steady state $(0,(r-\gamma)/r)$. Linearizing at the steady state $(0,0)$ gives
\begin{equation*}
J(0,0) = \begin{pmatrix} 1 & 0   \\ 0 & r - \gamma \end{pmatrix}.
\end{equation*}
By our assumption that $r < \gamma$, this steady state is a saddle. Linearizing the system at the steady state $(1,0)$ gives
\begin{equation*}
J(1,0) = \begin{pmatrix} -1 \ & \  -1 - \dfrac{\theta}{\alpha}  \\[10pt] 0 & \dfrac{\theta}{\alpha} - \gamma \end{pmatrix}.
\end{equation*}
From $J(1,0)$, we then conclude that $(1,0)$ is locally asymptotically stable if $\theta < \alpha \gamma$ and it is unstable if $\theta > \alpha \gamma$. 

If we now impose the additional assumption $\theta < \alpha \gamma$, then the system only contains two non-negative steady states: $(0,0)$ and $(1,0)$. To see this, we note that the left-hand side of equation (\ref{alg1}), as a function of $y$, is a convex parabola (since $r < 1$) with a positive constant term. The coefficient of the $y$ term is also positive, as
\begin{align*}
    \theta^2 + \alpha \theta + 2\alpha \gamma  - r \alpha \theta - \theta = \theta^2 + \alpha \gamma + \alpha \theta (1-r) + (\alpha \gamma - \theta) > \theta^2 + \alpha \gamma > 0.
\end{align*}
Therefore, the parabola has non-negative roots and system  (\ref{alg1}) - (\ref{alg2}) has no non-negative solutions. This shows that the only non-negative steady states are $(0,0)$ and $(1,0)$.

Next, we consider the case $\theta > \alpha \gamma$. In this case, an additional co-existence steady state, $(x_*, y_*)$, where $x_*, y_* > 0$ may be introduced if system (\ref{alg1}) - (\ref{alg2}) has positive solutions. It is clear to see that  the parabola on the left-hand side of equation (\ref{alg1}) is still convex but the constant term is now negative. Hence, this parabola has exactly one positive real root, $y_*$. Then $x_*$ may be obtained from equation (\ref{alg2}). In order for the steady state to be meaningful, we set $x_*$ must be positive, which is not the case for all values of the model parameters. We impose the following condition to ensure that $x_*$ is positive: 
\begin{equation}
    y_*^2 + (\alpha + \theta)y_* < 1, \label{existence_cond31}
\end{equation}
where
\begin{equation*}
    y_* = \dfrac{A + \sqrt{A^2 - 4 \alpha (\theta - r \theta + \gamma)(\alpha \gamma - \theta)}}{2(\theta - r \theta + \gamma)},
\end{equation*}
and
\begin{equation*}
    A = r \alpha \theta + \theta - \theta^2 - \alpha \theta - 2 \alpha \gamma.
\end{equation*}
We assume these conditions are satisfied so that the steady state $(x_*, y_*)$ exists and has positive coordinates. In fact, since $r < \gamma$, these conditions are satisfied as the existence of a stable (unique) positive steady state is ensured as a corollary of non-negativity of solutions, boundedness of solutions, and  Proposition \ref{prop_dulac}. 

As we will now show, the assumption $r < \gamma$ is not necessary for the stability of the steady state $(x_*, y_*)$ -- only its existence is necessary. If it exists, linearizing at this steady state gives the matrix
\begin{equation*}
    J(x_*, y_*) = \begin{pmatrix} -x_* \quad & -x_* - \dfrac{\theta \alpha x_*}{(\alpha + y_*)^2} \\[10pt] -ry_* + \dfrac{\theta y_*}{\alpha + y_*} & -ry_* - \dfrac{\theta x_* y_*}{(\alpha + y_*)^2} \end{pmatrix}.
\end{equation*}
We compute the determinant of this matrix:
\begingroup
\allowdisplaybreaks
\begin{align*}
    \det J(x_*, y_*) &= \left(rx_* y_* + \dfrac{\theta x_*^2 y_*}{(\alpha +y_*)^2} \right) -\left(rx_*y_* + \dfrac{\theta \alpha r x_* y_*}{(\alpha + y_*)^2} - \dfrac{\theta x_* y_*}{\alpha + y_*} - \dfrac{\theta^2 \alpha x_* y_*}{(\alpha + y_*)^3} \right) \\[20pt]
    &= \dfrac{\theta x_* y_*}{\alpha + y_*} + \dfrac{\theta x_* y_*}{(\alpha + y_*)^2} \left( x_* - ar + \dfrac{\alpha \theta}{\alpha + y_*} \right) \\[20pt]
    &= \dfrac{\theta x_* y_*}{(\alpha + y_*)^2} \left( \alpha + y_* + x_* - \alpha r + \dfrac{\alpha \theta}{\alpha + y_*} \right) > 0,
\end{align*}
\endgroup
where the last inequality follows since $r < 1$. 

The trace of $J(x_*,y_*)$ is
\begin{equation*}
    \text{tr} \ J(x_*,y_*) = -x_*-ry_*  - \dfrac{\theta x_* y_*}{(\alpha + y_*)^2} < 0.
\end{equation*}
Therefore, the steady state $(x_*,y_*)$ is locally asymptotically stable whenever it exists. It is therefore also globally asymptotically stable.

By Theorem \ref{theorem_pos_bounded}, Proposition \ref{prop_dulac}, and the Poincar\'e-Bendixson Theorem, the local asymptotic stability of the steady state $(1,0)$ implies the global asymptotic stability if $\theta < \alpha \gamma$. Similarly, if $\theta > \alpha \gamma$, then $(x_*, y_*)$ is globally asymptotically stable.

We summarize these results as follows.

\begin{theorem}
Consider system (\ref{e11}) - (\ref{e22}) over the region $\mathbb{R}_+^2$.
\begin{enumerate}
    \item If $\gamma > \max\{r,\theta/\alpha \}$, then the only two non-negative steady states of the system are $(x,y) = (0,0)$ and $(x,y) = (1,0)$. The steady state $(0,0)$ is a saddle and the steady state $(1,0)$ is a stable node. Furthermore, the steady state $(1,0)$ is globally asymptotically stable on $\mathbb{R}_+^2$. \\
    
    \item If $r < \gamma < \theta/\alpha$, then there exists an additional, positive, steady state, $(x_*, y_*)$. The steady states $(0,0)$ and $(1,0)$ are unstable (saddles) and $(x_*, y_*)$ is globally asymptotically stable on $\mathbb{R}_+^2$.
\end{enumerate}
\label{prop_dynamics_2D}
\end{theorem}

\begin{figure}[H]
\centering
	\begin{subfigure}[t]{0.45\textwidth}\centering
\includegraphics[width=\textwidth]{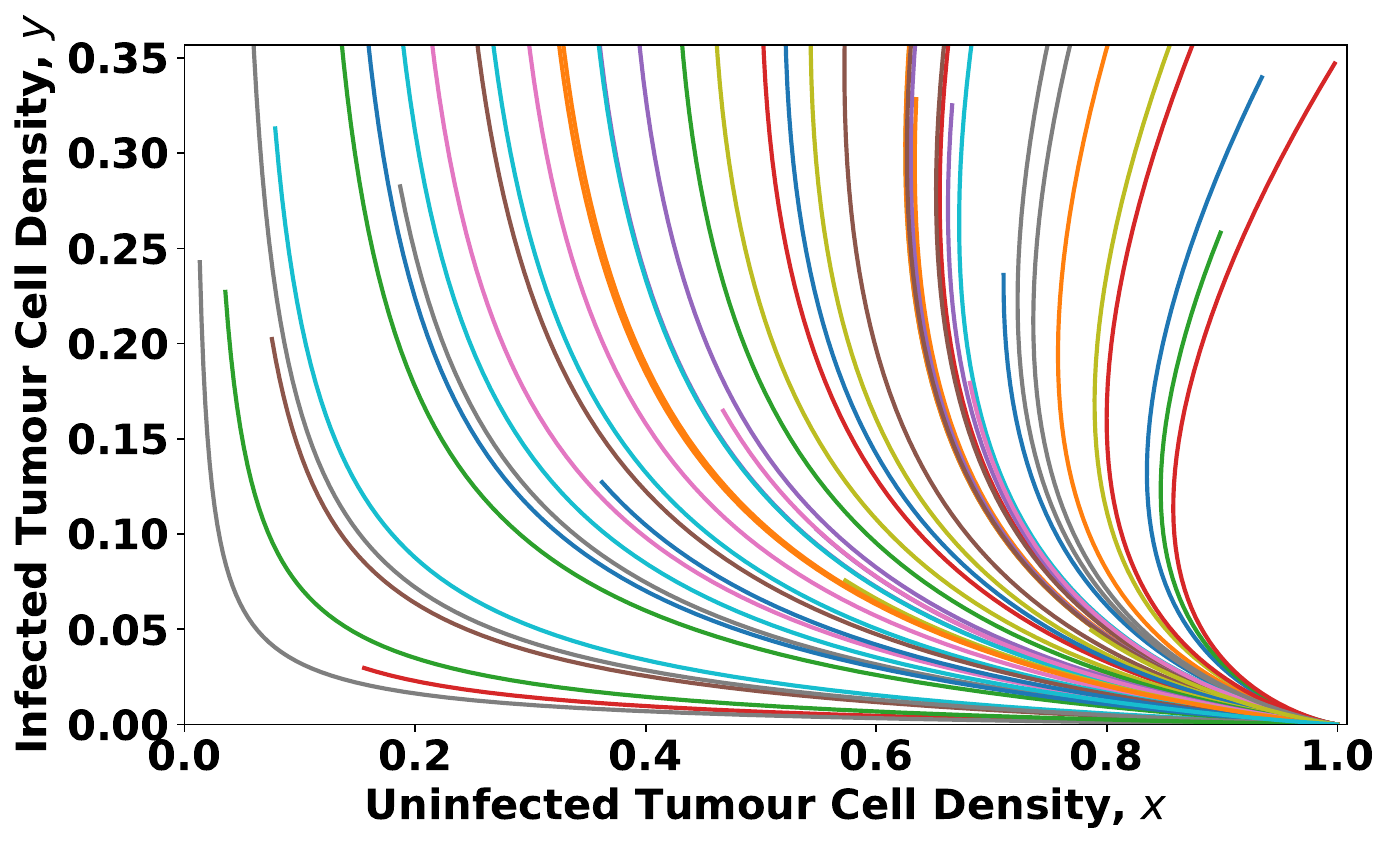}
         \caption{}\label{phase_portrait1a}
    \end{subfigure}
    \hspace{5mm}
\begin{subfigure}[t]{0.45\textwidth}\centering
\includegraphics[width=\textwidth]{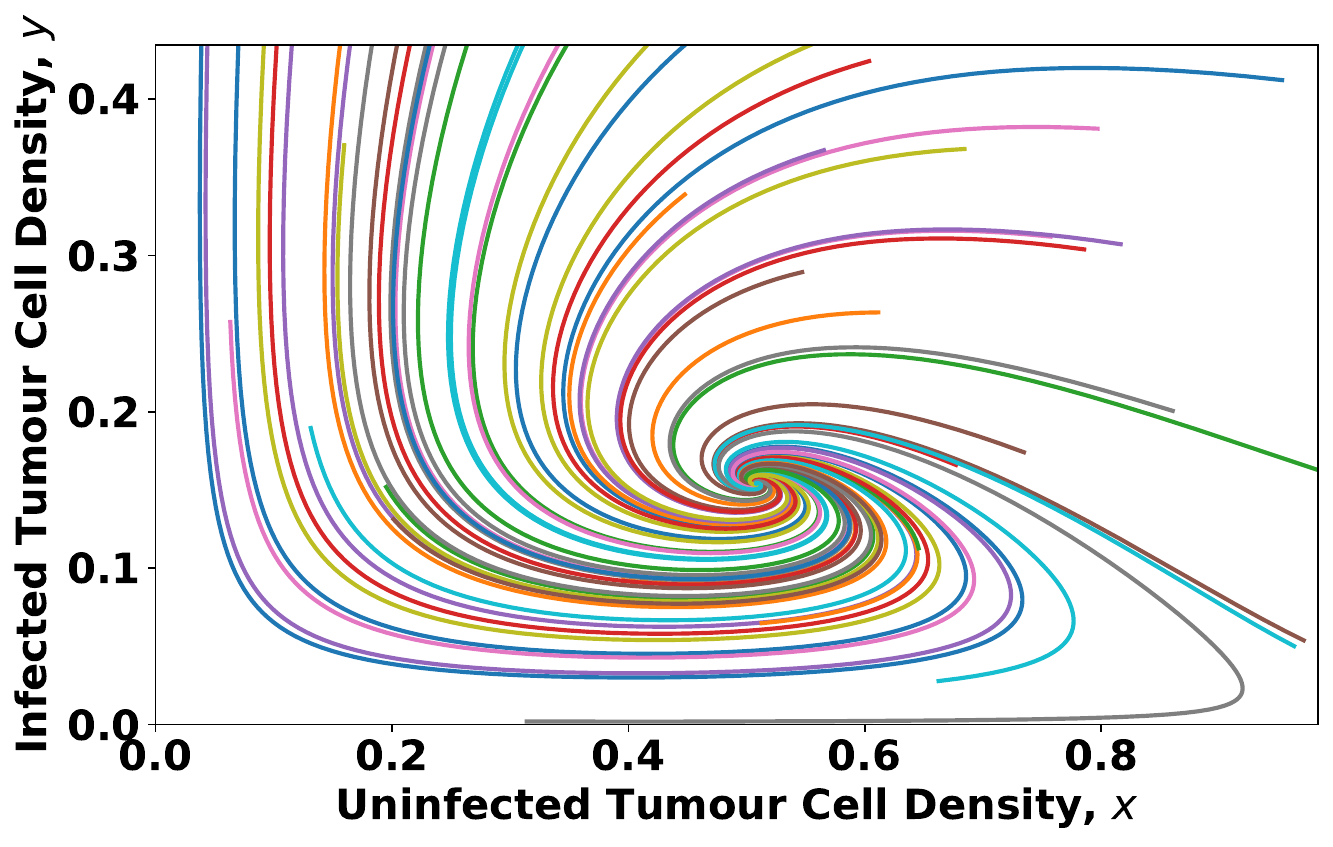}\\
 \caption{}\label{phase_portrait1b}
\end{subfigure}
\caption{The phase portrait of system (\ref{e11}) - (\ref{e22}) when $r < \gamma$. (a) When $\gamma > \max \{ r, \theta/\alpha \}$, the trajectories approach the steady state $(1,0)$ as $t \rightarrow \infty$. (b) When $r < \gamma < \theta/\alpha$, the trajectories approach the steady state $(x_*, y_*)$ as $t \rightarrow \infty$. }
\label{phase_portrait_1}
\end{figure}

We numerically illustrate Theorem \ref{prop_dynamics_2D} in Figure  \ref{phase_portrait_1}. The phase portraits in Figure \ref{phase_portrait_1} are produced with all parameters, except for $\alpha$, being assigned (after non-dimensionalization) based on the values in Table 2. Doing so gives the parameter values $\theta = 2.52908, r = 0.531107, \gamma = 1.29362.$ In Figure \ref{phase_portrait_1}(a), we set $\alpha = 10.0$ and in Figure \ref{phase_portrait_1}(b), we set $\alpha = 1.0$.

Clinically, the stability of $(1,0)$ is not a favourable result, representing a failure of the OV treatment. From Proposition \ref{prop_dynamics_2D}, we see that one condition which leads to this occurrence is the virus-induced death rate, $\gamma$, being made sufficiently large. This leads to an idea which will come up again in the case of oxygen dependence: Having a virus-induced death rate which is too large relative to the infection rate will decrease the efficacy of the OV. Instead, it is important to make sure that the tumour cells are not being killed faster than they are able to infect adjacent tumour cells. This suggests that when engineering an OV, it is important to achieve an appropriate balance between the infection rate and oncolysis rate of the virus. 

Note that the condition $\gamma > \theta/\alpha$ in Proposition (\ref{prop_dynamics_2D}) follows directly from setting $\mathcal{R}_0 < 1$ where $\mathcal{R}_0$ is the basic reproduction number. Following the approach of \cite{reproduction_number}, the basic reproduction number may be computed by using the next generation method. We omit the details here.

So far, we have considered the case where $r < \gamma$, i.e., when the growth rate of the infected tumour cells is bounded by their death rate. We have seen that total extinction of the uninfected tumour cells is not possible in this case. We now consider the case $r > \gamma$. In this case, we have an additional non-negative steady state, $(0,(r-\gamma)/r)$. This steady state may represent a semi-successful treatment outcome in the case $\gamma \approx r$. Hence, stability of this steady state is clinically preferable. 

We note first that if $r > \gamma$, the matrix $J(0,0)$ has two positive eigenvalues and hence, $(0,0)$ is an unstable node. We  assume that $\theta < \alpha \gamma$ in order to rule out the existence of a non-negative co-existence steady state. In this case, the eigenvalues of $J(1,0)$ remain negative and so $(1,0)$ remains a stable node. Linearizing system (\ref{e11}) - (\ref{e22}) at the steady state $(0,(r-\gamma)/r)$:
\begin{equation*}
    J\left(0,\dfrac{r-\gamma}{r} \right) = \begin{pmatrix} \dfrac{r \gamma (\alpha + \theta + 1) - (r^2\theta + \gamma^2)}{r(\alpha r + r - \gamma)} & 0 \\[20pt] \dfrac{(r-\gamma)(\theta + \alpha + \alpha r - r)}{\alpha r + r - \gamma} \quad & \gamma - r \end{pmatrix}.
\end{equation*}
Since this is a lower triangular matrix, the eigenvalues are the elements of the main diagonal. The eigenvalue $\gamma - r$ is negative since $r > \gamma$. The remaining eigenvalue is positive since 
\begin{align*}
    r \gamma (\alpha + \theta + 1) &= r (\alpha \gamma) + r \gamma \theta + r \gamma > (r^2)(\theta) + r \gamma \theta + (\gamma)(\gamma) > r^2 \theta + \gamma^2.
\end{align*}
Therefore, $(0,(r-\gamma)/r)$ is a saddle and hence unstable. Therefore, even in the case where $r > \gamma$, the only locally stable steady state is $(1,0)$. This also remains true if $r = \gamma$, as can be seen via direct substitution.

The steady state $(1,0)$ is unstable if $\theta > \alpha \gamma$ and $(0,0)$ is always unstable. If $r > \gamma$, then the steady state $(0,(r-\gamma)/r)$ is locally asymptotically stable if and only if
\begin{equation}
    \theta > \gamma \left( \dfrac{\alpha}{r-\gamma} + \dfrac{1}{r} \right).  \label{existence_cond32}
\end{equation}
This condition is obtained by requiring all the eigenvalues of $J(0,(r-\gamma)/r)$ to be negative. Note that since $r < 1$, condition (\ref{existence_cond32}) implies that $\theta > \alpha \gamma$. Since all solutions are non-negative and bounded, and closed orbits may not exist, it follows that violating condition (\ref{existence_cond32}) implies the existence and stability of the positive steady state $(x_*, y_*)$.

We summarize the results of the case $r > \gamma$ in the following theorem.

\begin{theorem}
Consider system (\ref{e11}) - (\ref{e22}) over the region $\mathbb{R}_+^2$. If $r > \gamma$, then: 
\begin{enumerate}
    \item The steady state $(x,y) = (1,0)$ is globally asymptotically stable on $\mathbb{R}_+^2$ if $$\theta < \alpha \gamma.$$ 
    \item The positive steady state $(x,y) = (x_*, y_*)$ is globally asymptotically stable on $\mathbb{R}_+^2$ if $$ \alpha \gamma < \theta < \gamma \left( \dfrac{\alpha}{r-\gamma} + \dfrac{1}{r} \right). $$
    \item The steady state $(x,y) = (0, (r-\gamma)/r)$ is globally asymptotically stable on $\mathbb{R}_+^2$ if $$ \theta > \gamma \left(\dfrac{\alpha}{r-\gamma} + \dfrac{1}{r} \right) .$$
\end{enumerate} \label{additional_important_theorem_chap3_1}
\end{theorem}

\begin{figure}[H]
\centering
	\begin{subfigure}[t]{0.45\textwidth}\centering
\includegraphics[width=\textwidth]{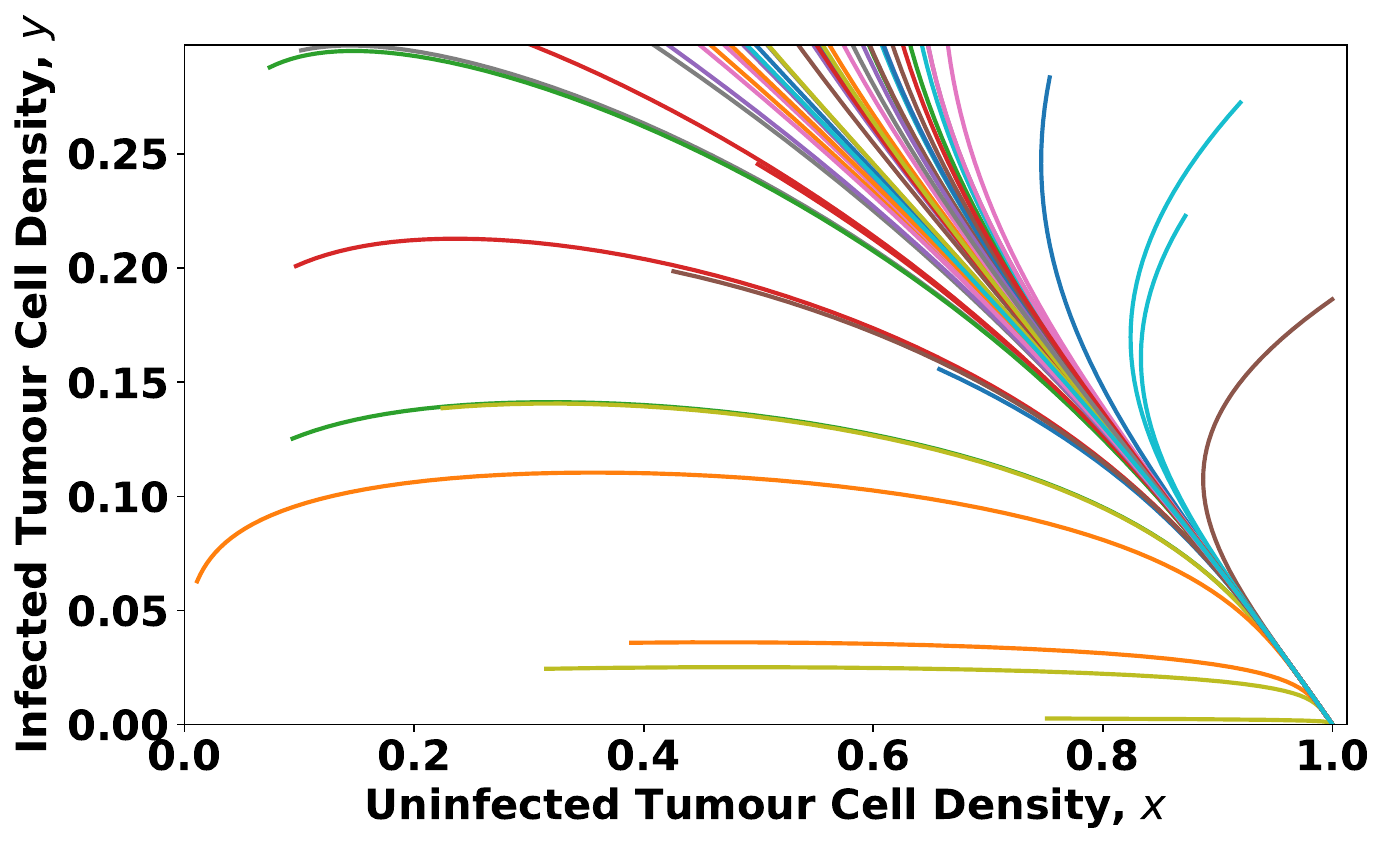}
         \caption{}
    \end{subfigure}
    \hspace{5mm}
\begin{subfigure}[t]{0.45\textwidth}\centering
\includegraphics[width=\textwidth]{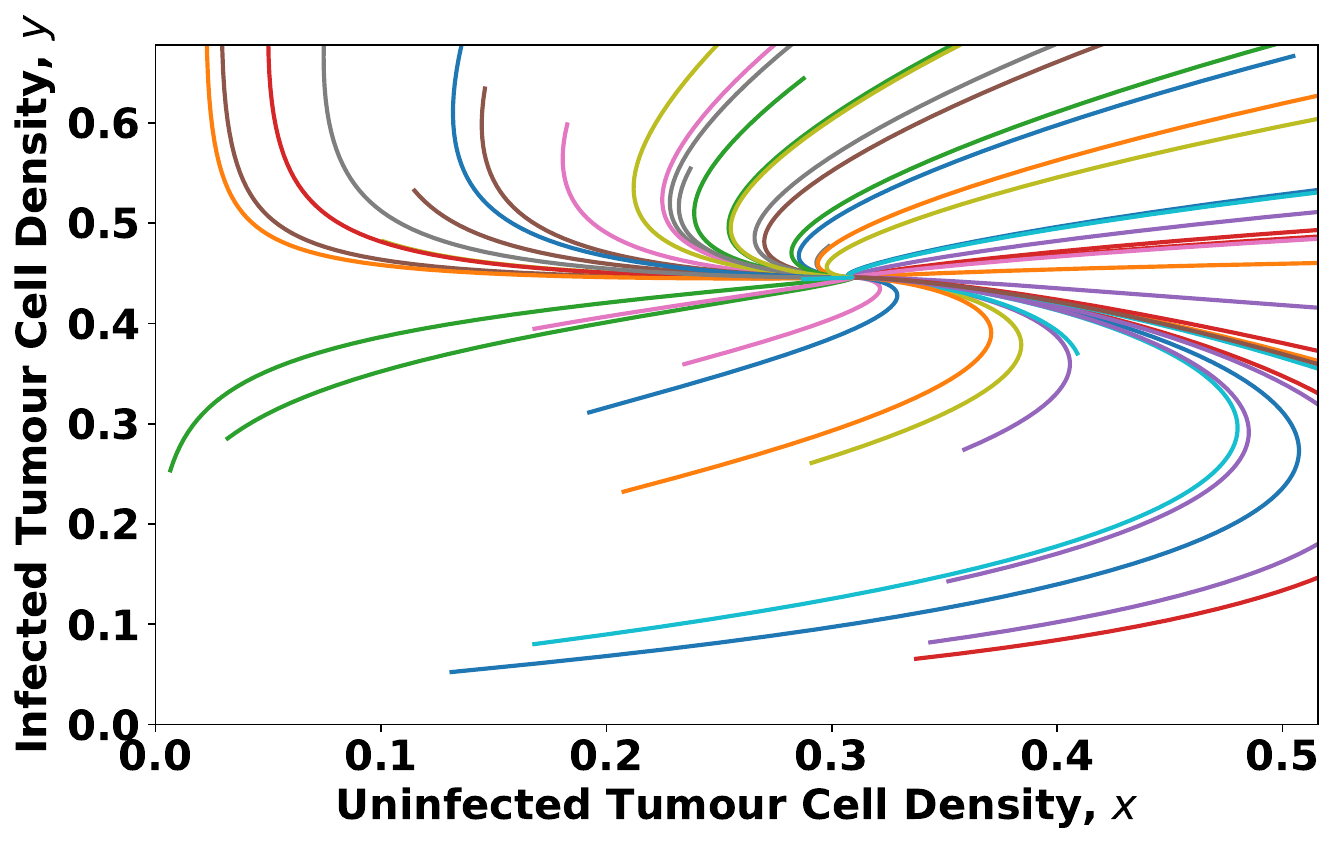}\\
 \caption{}
\end{subfigure}
\begin{subfigure}[t]{0.45\textwidth}\centering
\includegraphics[width=\textwidth]{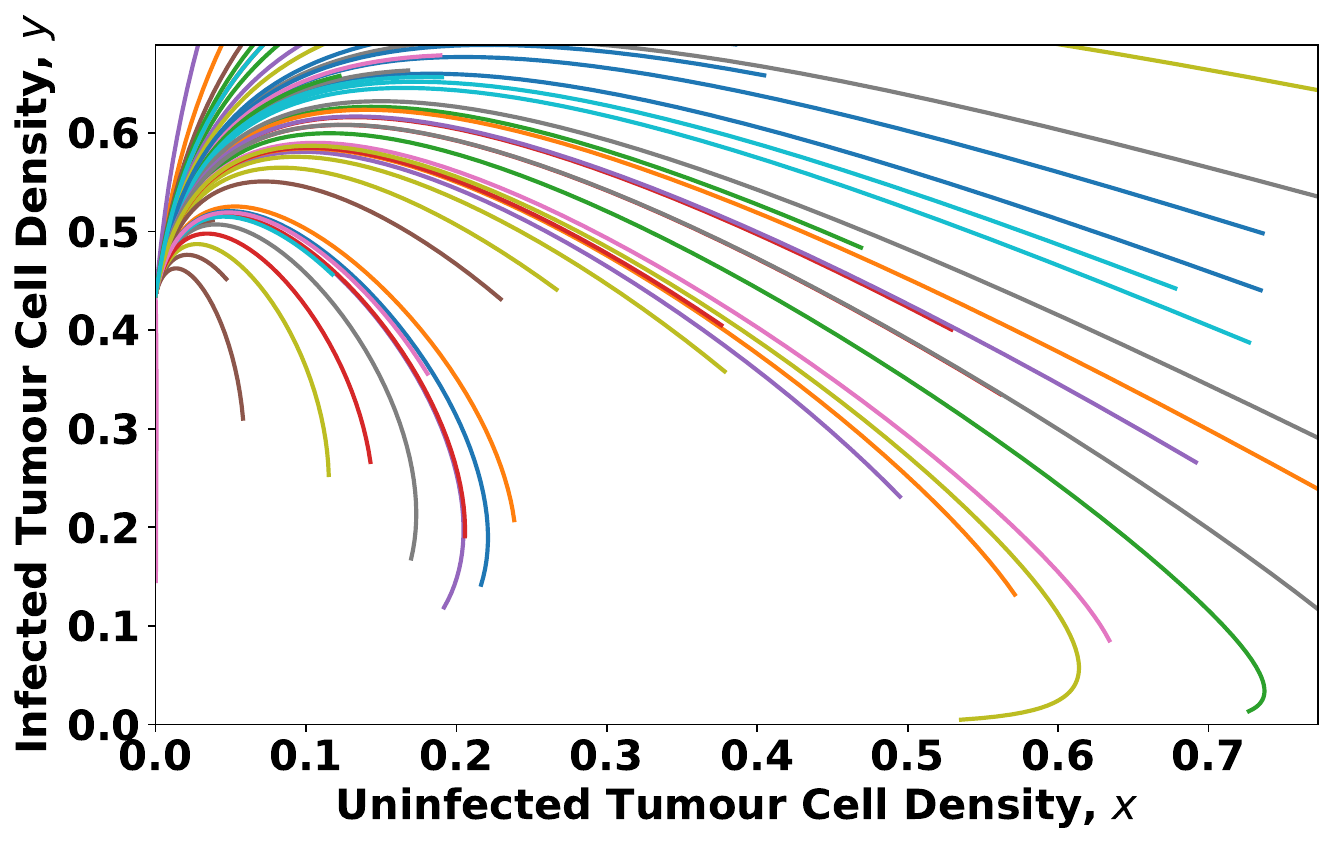}\\
 \caption{}
\end{subfigure}
\begin{subfigure}[t]{0.45\textwidth}\centering
\includegraphics[width=\textwidth]{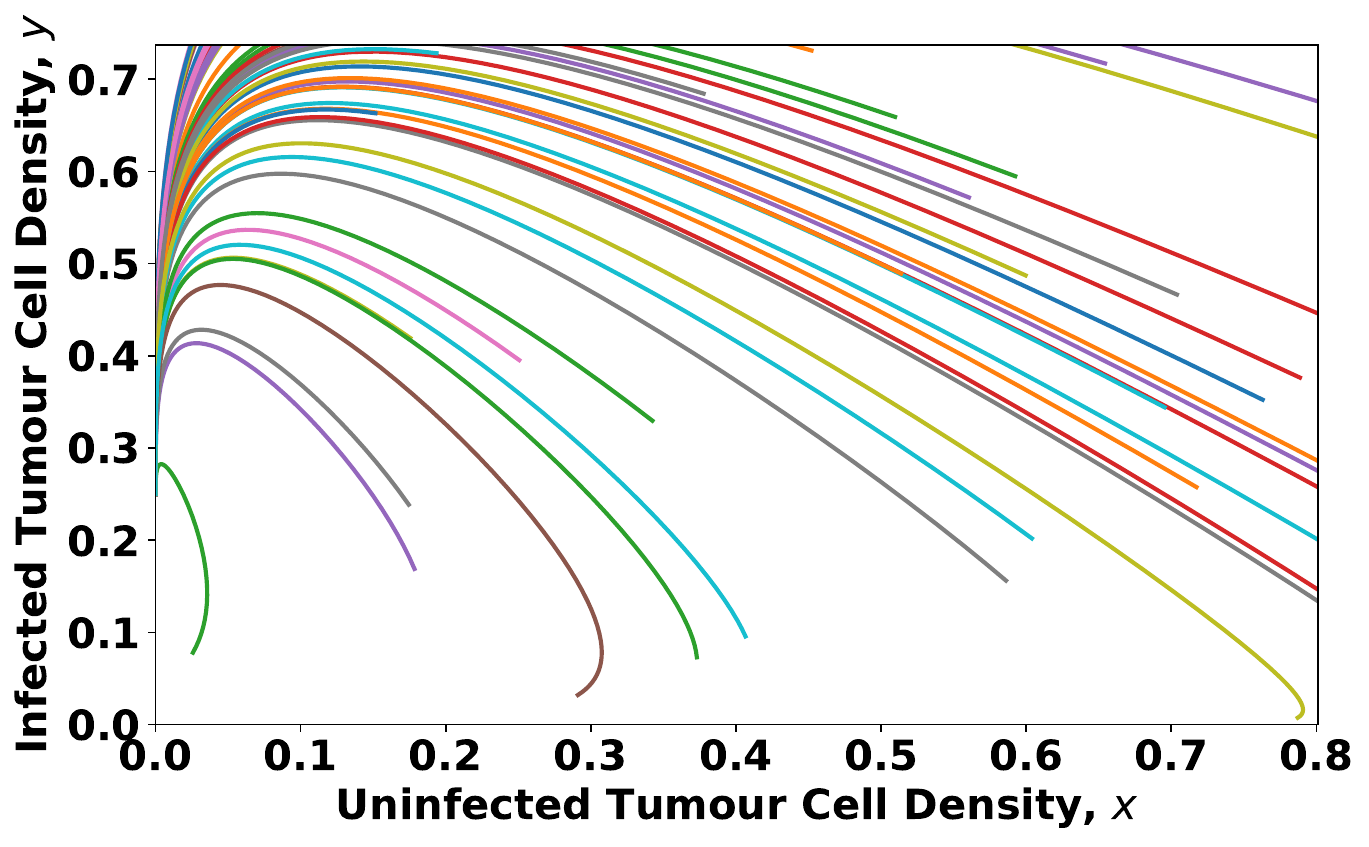}\\
 \caption{}
\end{subfigure}
\caption{The phase portrait of system (\ref{e11}) - (\ref{e22}) when $r > \gamma$. (a) When $\theta < \alpha \gamma$, the trajectories approach the steady state $(1,0)$ as $t \rightarrow \infty$. (b) When $\alpha \gamma < \theta < \gamma (\alpha/(r-\gamma) + 1/r) $, the trajectories approach the steady state $(x_*,y_*)$ as $t \rightarrow \infty$. (c) \& (d) When $ \theta > \gamma (\alpha/(r-\gamma) + 1/r) $, the trajectories approach the steady state $(0,(r-\gamma)/r)$ as $t \rightarrow \infty$. In (d), the value of the parameter $r$ is taken closer to $\gamma$ than in (c), resulting in decreased density of infected tumour cells.  }
\label{profzoussuggestedfigurechap3}
\end{figure}

We numerically illustrate Theorem \ref{additional_important_theorem_chap3_1} in Figure \ref{profzoussuggestedfigurechap3}. The phase portraits in this figure are produced by using the parameter values given in the code in Appendix A, except for the parameters $r, \theta,$ and $\gamma$. We set $\gamma = 0.3$. In Figure \ref{profzoussuggestedfigurechap3} (a), we set $r = 0.5311$ and  $\theta = 0.01$. In Figure \ref{profzoussuggestedfigurechap3} (b), we set $r = 0.5311$ and  $\theta = 0.3$. In Figure \ref{profzoussuggestedfigurechap3} (c), we set $r = 0.5311$ and  $\theta = 0.9$. In Figure \ref{profzoussuggestedfigurechap3} (d), we set $r = 0.4$ and  $\theta = 1.4$. 

Biologically, the case $\theta > \alpha \gamma$ corresponds to a low virus-induced death rate relative to the infection rate (since in practice, $\alpha$ is typically less than 1). This condition leads to a more clinically favourable outcome compared to the condition $\theta < \alpha \gamma$, as the uninfected tumour cell-dominant steady state becomes unstable. If we then consider the additional condition $r > \gamma$, then there exists an infected tumour cell-dominant steady, $(0,(r-\gamma)/\gamma)$, which corresponds to complete eradication of uninfected tumour cells. Biologically, this clinically favourable steady state exists when infected tumour cells can proliferate at a greater rate than they are destroyed by the virus. This (perhaps rather unintuitively) suggests that an OV should not be engineered to hinder the proliferation capability of the cancer cells and, in fact, a greater growth rate of the infected cancer cells can lead to improved clinical outcomes. The idea is to minimize $(r-\gamma)/r$ while also ensuring that the infected tumour cell-dominant steady state is stable, i.e., inequality (\ref{existence_cond32}) holds. The modelling suggests that the most potent OV is one with a high infection rate, low oncolysis rate, and that minimally inhibits the proliferation rate of the cancer cells. By taking $\gamma \rightarrow r^{-}$, we have $y \rightarrow 0$ as $t \rightarrow \infty$ as long as $\theta$ still satisfies condition (\ref{existence_cond32}). While this might lead to the naive assumption of simply engineering a virus which has a very large infection rate compared to the proliferate rate of tumour cells, this type of OV may also be associated with increased toxicity \cite{simpson}, adding another layer of complexity.

\begin{figure}[h]
\begin{center}
\begin{subfigure}[t]{0.45\textwidth}\centering
\begin{tikzpicture}[xscale=10.0,yscale=5.0]
\draw[line width = 0.4mm, ->] (0,0) -- (0.6,0) node[right] {$\boldsymbol{{\gamma}}$};
\draw[line width = 0.4mm, ->] (0,0) -- (0,1.1) node[left] {$\boldsymbol{{\theta}}$};
\draw[ domain = 0:(0.21/0.3954-0.075),line width = 0.4mm, smooth, ->, variable =\N] plot ({\N},{\N + 0.1*\N/( (0.21/0.3954) -\N)}) node[above]{$\boldsymbol{{\theta = \gamma \left(\dfrac{\alpha}{r - \gamma} + \dfrac{1}{r} \right)}}$};
\draw[line width = 0.4mm, -] (0.21/0.3954-0.05,0.02) -- (0.21/0.3954-0.05,-0.02) node[below] {$\boldsymbol{{r}}$};
\draw[line width = 0.4mm, -] (0,0) -- (0,0) node[below] {$\boldsymbol{{0}}$};
\fill[blue, fill opacity = 0.2] (0,1)--(0,0) -- plot[domain=0:(0.21/0.3954-0.08), smooth] (\x,{\x + 0.1*\x/( (0.21/0.3954) -\x)}) -- (0.21/0.3954-0.08,1) -- cycle;
\node[circle,fill,fill=black!45!green,label={above: small $\theta$}] at (0.3,0.1) {};
\node[circle,fill,fill=black!20!red,label={above: large $\theta$}] at (0.3,0.8) {};
\end{tikzpicture}
\caption{}
\end{subfigure}
\hspace{5mm}
\begin{subfigure}[t]{0.5\textwidth}\centering
\includegraphics[width=\textwidth]{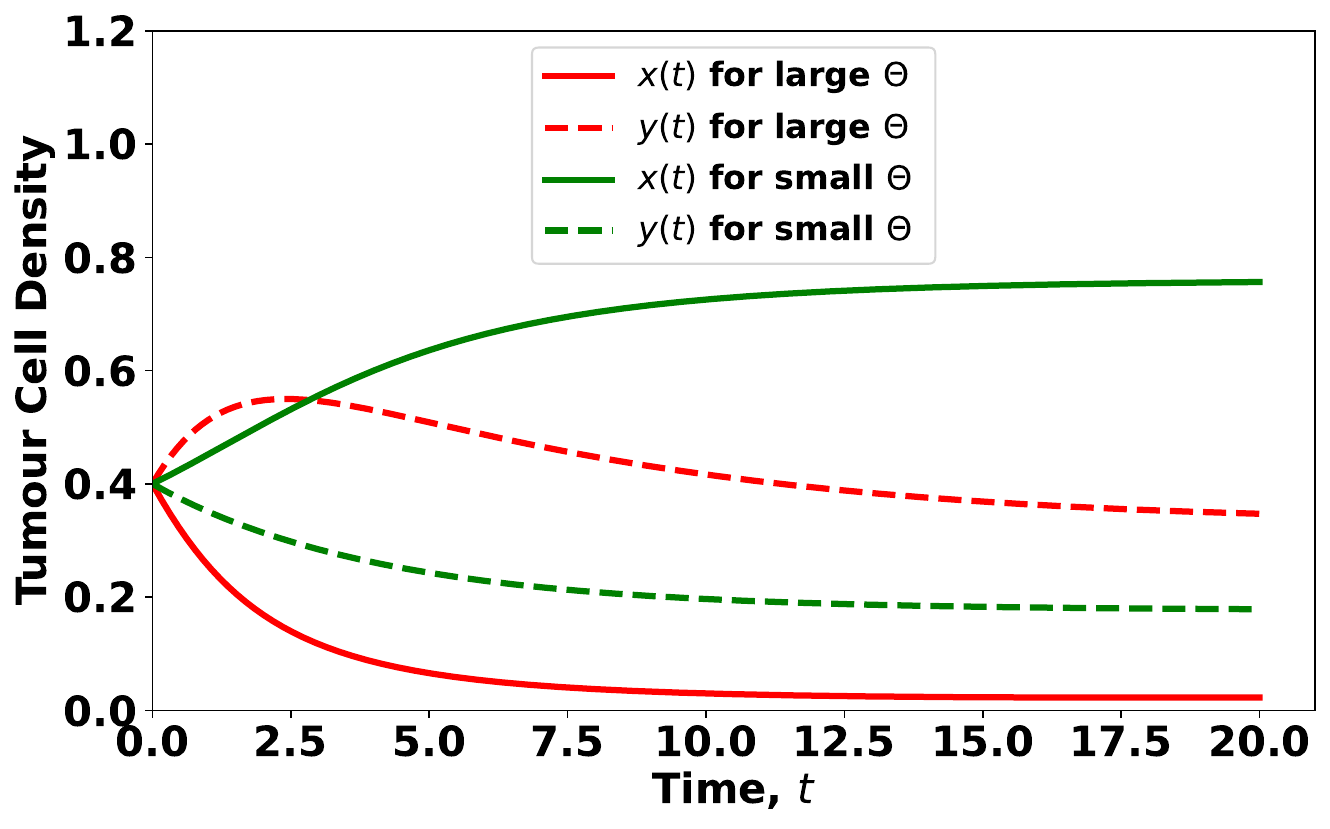}\\
 \caption{}\label{}
\end{subfigure}
\end{center}
\caption{{(a) If the pair $(\gamma_*,\theta_*)$ belongs to the blue region of the $\gamma \theta$-plane, then the infected tumour cell-dominant steady state is stable. (b) The red curves are the solution curves when inequality (\ref{existence_cond32}) is satisfied. The green curves are the solution curves when the condition is not satisfied. In the latter case, the solutions converge to the positive steady state.}}
\label{fig:chap3_extra01}
\end{figure}
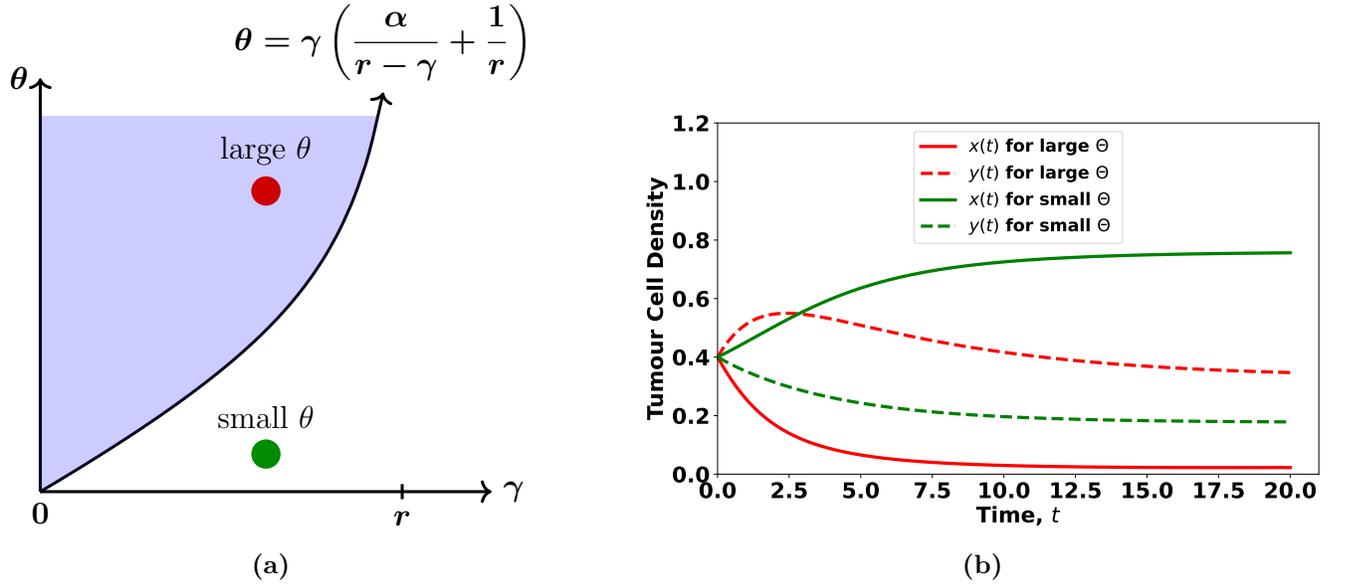

Figure \ref{fig:chap3_extra01} (a) gives guidance on how to choose the infection rate, $\theta$, given the oncolysis rate, $\gamma$. It can be seen in Figure \ref{fig:chap3_extra01} (b) that if the infection rate is too small, the tumour cell densities will converge to the positive steady state. On the other hand, if $\theta$ is large enough, then all of the tumour cells are eventually infected by the virus. 

We summarize the existence and stability results of this section in Table 1.

\begin{table}[H]
    \centering
    \caption{Conditions for Existence and Stability of Steady States of System (\ref{e11}) - (\ref{e22})}
  \resizebox{18cm}{!}{
  \begin{tabular}{c|c|c}
\hline
\hline
\textbf{Steady State} & \textbf{Existence Condition(s)} & \textbf{Global Asymptotic Stability Condition} \\\hline 
 $(0,0)$ &  Always & Unstable \\[5pt]
 $(1,0)$ &  Always & $\theta < \alpha \gamma$ \\[12pt]
 $\left(0,\dfrac{r-\gamma}{r} \right)$ & $r > \gamma$  &  $\theta > \gamma \left(\dfrac{\alpha}{r-\gamma} + \dfrac{1}{r} \right)$ \\[18pt]
 $(x_*, y_*)$ & $\alpha \gamma < \theta < \gamma \left( \dfrac{\alpha}{r-\gamma} + \dfrac{1}{r} \right)$ & Existence \\[10pt] \hline
  
\end{tabular}}
\end{table}

\subsection{Dynamics of the local model -- case II: oxygen dependence}

We now perform a local stability analysis of the relevant  steady states of system (\ref{e1}) - (\ref{e3}). We begin by computing the Jacobian matrix of this system,
\begin{equation}
J(x,y,z) = \begin{pmatrix} 1 - 2x - y - \dfrac{\theta(z)y}{\alpha + y}  & -x - \dfrac{\alpha \theta(z)x}{(\alpha + y)^2} & -\dfrac{\theta'(z)xy}{\alpha + y}  \\[18pt] -ry - \dfrac{\theta(z)y}{\alpha + y} & r - rx - 2ry + \dfrac{\alpha \theta(z) x}{(\alpha + y)^2} - \gamma(z) \ \ & \dfrac{\theta'(z)xy}{\alpha + y} - \gamma'(z)y  \\[18pt] -q_1z & -q_2z & -\beta - q_1x - q_2 y\end{pmatrix}. \label{Jacobian}
\end{equation}
We first consider the simplest steady state, the \textit{tumour-free} steady state, $(x,y,z) = (0,0,1)$. Linearizing the system about this point gives
\begin{equation}
J(0,0,1) = \begin{pmatrix} 1 & 0 & 0 \\ 0 & r - \gamma(1) & 0 \\ -q_1 & -q_2 & -\beta  \end{pmatrix}, \label{tumour-free}
\end{equation}
which is a lower triangular matrix with eigenvalues $1, r - \gamma(1), -\beta$. Since this matrix will always have a positive eigenvalue, the tumour-free steady state is unstable. The maximum dimension of its stable manifold is 2, which occurs if and only if $r < \gamma(1)$. This corresponds to the fact that if the virus-induced death rate of tumour cells, $\gamma$, is sufficiently large, then there will be larger domain of initial conditions for which the solution will converge to the tumour-free steady state.
Next, we consider the case where the uninfected tumour cells dominate, i.e., $x = 1$ and $y = 0$. In this case, we have the following steady state, which corresponds to the failure of OV treatment:
\begin{equation}
(x,y,z) = \left( 1, 0, z^*  \right), \quad \text{where} \quad z^* := \dfrac{\beta}{\beta + q_1}.
\end{equation}
Linearizing the system at this steady state gives
\begin{equation}
J\left( 1, 0, z^* \right) = \begin{pmatrix} -1 & -1 - \dfrac{\gamma(z^*)}{\alpha} & 0 \\[13pt]  0 & \dfrac{\theta(z^*)}{\alpha} - \gamma(z^*)  & 0 \\[13pt] -q_1 z^* & -q_2 z^* & -\beta - q_1 \end{pmatrix}, \label{tumour-dominates}
\end{equation}
which has eigenvalues
\begin{equation}
\lambda_{1}^u = -1, \quad \lambda_2^u = \dfrac{\theta(z^*)}{\alpha} - \gamma(z^*), \quad \lambda_3^u = -\beta - q_1.
\end{equation}
Considering the conditions for which these eigenvalues are all negative gives the following proposition.

\begin{proposition}
The tumour-dominant steady-state, $(x,y,z) = (1,0,z^*)$, is locally asymptotically stable if $\theta(z^*) < \alpha \gamma(z^*)$. \label{prop2}
\end{proposition}

From a clinical perspective, the local asymptotic stability of the tumour-dominant steady-state is an unfavourable result. Biologically, this occurs when the infection rate of tumour cells by the OV is too low compared to the virus-induced death rate. This leads to an important insight: engineering a virus which can destroy tumour cells at a fast rate is not useful if the infection rate is too low. It is important to have a virus which is sufficiently effective at infecting cancer cells - not just destroying them. The inequality in Proposition \ref{prop2} can give an estimate on how large these rates should be for a useful OV. 

We are also interested in the existence of an uninfected tumour cell-free steady, i.e., one of the form $(0,y_*, z_*)$. From system (\ref{e1}) - (\ref{e3}), it can been seen that such a solution may be determined by solving the system
\begin{align}
    r(1-y) - \gamma(z) &= 0, \label{alg_sys7} \\
    \beta(1-z) - q_2 yz &= 0. \label{alg_sys8}
\end{align}
This system may have no solutions, one solution, or multiple solutions depending on the properties of the oncolysis function, $\gamma(z)$. Stability of this steady state is favourable and hence, we impose the additional condition $\gamma_\infty < r$ so as to ensure the existence of a  positive solution of system (\ref{alg_sys7}) - (\ref{alg_sys8}). Notice that this condition is similar to the existence of the uninfected tumour cell-free steady state condition in Subsection 3.1. Moreover, it should also be noted that $0 \leq y_*, z_* \leq 1$.

Linearizing at $(0,y_*,z_*)$ gives the matrix
\begin{equation}
    J(0,y_*,z_*) = \begin{pmatrix}  1 - y_* - \dfrac{\theta(z_*)y_*}{\alpha + y_*}  & 0 & 0  \\[18pt] -ry_* - \dfrac{\theta(z_*)y_*}{\alpha + y_*} & -ry_* \ \ &  - \gamma'(z_*)y_*  \\[18pt] -q_1z_* & -q_2z_* & -\dfrac{\beta}{z_*} \end{pmatrix}.
\end{equation}
The eigenvalues of this matrix are
\begin{equation*}
    \lambda_1^n = 1 - y_* - \dfrac{\theta(z_*)y_*}{\alpha + y_*}, \quad \lambda_{2,3}^n = \dfrac{-(\beta +  y_* z_*) \pm \sqrt{(\beta + ry_*z_*)^2 - 4z_* (\beta r y_* - \gamma'(z_*)q_2 y_* z_*^2)}}{2z_*}.
\end{equation*}
It is clear that all of these eigenvalues have no imaginary part. Hence, $(0,y_*,z_*)$ is either a stable node or a three-dimensional saddle. The former case is preferable, as all tumour cells will eventually be infected as $t \rightarrow \infty$. This occurs when the eigenvalues are all negative, leading to the following proposition.

\begin{proposition}
Consider the steady state $(0,y_*,z_*)$, where $y_*$ and $z_*$ satisfy the equations
\begin{equation}
    \gamma(z_*) = r \left( 1 + \dfrac{\beta}{q_2} -\dfrac{\beta}{q_2 z_*} \right), \quad y_* = \dfrac{\beta}{q_2} \cdot \dfrac{1-z_*}{z_*}. \label{diagram_equations_chap3_1}
\end{equation}
If $\gamma_\infty < r$, then such $y_*$ and  $z_*$ exist and $0 \leq y_* \leq 1$, $\beta/(\beta + q_2) \leq z_* \leq 1$. Moreover, the steady state $(0,y_*,z_*)$ is locally asymptotically stable if and only if
\begin{align}
    \theta (z_*) > \dfrac{(1-y_*)(\alpha + y_*)}{y_*} \quad \text{and} \quad \gamma'(z_*) < \dfrac{\beta r}{q_2 z_*^2}. 
\end{align}\label{prop_gamma_prime}
\end{proposition}
\begin{remark}
The condition $\gamma_\infty < r$ is a sufficient condition for the existence of the steady state $(0,y_*,z_*)$. A necessary and sufficient condition for the existence of this steady state is $\gamma_\infty < r(1+\beta/q_2)$. The latter condition, however, does not guarantee that $z_* \leq 1$.
\end{remark}

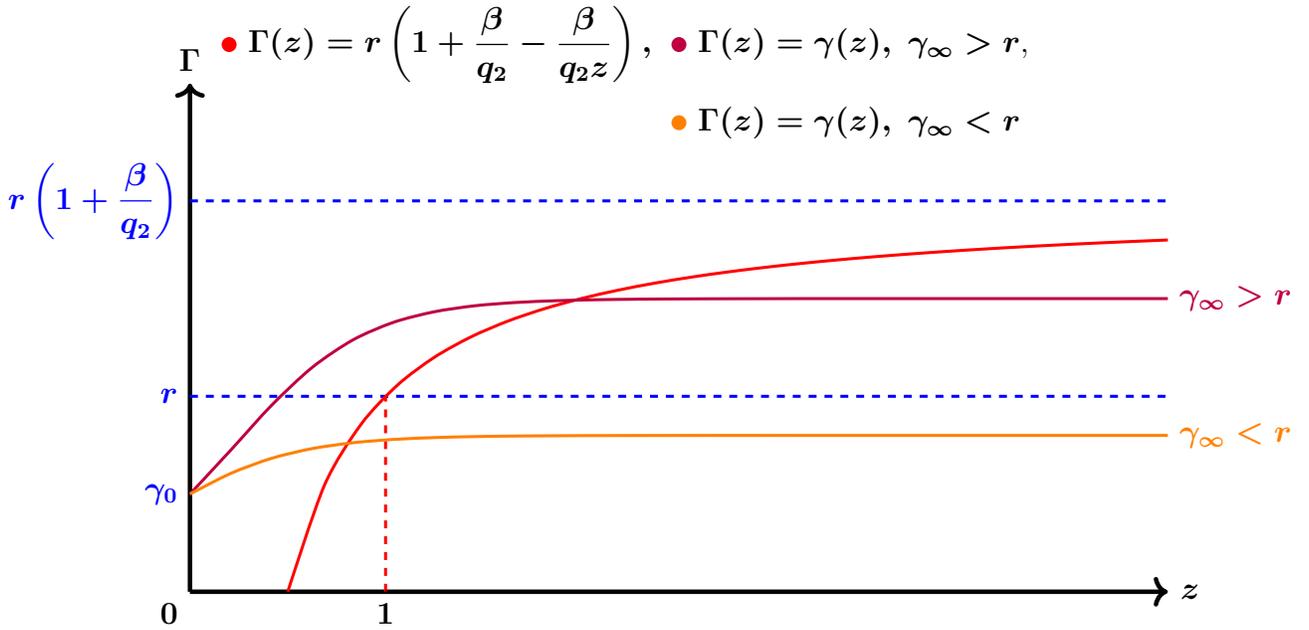
\begin{figure}[H]
\begin{center}
\begin{tikzpicture}[xscale=2.6,yscale=5.2]
\draw[line width = 0.6mm, ->] (0,0) -- (5,0) node[right] {$\boldsymbol{z}$};
\draw[line width = 0.6mm, ->] (0,0) -- (0,1.3) node[above] {$\boldsymbol{\Gamma}$};
\draw[blue,line width = 0.4mm, dashed, -] (0,1) -- (5,1) node[left] {};
\draw[blue,line width = 0.4mm, dashed, -] (0,0.5) -- (5,0.5) node[left] {};
\draw[blue,line width = 0.4mm, dashed, -] (0,1) -- (0,1) node[left] {$\boldsymbol{r \left( 1 + \dfrac{\beta}{q_2} \right)}$};
\draw[blue,line width = 0.4mm, dashed, -] (0,0.5) -- (0,0.5) node[left] {$\boldsymbol{r}$};
\draw[blue,line width = 0.4mm, dashed, -] (0,0.25) -- (0,0.25) node[left] {$\boldsymbol{\gamma_0}$};
\draw[red,line width = 0.4mm, dashed, -] (1,0) -- (1,0.5) node[left] {$\boldsymbol{}$};
\draw[line width = 0.4mm, dashed, -] (1,0) -- (1,0) node[below] {$\boldsymbol{1}$};
\draw[line width = 0.4mm, dashed, -] (0,0) -- (0,0) node[below left] {$\boldsymbol{0}$};
\draw[red, smooth, domain = 0.5:5,line width = 0.4mm, variable =\N] plot ({\N},{0.5*(2 - 1/\N)});
\draw[purple, smooth, domain = 0:5,line width = 0.4mm, variable =\N] plot ({\N},{0.25*0.75/(0.25 + 0.5*exp(-3*\N))})node[right]{$\boldsymbol{{\gamma_\infty > r}}$};
\draw[orange, smooth, domain = 0:5,line width = 0.4mm, variable =\N] plot ({\N},{0.25*0.4/(0.25 + 0.15*exp(-3*\N))}) node[right]{$\boldsymbol{{\gamma_\infty < r}}$};
\node[circle,inner sep = 2pt,fill,fill=red,label={right: $\boldsymbol{\Gamma(z) = r \left( 1 + \dfrac{\beta}{q_2} - \dfrac{\beta}{q_2 z} \right),}$}] at (0.2,1.4) {};
\node[circle,inner sep = 2pt,fill,fill=purple,label={right: $\boldsymbol{\Gamma(z) = \gamma(z), \ \gamma_\infty > r},$}] at (2.5,1.4) {};
\node[circle,inner sep = 2pt,fill,fill=orange,label={right: $\boldsymbol{\Gamma(z) = \gamma(z), \ \gamma_\infty < r}$}] at (2.5,1.2) {};
\end{tikzpicture}
\end{center}
\caption{{The $z$-coordinate of the intersection of the red curve with the purple curve gives the oxygen concentration at the steady state, $z_*$, if $ r < \gamma_\infty < r(1+\beta/q_2)$. The intersection of the red curve with the orange curve gives this steady state if $\gamma_\infty < r$. This latter case guarantees $z_* < 1$. Note that the equation of the red curves comes from (\ref{diagram_equations_chap3_1})  in Proposition \ref{prop_gamma_prime}.}}
\label{fig:chap3_gamma_graph1}
\end{figure}

Proposition \ref{prop_gamma_prime} gives some important conditions for constructing an effective OV. The condition $\gamma_\infty < r$, similarly to Subsection 3.1, gives a sufficient existence condition. The first stability condition is consistent with our previous results: namely, a sufficiently large infection rate is an important factor of OV efficacy. The second stability condition is perhaps more interesting: An oncolysis rate which grows \textit{slowly} in response to increases in oxygen concentration of the tumour microenvironment.

We now consider system (\ref{e1}) - (\ref{e3}) under certain parameter conditions and establish a global stability result concerning the tumour-dominant steady state, $(1,0,z^*)$. In particular, we consider the case $q_1 = 0$ for the sake of mathematical tractability. Biologically, this corresponds to tumour cells which are unable to consume oxygen. While this condition does not typically represent a biologically realistic situation, it may be considered a \textit{best-case scenario}, as less oxygen is consumed and therefore, more oxygen is available to increase the efficacy of the OV. 

We begin by proving an auxiliary result for which we do not need the assumption $q_1 = 0$. Consider the following region in the positive octant in $\mathbb{R}^3$:
\begin{equation}
    \Omega := \left\{ (x,y,z) \in \mathbb{R}^3 : x \geq 0, \ y \geq 0, \ x + y \leq 1, \ 0 \leq z \leq 1 \right\} 
\end{equation}
The idea is to show that this region defines a so-called \textit{trapping region} from which no solution trajectories of system (\ref{e1}) - (\ref{e3}) may exit. We state this in the following lemma.

\begin{lemma}
The region $\Omega \subset \mathbb{R}_+^3$ is a positively invariant set for system (\ref{e1}) - (\ref{e3}). \label{lemma331}
\end{lemma}

\begin{proof}

Let $(x(t),y(t),z(t))$ denote a solution of system (\ref{e1}) - (\ref{e3}) with initial condition in $\mathcal{U}$. Proving this lemma is equivalent to showing that $\Omega$ defines a trapping region for all $t \geq 0$. First note that by Theorem \ref{theorem_pos_bounded}, $x(t), y(t), z(t) \geq 0$ for all $t \geq 0$.  If the trajectory were to exit the region, then by continuity, it would cross either the $z = 1$ boundary or the plane $x+y = 1$ at some time $t^*$. Assume that the trajectory crosses $z=1$ at time $t^*$. Then from equation (\ref{e3}), $z'(t^*) = -q_1 x  - q_2 y \leq 0$.  Therefore, the vector field at this boundary point does not point in the positive $z$ direction, contradicting the assumption since the trajectory may not exit through the $z = 1$ plane. Hence,  we have shown that $z \leq 1$.

\begin{figure}[H]
\begin{center}
\begin{tikzpicture}[xscale=3.0,yscale=3.0]
\draw[line width = 0.6mm, ->] (-0.8,-0.8) -- (-1,-1) node[left] {$\boldsymbol{x}$};
\draw[blue,line width = 0.4mm, dashed, -] (0,0) -- (-0.8,-0.8) node[left] {};
\draw[line width = 0.6mm, ->] (1.12,0) -- (2,0) node[right] {$\boldsymbol{y}$};
\draw[blue,line width = 0.4mm, dashed, -] (0,0) -- (1.12,0) node[left] {};
\draw[line width = 0.6mm, ->] (0,1.12) -- (0,1.8) node[above] {$\boldsymbol{z}$};
\draw[blue,line width = 0.4mm, dashed, -] (0,0) -- (0,1.12) node[left] {};
\draw[blue, line width = 0.6mm, -] (-0.8,-0.8) -- (1.12,0) node[above] {};
\draw[blue, line width = 0.6mm, -] (-0.8,-0.8) -- (-0.8,1.12-0.8) node[above] {};
\draw[blue, line width = 0.6mm, -] (1.12,0) -- (1.12,1.12) node[above] {};
\draw[blue, line width = 0.6mm, -] (0,1.12) -- (1.12,1.12) node[above] {};
\draw[blue, line width = 0.6mm, -] (0,1.12) -- (-0.8,1.12-0.8) node[above] {};
\draw[blue, line width = 0.6mm, -] (1.12,1.12) -- (-0.8,1.12-0.8) node[above] {};
\draw[line width = 0.6mm, -] (-0.8,-0.8) -- (-0.8,-0.8) node[below] {$\boldsymbol{1}$};
\draw[line width = 0.6mm, -] (1.12,0) -- (1.12,0) node[below] {$\boldsymbol{1}$};
\draw[line width = 0.6mm, -] (0,1.12) -- (0,1.12) node[above right] {$\boldsymbol{1}$};
\draw[line width = 0.6mm, -] (0,0) -- (0,0) node[below] {$\boldsymbol{O}$};
\draw[line width = 0.6mm, -] (1.1,1.1) -- (1.1,1.1) node[above right] {$\boldsymbol{  \color{blue}{\Omega} }$};
\node[circle,fill,fill=black!20!red,label={above left: $\boldsymbol{(1,0,1)}$}] at (-0.8,1.12-0.8) {};
\end{tikzpicture}
\end{center}
\caption{{The trapping region $\Omega$. When $q_1 = 0$, all solutions of system (\ref{e1}) - (\ref{e3}) with initial conditions in this region converge to the steady state $(x,y,z) = (1,0,1)$ as $t \rightarrow \infty$.}}
\label{fig:chap3_invariantsubset}
\end{figure}
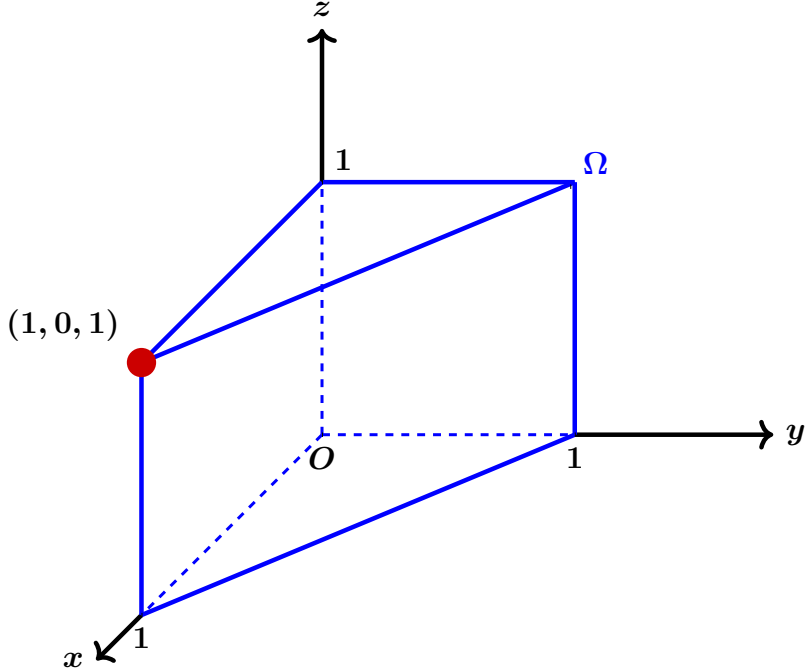

Now we need only establish that no trajectory may exit the region through the plane $x + y = 1$. We do this by showing that the vector field on this plane points into the region $\Omega$. On the plane $x+y = 1$, the sum of equations (\ref{e1}) and (\ref{e2}) is
\begingroup
\allowdisplaybreaks
\begin{align*}
    \dfrac{\mathrm{d}}{\mathrm{d}t} [x(t) + y(t)] &= (x + ry)(1-x-y) - \gamma (z) y \\
    &= (x+ry)(0) - \gamma(z)y < 0.
\end{align*}
\endgroup
Hence, $y'(t) < -x'(t)$ which implies that, by chain rule, $\mathrm{d}y/\mathrm{d}x < -1$. Therefore, the vector field on the plane $x+y = 1$ points into the region and no trajectory may exit through this plane.
We conclude that no trajectory contained in the region $\Omega$ may exit this region, completing the proof. 
\end{proof}

We are now in a position to give the theorem on global stability of the steady state $(1,0,z^*)$. Since we consider the case $q_1 = 0$, we have $z^* = \beta/(\beta+q_1) = 1$. Hence, the steady state becomes $(1,0,1)$. We then have the following theorem:

\begin{theorem}
Consider system (\ref{e1}) - (\ref{e3}) when $q_1 = 0$ and $q_2 > 0$. If   $\theta(z) <  \alpha \gamma (z)$ for $z \in [0,1]$, then $(x,y,z) = (1,0,1)$ is globally asymptotically stable on $\Omega$. \label{lyapunov_theorem1}
\end{theorem}

\begin{proof}
Since $[0,1]$ is a compact interval, we can choose   $\varepsilon \in (0,1)$ such that $\theta(z) < \varepsilon \alpha \gamma(z)$ for all $z \in [0,1]$. Define the function $V : \text{Int} \ {\Omega} \cup \{ (1,0,1) \}  \rightarrow \mathbb{R}$ as
\begin{equation}
    V(x,y,z) := x -\ln(x) - 1 + y + \dfrac{(1-\varepsilon)\gamma_0}{2q_2} \left[z - \ln\left( {z} \right) - 1\right].
\end{equation}
We can show that $V$ defines a Lyapunov function. $V$ is positive definite as $V(1,0,1) = 0$ and $V(x,y,z) > 0$ for all $(x,y,z) \in \text{Int} \ \Omega$. Taking the time derivative of $V$ gives 
\begingroup
\allowdisplaybreaks
\begin{align*}
    \dfrac{\mathrm{d}V}{\mathrm{d}t} &= \left(\dfrac{x-1}{x}\right) \left[ x(1-x-y) - \dfrac{\theta(z)xy}{\alpha + y} \right] + ry(1-x-y) + \dfrac{\theta(z)xy}{\alpha + y} - \gamma(z)y + \dots \\[20pt]
    & \hspace{90mm} \dots + \dfrac{(1-\varepsilon)\gamma_0}{2q_2} \cdot \left(\dfrac{z-1}{z} \right)\left[ \beta (1-z)  - q_2 yz \right] \\[20pt]
    &= -(1-x-ry)(1-x-y) - \dfrac{\theta(z)y(x-1)}{\alpha + y}  + \dfrac{\theta(z)xy}{\alpha + y} - \gamma(z)y + \dfrac{(1-\varepsilon)\gamma_0}{2q_2} \left[ -\dfrac{\beta}{z} (1-z)^2 - q_2 y (z-1) \right] \\[20pt]
    &= -(1-x-ry)(1-x-y) + \dfrac{\theta(z)y}{\alpha + y} - \varepsilon \gamma(z)y - (1-\varepsilon) \gamma(z)y  + \dfrac{(1-\varepsilon)\gamma_0}{2q_2} \left[ -\dfrac{\beta}{z} (1-z)^2 - q_2 y (z-1) \right] \\[20pt]
    &\leq -(1-x-ry)(1-x-y) + \dfrac{[\theta(z) - \varepsilon \alpha \gamma(z) ]y}{\alpha+y} - (1-\varepsilon) \gamma(z)y + \dfrac{(1-\varepsilon)\gamma_0 y}{2}.
\end{align*}
\endgroup
Note that $\dot{V}(1,0,1) = 0.$ Next, since $r < 1$ and $x+y \leq 1$ by Lemma \ref{lemma331},  we have $1-x-ry \geq 0$. Hence, $-(1-x-ry)(1-x-y) \leq 0$. By our assumption, it follows that $\theta (z) - \varepsilon \alpha \gamma(z) < 0$ since $z$ remains in $[0,1]$ by Lemma \ref{lemma331}. Finally, since $\gamma(z) \geq \gamma_0$, it holds that
\begin{equation*}
    - (1-\varepsilon) \gamma(z)y + \dfrac{(1-\varepsilon)\gamma_0 y}{2} < 0.
\end{equation*}
Therefore, $\dot{V} < 0$ on $\text{Int} \ \Omega$. By LaSalle's invariance principle, we conclude that the tumour-dominant steady state $(x,y,z) = (1,0,1)$ is globally asymptotically stable.  
\end{proof}

\begin{remark}
Theorem \ref{lyapunov_theorem1} assumes that $q_2$ is a positive constant. If $q_2 = 0$, establishing global stability is trivial as equations (\ref{e1}) and (\ref{e2}) are decoupled from equation (\ref{e3}) and global asymptotic stability of the tumour-dominant steady state of the resulting two-variable system follows from Theorem \ref{prop_dynamics_2D}.
\end{remark}

As previously stated, the condition $q_1 = 0$ in Theorem \ref{lyapunov_theorem1} biologically represents a best-case scenario in which the uninfected tumour cells are unable to consume oxygen, leading to a more effective adenovirus due to increased oxygen concentration in the tumour microenvironment. In practice $\alpha < 1$ and so the condition $\theta (z) < \alpha \gamma(z)$ reflects a virus which has a significantly larger oncolysis rate compared to its infection rate. This is analogous to the condition required in Theorem \ref{prop_dynamics_2D}, providing further evidence that a very high oncolysis rate is not a favourable characteristic of an oncolytic adenovirus. 

While we do not analytically consider the case $q_1 > 0$, the numerical simulations in Section 4 lead us to conjecture that the steady state $(1,0,z^*)$ remains globally asymptotically stable in this case, under the condition $\theta(z) < \alpha \gamma(z)$.
 
It is clear that the relationship between the functions $\theta(z)$ and $\gamma(z)$ is an important factor in the dynamics of the system. Biologically, if the infection rate is too low relative to the virus-induced death rate, infected tumour cells may die faster than they are able to infect a sufficient number of uninfected cells, hence leading to an uninfected tumour cell-dominant steady state. On the other hand, if the virus-induced death rate is too low, not enough tumour cells will die for the OV to be an effective therapeutic agent. The interplay between these functions and their effect on the OV efficacy is one of the topics of the next section.

\section{Numerical simulations: local model}

In this section, we perform numerical simulations of the local model. We perform the simulations using system (\ref{3eq1}) - (\ref{3eq3}). The units of $u$ and $n$ are cells/mm$^3$ and the units of $c$ are millimolars (mM). Unless otherwise stated, we set the initial conditions to be $u_0= 10000 $ and $n_0 = 100$, as in \cite{parameter_esimates_init_cond}. Similarly to  \cite{oxygen_background_concentration}, we set $c_0 = 4.3751$.

\begin{table}[H]
    \centering
    \caption{Local Model Parameters}
\label{tab1}
  \resizebox{18cm}{!}{
  \begin{tabular}{l|l|l|c}
\hline
\hline
 \textbf{Parameter} &  \textbf{Parameter Name} & \textbf{Value} & \textbf{Reference} \\\hline 
$ r_1 $      &  Growth Rate of Uninfected Tumour Cells & $0.3954$ day$^{-1} $ & \cite{Melanoma_growth_rate} \\\hline
$r_2 $      & Growth Rate of Infected Tumour Cells &  $0.21$ day$^{-1}$  & Estimated ($r_2 < r_1$) \\\hline

$ K $ &  Tumour Carrying Capacity &  $1.0 \times 10^6$ cells/mm$^3$ & \cite{parameter_esimates_init_cond} \\\hline

$ \alpha $ & Hill Constant   & $1.0 \times 10^5$ cells/mm$^3$  & \cite{parameter_esimates_init_cond} \\\hline

$ \phi $ & Oxygenation Rate   & $1.0 \times 10^4$ mM day$^{-1} $ & \cite{oxygen_background_concentration}  \\\hline
$ \beta $ & Oxygen Consumption Rate of Healthy Surrounding Cells     & $5.0976$ day$^{-1} $ & \cite{oxygen_consumptionby_healthy_cells} \\\hline
$ q_1 $ &  Oxygen Consumption Rate of Uninfected Tumour Cells & $5.47 \times 10^{-5}$ mm$^3$ cells$^{-1}$ day$^{-1} $ & \cite{oxygen_consumption_rate} \\\hline

$ q_2 $      & Oxygen Consumption Rate of Infected Tumour Cell   & $2.735 \times 10^{-5}$ mm$^3$ cells$^{-1}$ day$^{-1} $ & \cite{oxygen_consumption_rate}  \\\hline

$ \gamma $      &  Virus-Induced Death Rate  & 0.5115 day$^{-1} $& \cite{parameter_esimates_init_cond} \\\hline

$ \theta$ & Infection Rate     & $1.0$ day$^{-1} $ & \cite{parameter_esimates_init_cond}  \\\hline

\end{tabular}
}
\end{table}

\begin{figure}[H]
    \centering
    \includegraphics[scale=0.4]{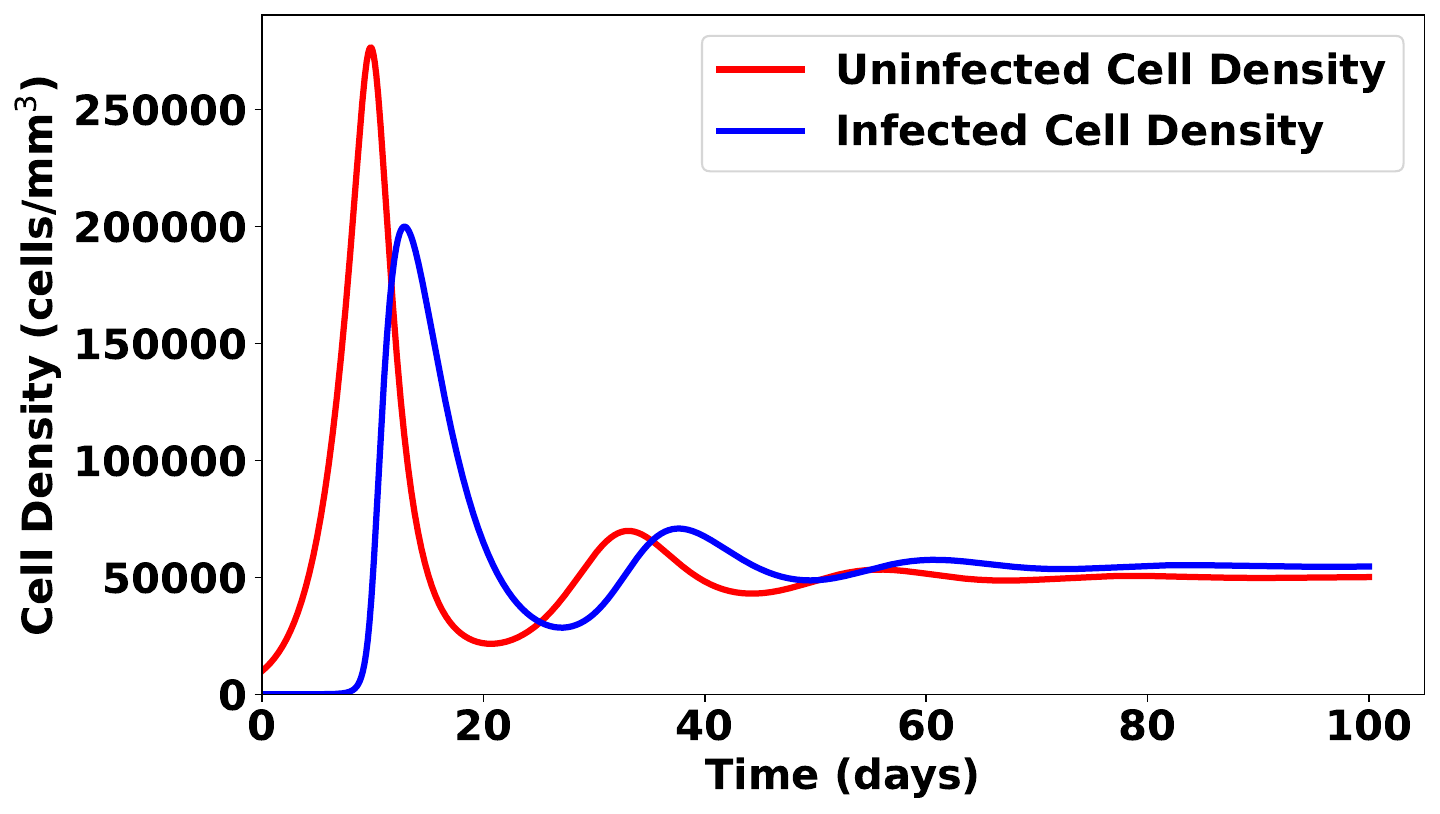}
    \caption{Tumour cell density dynamics: constant $\theta$ and $\gamma$}.
    \label{graph00}
\end{figure}

Table 2 gives the parameters value which we use in the case where $\gamma$ and $\theta$ are constants, rather than functions of oxygen concentration. In this case, plotting the tumour cell densities gives Figure \ref{graph00}. If we consider this to be the standard case, we can test the effect of including oxygen dependence of the functions $\theta$ and $\gamma$.

\begin{figure}[H]
\centering
	\begin{subfigure}[t]{0.47\textwidth}\centering
\includegraphics[width=\textwidth]{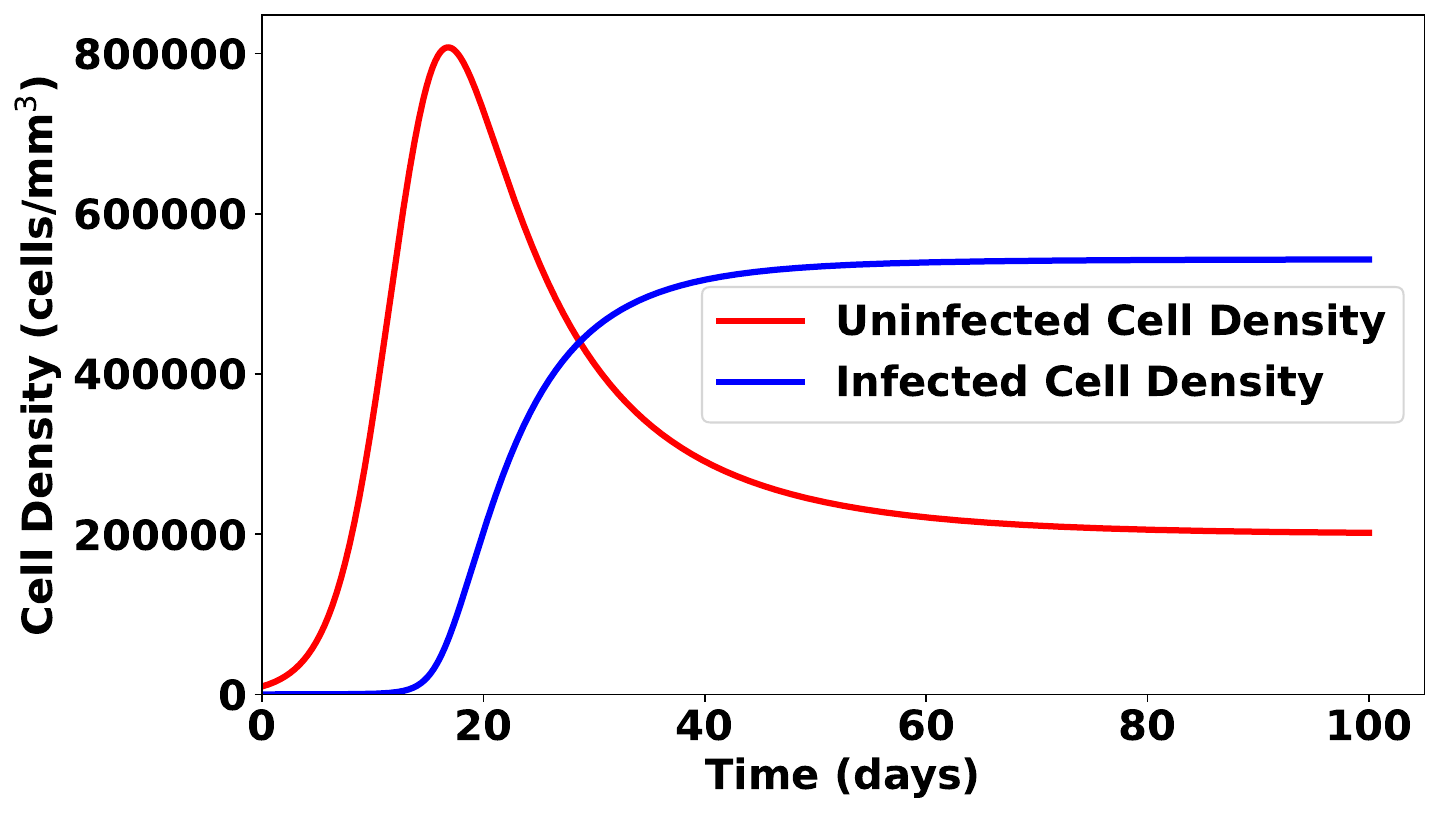}
         \caption{}\label{fig2_a}
    \end{subfigure}
    \hspace{5mm}
\begin{subfigure}[t]{0.45\textwidth}\centering
\includegraphics[width=\textwidth]{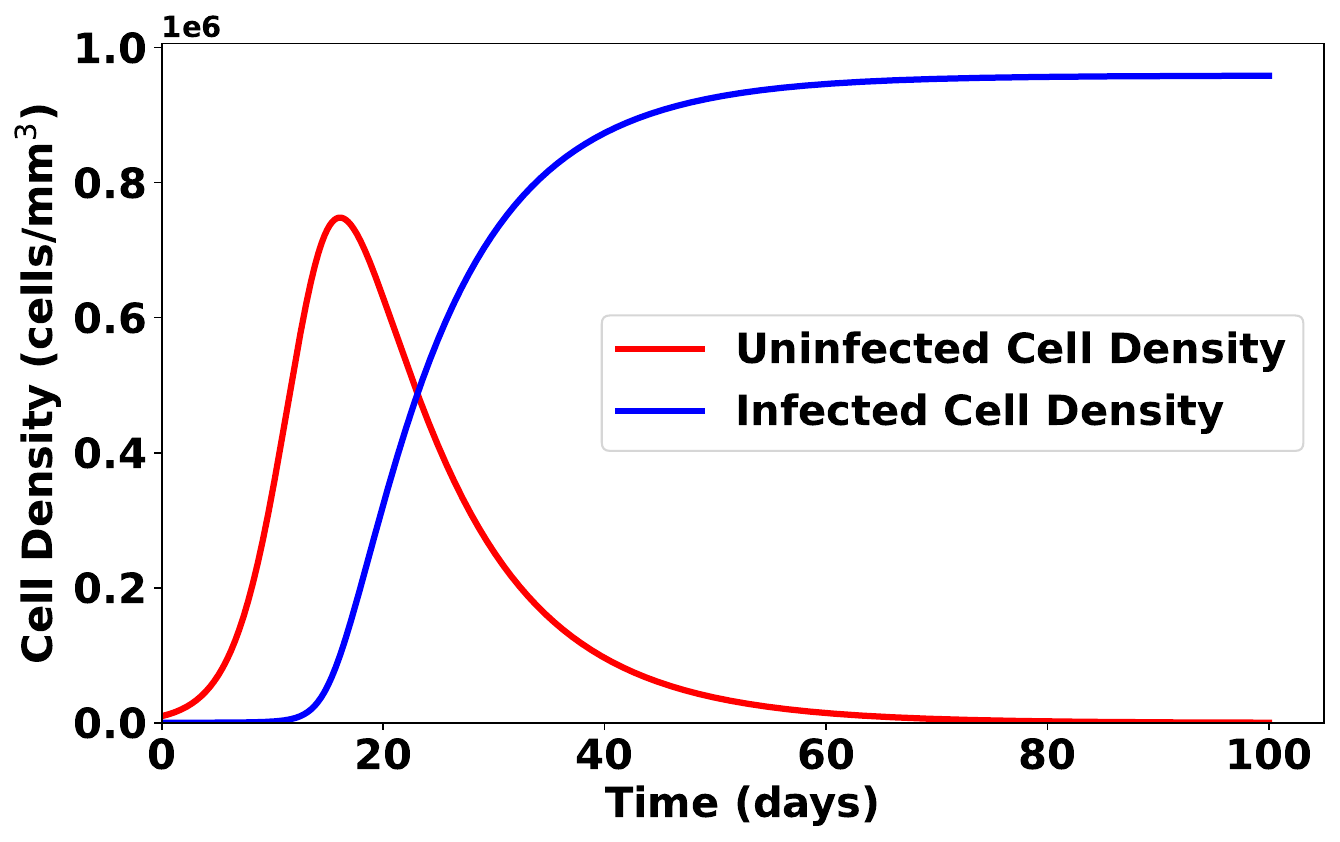}\\
 \caption{}\label{fig2_b}
\end{subfigure}
\caption{Tumour cell density dynamics in the case where $\theta(c) > (\alpha/K)\gamma(c)$ for all $c \geq 0$. (a) In this case, we have $\theta_0 = 0.1, \theta_\infty = 0.12, k_\theta = 0.08, \gamma_0 = 0.05115, \gamma_\infty = 0.09115, k_\gamma = 0.08$. The infected tumour cells dominate in the long run. This represents a relatively favourable response to the oncolytic virotherapy. (b) In this case, we have $\theta_0 = 0.1, \theta_\infty = 0.12, k_\theta = 0.08, \gamma_0 = 0.005115, \gamma_\infty = 0.009115, k_\gamma = 0.008$. When the virus-induced death rate is too low, infected cells still dominate but will ultimately approach a larger value at the steady state compared to the previous case. }
\label{graph01}
\end{figure}

\begin{figure}[H]
\centering
	\begin{subfigure}[t]{0.45\textwidth}\centering
\includegraphics[width=\textwidth]{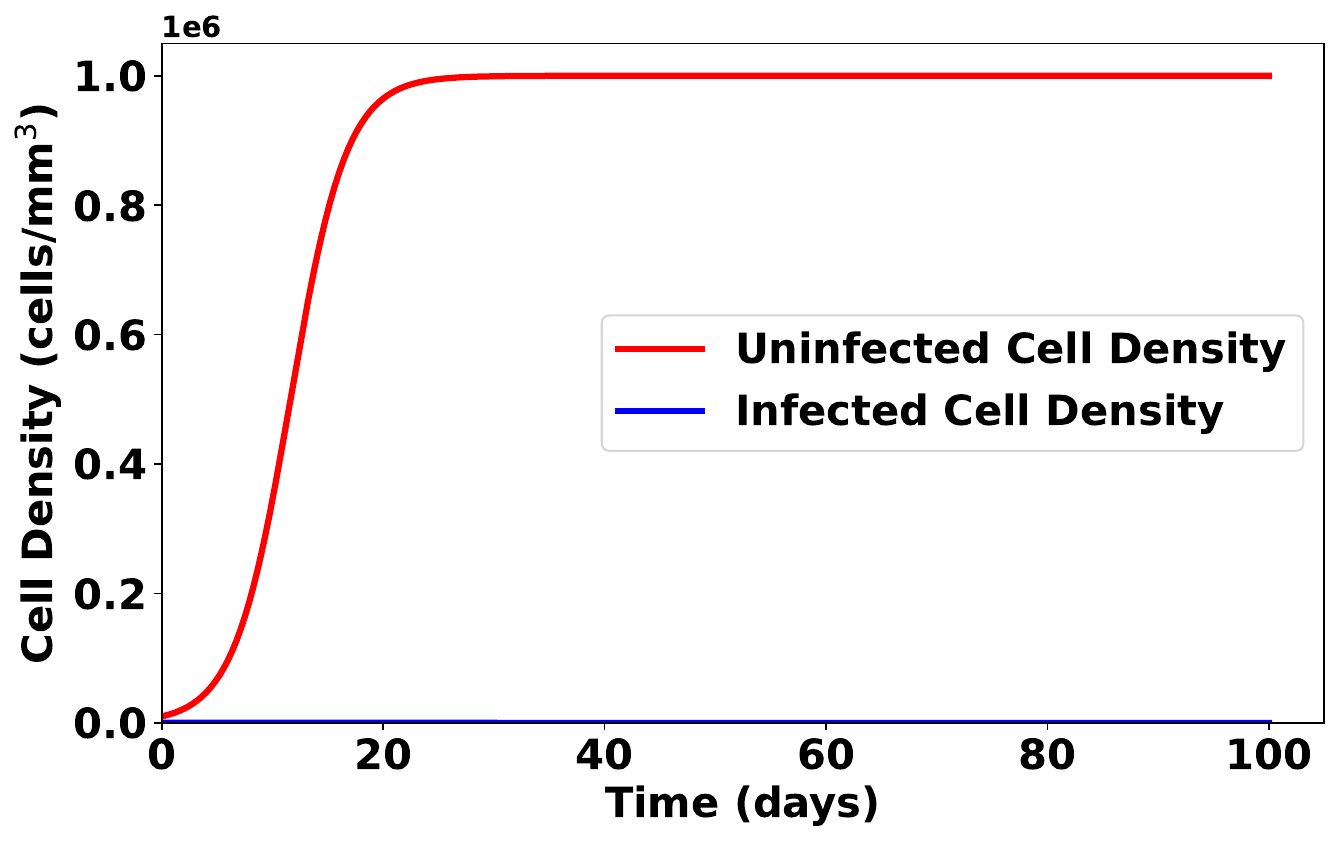}
         \caption{}\label{fig2_a}
    \end{subfigure}
    \hspace{5mm}
\begin{subfigure}[t]{0.45\textwidth}\centering
\includegraphics[width=\textwidth]{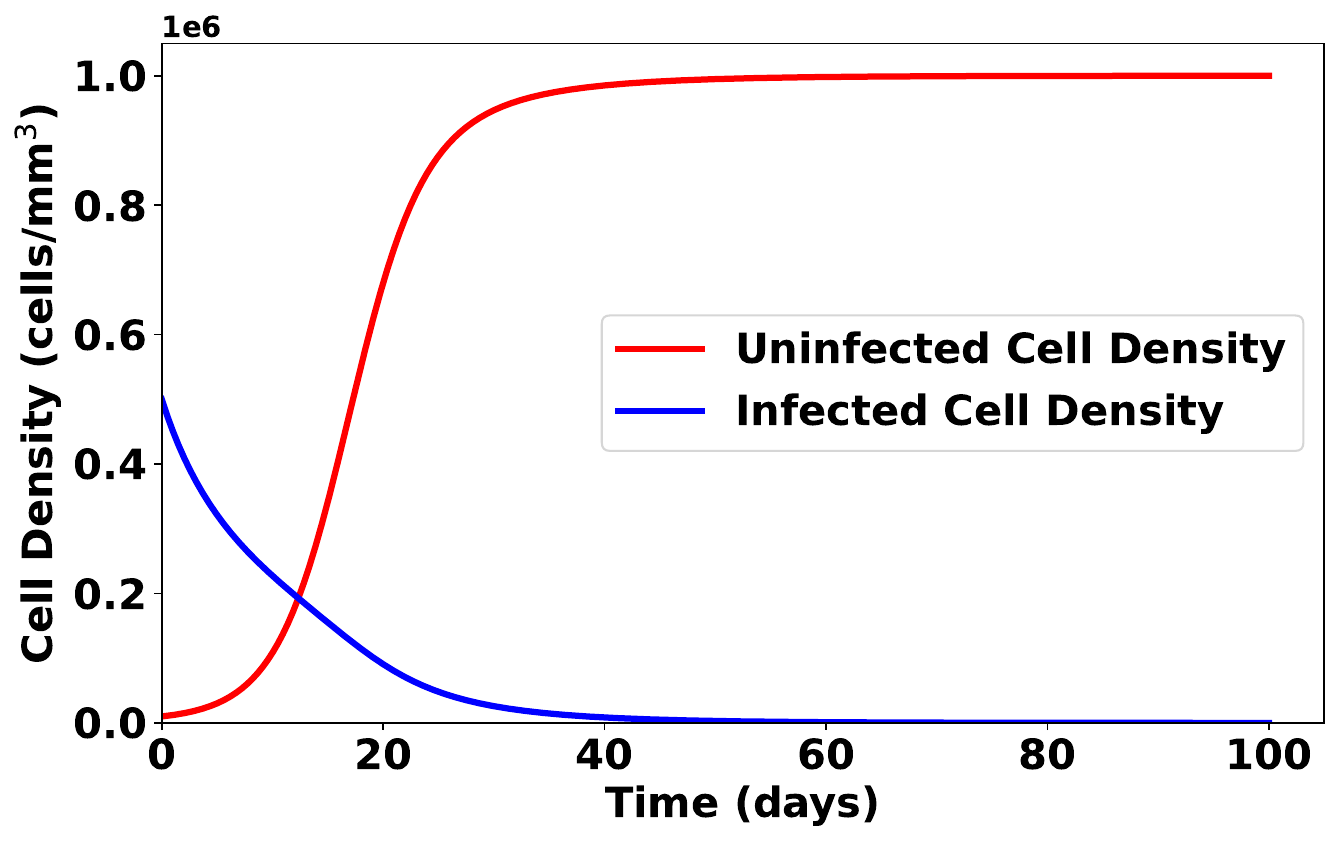}\\
 \caption{}\label{fig2_b}
\end{subfigure}
\caption{Tumour cell density dynamics in the case where $\theta(c) > (\alpha/K)\gamma(c)$ for $0 \leq c < c^*$ and  $\theta(c) < (\alpha/K)\gamma(c)$ for $c > c^*$. In this case, we have $\theta_0 = 0.01, \theta_\infty = 0.012, k_\theta = 0.008, \gamma_0 = 0.05115, \gamma_\infty = 0.2115, k_\gamma = 0.08$. (a) $n_0 = 100$. (b) $n_0 = 0.5 \times 10^6$.}
\label{graph02}
\end{figure}

\begin{figure}[H]
\centering
	\begin{subfigure}[t]{0.45\textwidth}\centering
\includegraphics[width=\textwidth]{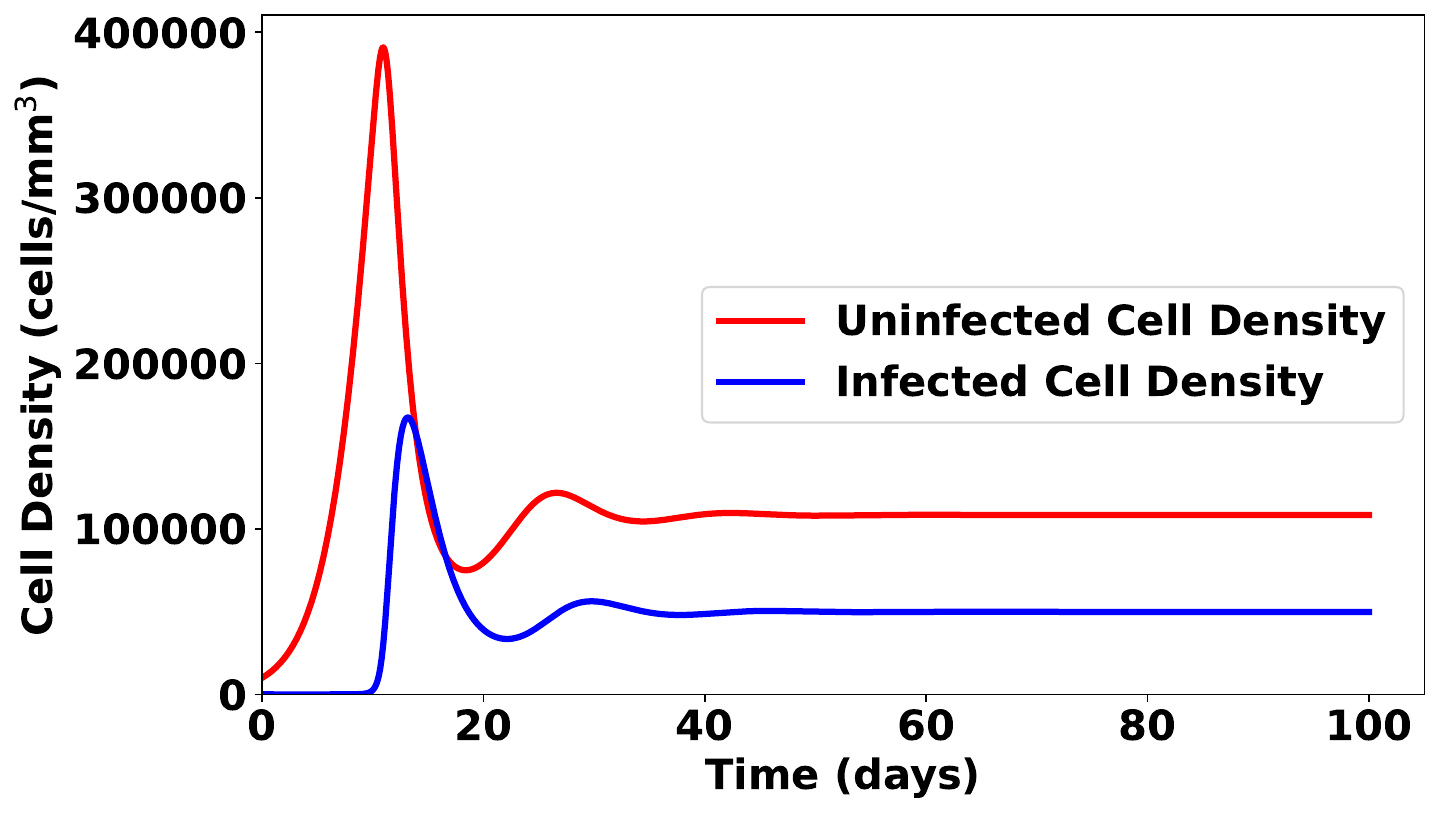}
         \caption{}\label{fig5_a}
    \end{subfigure}
    \hspace{5mm}
\begin{subfigure}[t]{0.45\textwidth}\centering
\includegraphics[width=\textwidth]{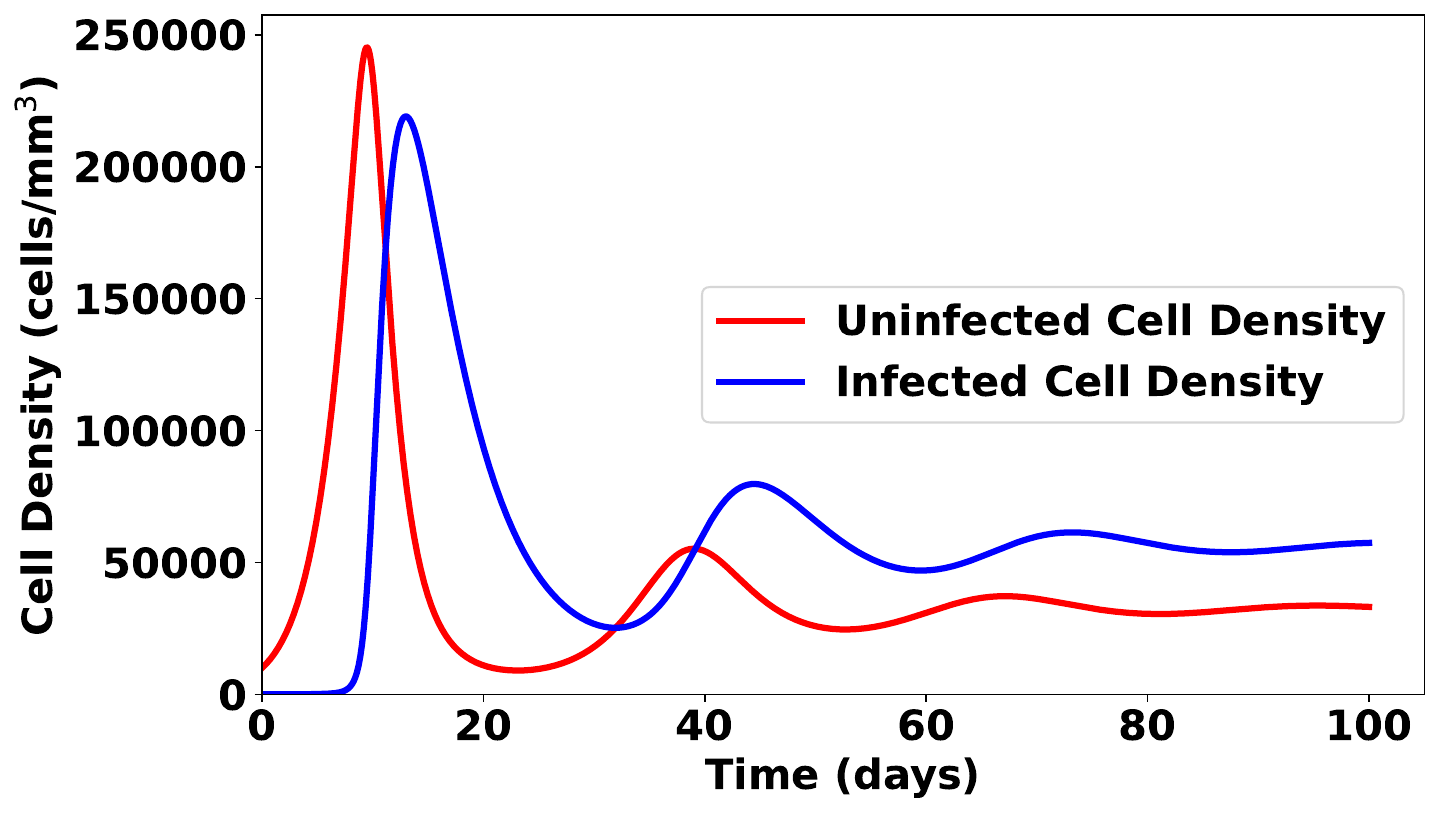}\\
 \caption{}\label{fig5_b}
\end{subfigure}
 \caption{Tumour cell density dynamics in the case where $\theta(c) < (\alpha/K)\gamma(c)$ for $0 \leq c < c^*$ and  $\theta(c) > (\alpha/K)\gamma(c)$ for $c > c^*$.  In this case, we have $\theta_0 = 5.115 \times 10^{-3}, \theta_\infty = 1.0, k_\theta = 0.08$ (a) $\gamma_0 = 0.2, \gamma_\infty = 0.4, k_\gamma = 0.08$. (b) $\gamma_0 = 0.7, \gamma_\infty = 0.9, k_{\gamma} = 0.08$.}
\label{graph03}
\end{figure}

\begin{figure}[H]
\centering
	\begin{subfigure}[t]{0.45\textwidth}\centering
\includegraphics[width=\textwidth]{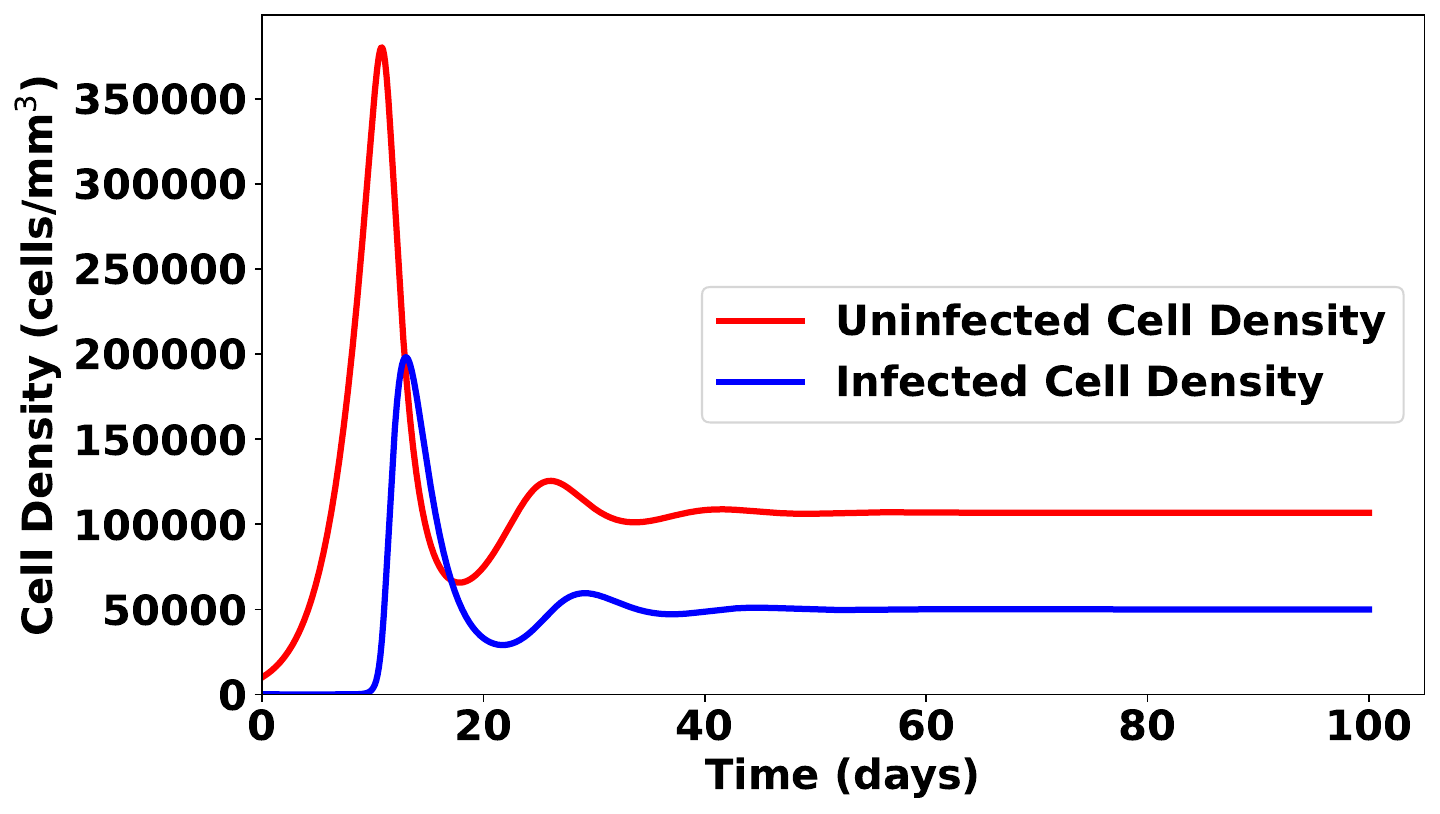}
         \caption{}\label{fig4_a}
    \end{subfigure}
    \hspace{5mm}
\begin{subfigure}[t]{0.45\textwidth}\centering
\includegraphics[width=\textwidth]{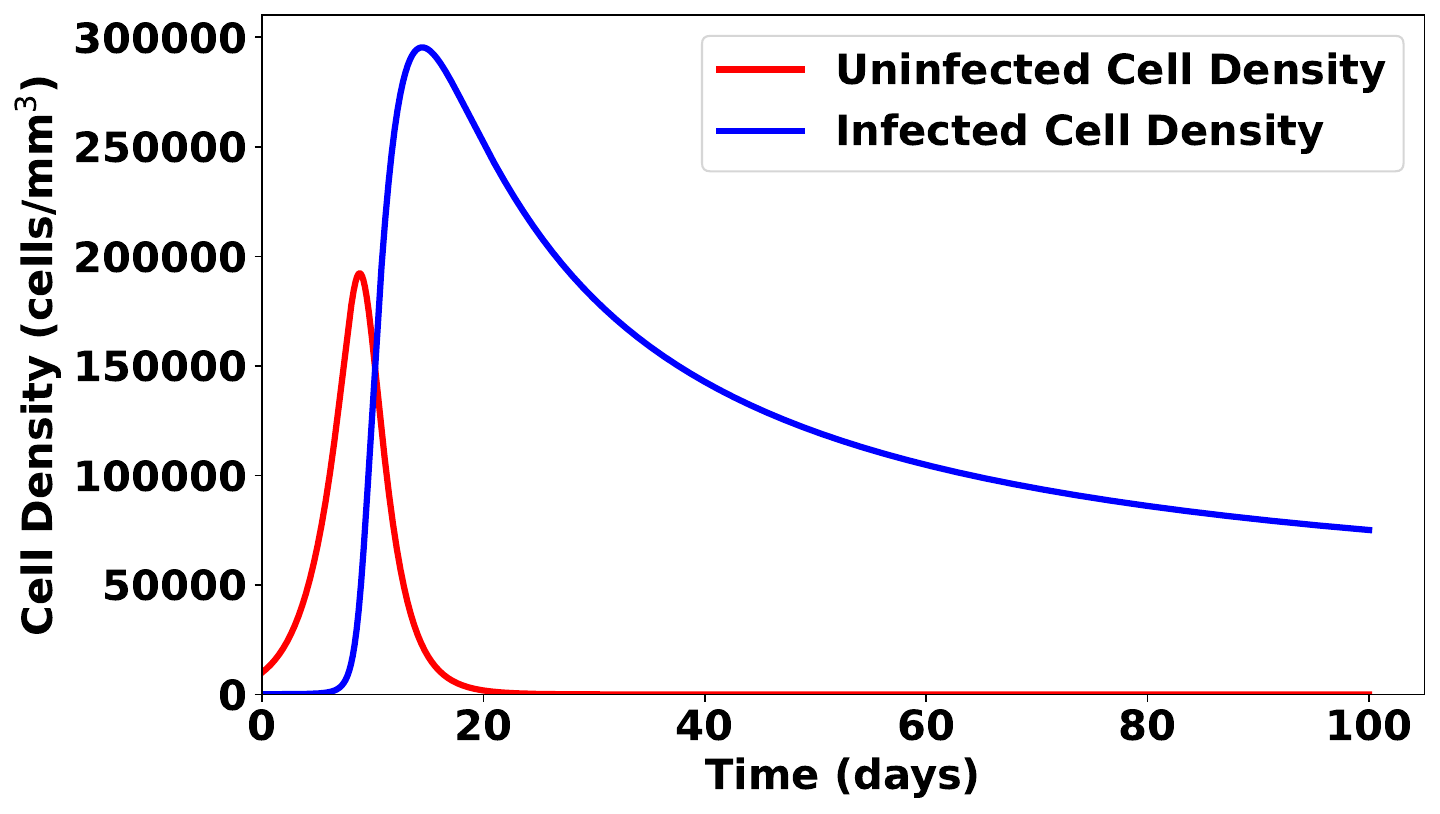}\\
 \caption{}\label{fig4_b}
\end{subfigure}
 \caption{Tumour cell density dynamics in the case where $\theta(c) < (\alpha/K)\gamma(c)$ for $0 \leq c < c^*$ and  $\theta(c) > (\alpha/K)\gamma(c)$ for $c > c^*$. In this case, we have $\gamma_0 = 0.1, \gamma_\infty = 0.9, \theta_0 = 5.115 \times 10^{-3}, \theta_\infty = 1.0, k_\theta = 0.08$ (a) $ k_\gamma = 0.008$. (b) $\gamma_0 = 0.09, \gamma_\infty = 0.2,  k_{\gamma} = 0.01$. }
\label{graph04}
\end{figure}

\begin{figure}[H]
    \centering
    \includegraphics[scale=0.35]{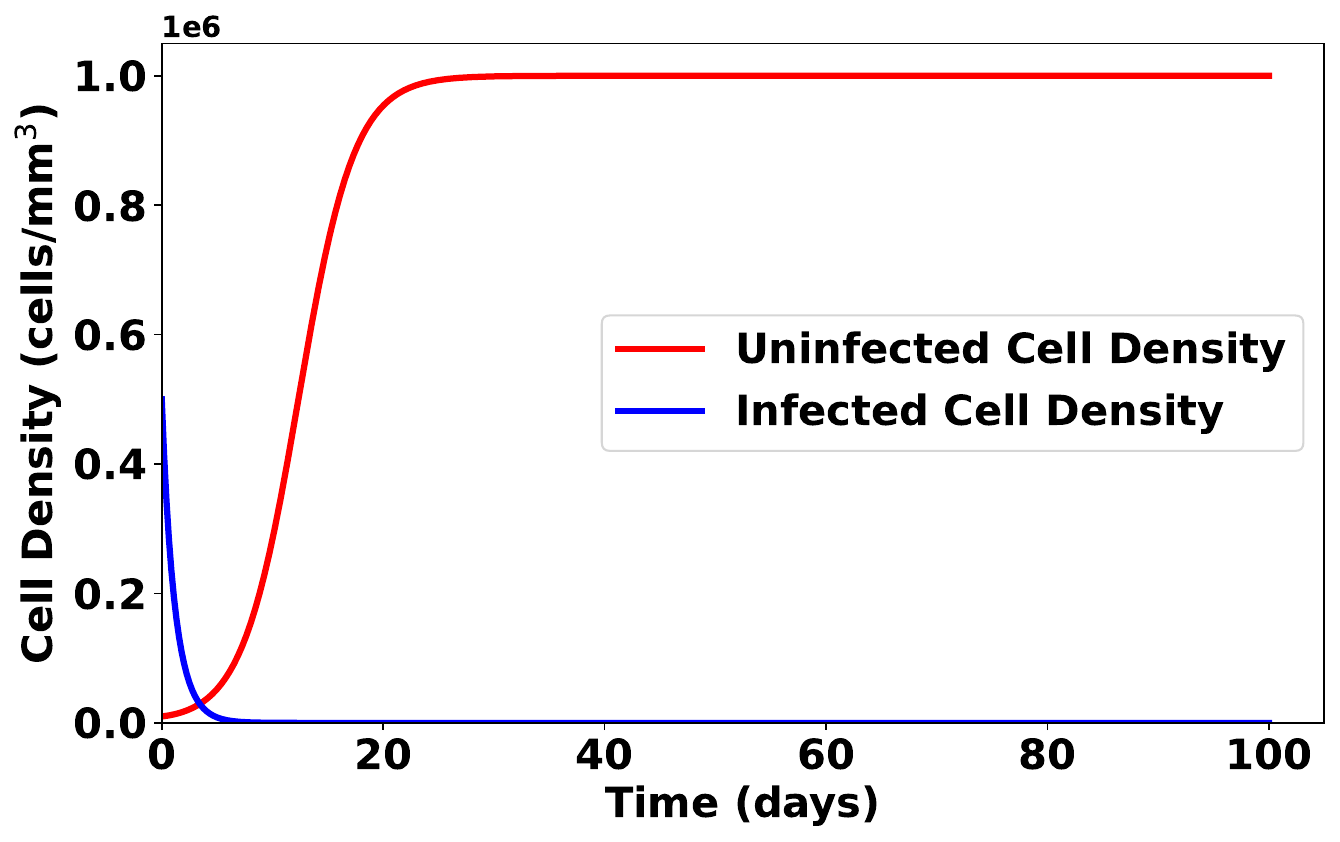}
    \caption{Tumour cell density dynamics in the case where $\theta(c) < (\alpha/K)\gamma(c)$ for all $c \geq 0$. In this case, we have $\gamma_0 = 0.3, \gamma_\infty = 1.0,k_\gamma = 0.8, \theta_0 = 0.005115, \theta_\infty = 0.02115, k_\theta = 0.8$. In this case, the virus-induced death rate is significantly greater than the infection rate. In this case, we set $n_0 = 5.0 \times 10^5$.}
    \label{graph05}
\end{figure}

In our simulations, we set $\theta$ and $\gamma$ to be sigmoid functions of $c$. In particular, we have
\begin{equation}
\theta(c) = \dfrac{\theta_\infty \theta_0}{\theta_0 + (\theta_\infty - \theta_0)e^{-k_\theta c}}, \quad
\gamma(c) = \dfrac{\gamma_\infty \gamma_0}{\gamma_0 + (\gamma_\infty - \gamma_0)e^{-k_\gamma c}}. \label{theta_gamma_equations1}
\end{equation}
We consider how different parameter values $\theta_0, \theta_\infty, \gamma_0, \gamma_\infty, k_\theta, k_\gamma$ impact the efficacy of the OV. Guided by Proposition \ref{prop2}, we choose these parameters such that we consider $\theta(c) < (\alpha/K)\gamma(c)$, $\theta(c) > (\alpha/K)\gamma(c)$, etc. We plot these results in Figures \ref{graph00} - \ref{graph05}. 

In Figure \ref{graph00}, we consider the case in which $\theta$ and $\gamma$ are the constants given in Table 2 rather than functions of the oxygen concentration. This is our first numerical exposure to a result which will be echoed throughout this subsection: higher infection rates relative to virus-induced death rates tend to lead to more favourable clinical results. In this case, the tumour cell densities both settle to a steady state well below the carrying capacity, suggesting some inhibition of the growth of the tumour cells.

In Figure \ref{graph01}, we have the case of a high infection rate relative to the virus-induced death rate. The assumption of Proposition \ref{prop2} is not satisfied and, unsurprisingly, the uninfected cell density is driven below the infected cell density, asymptotically. This case potentially represents a favourable result since in Figure \ref{graph01} (a), as the tumour cell density approaches a positive stable steady state value below the carrying capacity. In particular, the uninfected tumour cell density remains significantly lower than the infected tumour cell density. This illustrates the importance of the infection rate being sufficiently large. On other hand, as in Figure \ref{graph01} (b), having the virus-induced death rate be \textit{too low} leads to an unfavourable result in which all the tumour cells are infected but they nevertheless approach a value \textit{near} the carrying capacity -- note that they do not approach the carrying capacity in the case depicted by the figure. This illustrates the delicate balance between viral infection and virus-induced mortality. Furthermore, we note the differences between Figure \ref{graph01} and Figure \ref{graph00}: In both cases, $\theta > (\alpha/K)\gamma$, yet the dynamics are qualitatively different. This difference is a result of Figure \ref{graph01} depending on oxygen concentration; a consideration not made in Figure \ref{graph00}. 

On the other hand, Figure \ref{graph02} shows a clinically unfavourable result. Namely, the uninfected tumour cell density approaches the carrying capacity value while the infected tumour cells die out. In this case, treatment via OV has failed. This occurs when the $\theta$ and $(\alpha/K)\gamma$ curves intersect at some oxygen value, $c^*$. The outcome of the numerics, in this case, directly follows from Proposition \ref{prop2}. Regardless of the initial density of OV injection, $n_0$, (i.e., Figure \ref{graph02} (a) vs. Figure \ref{graph02} (b)) the asymptotic behaviour is the same. Biologically, this gives the following insight: in hypoxic environments, having very low lysis capabilites of the OV yields failure of the treatment regardless of initial density of the OV injection. It is worthwhile to note that the $z^*$ from the steady state considered in Proposition \ref{prop2} is \textbf{NOT} related to the quantity $c^*$, the $c$-coordinate of the intersection point of $\theta$ and $\gamma$.

Figure \ref{graph03} represents the reverse case of Figure \ref{graph02}, in which the inequalities are reversed and the results are clinically more favourable. This once again illustrates the importance of the virus-induced death rate in hypoxic environments and also shows the importance of the infection rate in oxygen-rich environments. Moving from Figure \ref{graph03} (a) to Figure \ref{graph03} (b), we increase the virus-induced death rate, while still maintaining a high viral  infection rate in oxygen rich conditions. This further supports the idea of achieving a balance between infection rates and lysis capacity as an OV engineering consideration.

In Figure \ref{graph04}, for low values of oxygen concentration (i.e., hypoxic environments) the infection rate is significantly less than the virus-induced death rate, whereas for high values  of  oxygen concentration, the virus-induced death rate is reduced. The figure shows that this case also represents a favourable clinical outcome represented by the dampening oscillations in Figure \ref{graph04} (a). Asymptotically, the tumour cell density approaches a positive steady state value well below the carrying capacity. This (once again) suggests the following insight: in hypoxic environments, it is important that the OV is more efficient at killing cancer cells than infecting them. However, if the oxygen concentration should be large, the OV must be more efficient at infecting tumour cells than inducing their death. In Figure \ref{graph04} (b), we decrease the growth rate of the $\gamma(c)$ function, leading to near-extinction of all tumour cells. Biologically, this represents an OV which has greater tumour-destroying capabilities over a lesser range of lower oxygen concentrations. Another interpretation is that it would be favourable for the infection rate to surpass the virus-induced death rate at lesser oxygen concentrations as long as the virus-induced death rate does initially dominates under \textit{extremely} hypoxic conditions. Such a virus must be engineered to initially be extremely potent at destroying tumour cells when there is almost no oxygen available in the tumour microenvironment but must quickly be able to adapt by having a much greater infection rate if the available oxygen concentration should increase. These results are consistent with Proposition \ref{prop_gamma_prime}.

Figure \ref{graph05} shows the case where the infection rate is very low compared to the virus-induced death rate. In this case, the uninfected tumour cell density dominates and approaches the carrying capacity. This result agrees with Proposition \ref{prop2}. This further supports the idea of a delicate balance between how effective the virus is at infected cancer cells and how potent the virus is at inducing death of tumour cells. In particular, we must ensure that the death rate is not too large compared to the infection rate.

These cases illustrate the following point which must be considered when engineering the OV: Having a virus too efficient at destroying and not efficient enough at infecting is not recommended. Perhaps equally importantly, we must also consider the oxygen conditions (i.e., hypoxia)  when engineering the OV as the functionality of the virus also depends on whether or not the tumour microenvironment is hypoxic.

\section{Regional oncolytic virotherapy model}

In this section, we extend our model to the regional setting by considering the case of lymph node invasion by the tumour cells. {Since movement through the lymphatic system is one of the main methods through which melanoma tumour cells may spread, it is of vital importance to consider lymph nodes as part of the model system. In the context of hypoxia, there is evidence which suggests that hypoxic conditions contribute to the upregulation of uPAR (a receptor on the surface of melanoma cells), leading to lymph node metastasis of the tumour cells \cite{Lymph_Hypoxia_Paper}. Hence, it is both of mathematical and biological interest to capture the dynamics of regional (i.e., lymphatic) spread of tumour cells, as well as OV efficacy, under various oxygen conditions.} 

As the thickness of the melanoma tumour increases, there is an increased probability of the tumour spreading to nearest lymph nodes \cite{lymphprob1}. We model a network of lymph nodes as a one-dimensional lattice, where each node represents a lymph node and the edges represent lymphatic vessels.  

\tikzset{int/.style={draw, line width = 1mm, minimum size=8em}}
\usetikzlibrary{positioning}

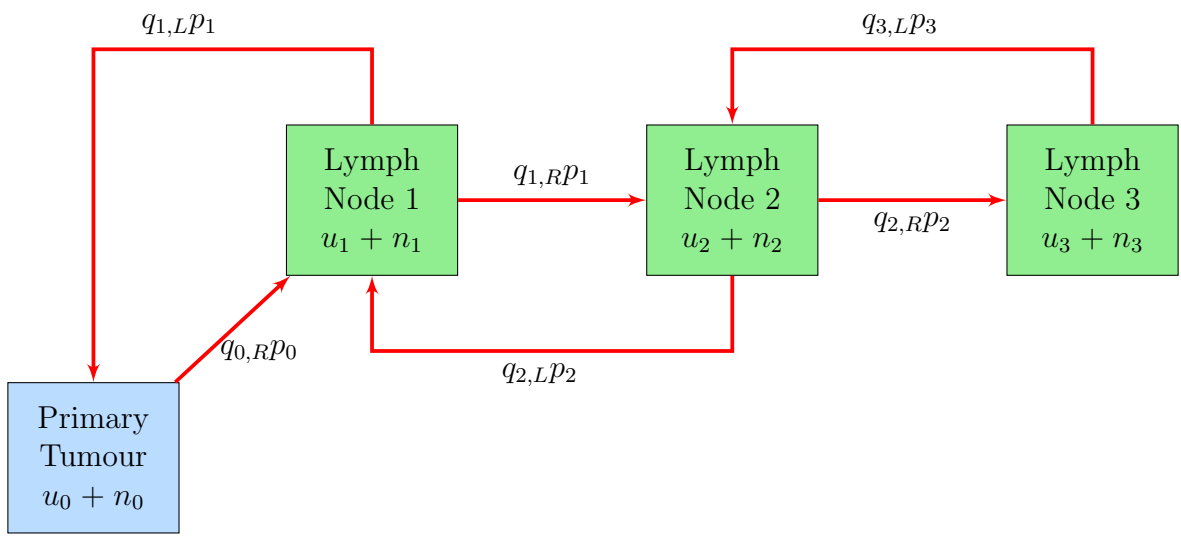
\begin{figure}[H]
    \centering
    \tikzstyle{block1}=[draw, fill=boxcolor_tumour, minimum size=2em, text width = 2.5cm, align = center, minimum height = 2.5cm]
\tikzstyle{block2}=[draw, fill=boxcolor_lymph, minimum size=2em, text width = 2.5cm, align = center, minimum height = 2.5cm]
\tikzstyle{arrow}=[->, red, text = black, line width=0.5mm]

 \tikzstyle{block1}=[draw, fill=boxcolor_tumour, minimum size=2em, text width = 2.0cm, align = center, minimum height = 2.0cm]
\tikzstyle{block2}=[draw, fill=boxcolor_lymph, minimum size=2em, text width = 2.0cm, align = center, minimum height = 2.0cm]
\tikzstyle{arrow}=[->, red, text = black, line width=0.5mm]

\begin{tikzpicture}[auto, >=latex']
    \node [block1] (t) {Primary Tumour \\ $u_0 + n_0$};
    \node [block2, above right = 2cm of t] (l1)  {Lymph Node 1 \\ $u_1 + n_1$};
	\node [block2, right = 2.5cm of l1] (l2) {Lymph Node 2 \\ $u_2 + n_2$};
	\node [block2, right = 2.5cm of l2] (l3) {Lymph Node 3 \\ $u_3 + n_3$};

    \draw [arrow] (t) -- node [left] {$$} node [right,pos=0.1mm] {$q_{0,R} p_0$} (l1);
\draw [arrow] (l1) -- node [left] {$$} node [above] {$q_{1,R} p_1$} (l2);
\draw [arrow] (l1.north) |- ++(0cm,+1.0cm) -- node [above,pos=0.3mm] {$q_{1,L} p_1$} ++(-3cm,0cm) -| (t);
\draw [arrow] (l2.south) |- ++(0cm,-1.0cm) -- node [below,pos=0.3mm] {$q_{2,L} p_2$} ++(-3cm,0cm) -| (l1);
\draw [arrow] (l2) -- node [left] {$$} node [below] {$q_{2,R} p_2$} (l3);
\draw [arrow] (l3.north) |- ++(0cm,+1.0cm) -- node [above, pos = 0.3mm] {$q_{3,L} p_3$} ++(-3cm,0cm) -| (l2);

\end{tikzpicture}
    \caption{The first three lymph nodes in a network. Each lymph node represents a different node (green) in the lattice. The tumour cells begin in the primary tumour (blue) and can travel through the lymphatic network, with some probability of spreading.}
    \label{graph0}
\end{figure}

The initial concentration of tumour cells at each node is set to $0$. Let $i$ denote the $i^{\text{th}}$ node from the primary tumour for $i = 1, 2, 3, \dots, \ell$ and let $i = 0$ denote the primary tumour. That is, $i = 0$ corresponds to the local case presented in Section 3.2. Note that $u_0, n_0$ and $c_0$ no longer represent initial conditions, but rather the primary tumour. The tumour cells may either travel to the left or to the right of their current position. We assume that the probability of tumour cells spreading to the adjacent lymph nodes depends on the density of the tumour cells at the given node and, hence, on the sum $u_i(t) + n_i(t)$.  Let $P_i(u_i+n_i)$ be the spreading rate of some fraction of tumour cells away from node $i$ to an adjacent lymph node in the network. This fraction of cells which leaves a given node is dependent on the tumour cell density at the given node i.e., $u_i(t) + n_i(t)$.  

The probability of the cells travelling left is $q_{i,L}$ and the probability of travelling right is $q_{i,R}$, where $q_{i,L} + q_{i,R} = 1$ for $i = 1, 2, \dots, \ell-1$. Moreover, $q_{0,R} = 1$. While it has been observed that lymph typically flows only in one direction \cite{zimmerman_lymphfact}, we allow for the possibility of some tumour cells to travel in the reverse direction. We assume that the probability of tumour cells reversing direction is low and therefore, we consider $q_{i,R} >> q_{i,L}$ in the  numerical simulations. 

On each node in the network, we have a system of ODEs which describes the number of tumour cells and the oxygen concentration. We use system (\ref{e1}) - (\ref{e3})  to model the dynamics of the tumour cells at each individual node. To this end, we propose the following system:
\begingroup
\allowdisplaybreaks
\addtolength{\jot}{7pt}
\begin{align}
\dfrac{\mathrm{d}u_0}{\mathrm{d}t} &= r_1 u_0 \left( 1- \dfrac{u_0+n_0}{K_0}  \right) - \dfrac{\theta (c_0) n_0u_0}{\alpha_0 + n_0} - q_{0,R} P_0 u_0  +  q_{1,L} P_1 u_{1},\label{eqlattice1}\\
\dfrac{\mathrm{d}n_0}{\mathrm{d}t} &= r_2 n_0 \left( 1 - \dfrac{u_0+n_0}{K_0}  \right) + \dfrac{\theta (c_0) n_0u_0}{\alpha_0 + n_0} - \gamma (c_0) n_0 -  q_{0,R} P_0 n_0 +  q_{1,L} P_1 n_1, \label{eqlattice2}  \\
\dfrac{\mathrm{d}u_i}{\mathrm{d}t} &= r_1 u_i \left( 1- \dfrac{u_i+n_i}{K_i}  \right) - \dfrac{\theta (c_i) n_iu_i}{\alpha_i + n_i} +  q_{i-1,R} P_{i-1} u_{i-1} +  q_{i+1,L} P_{i+1} u_{i+1} -  P_i u_{i}, \label{eqlatticei1}\\
\dfrac{\mathrm{d}n_i}{\mathrm{d}t} &= r_2 n_i \left( 1 - \dfrac{u_i+n_i}{K_i}  \right) + \dfrac{\theta (c_i) n_iu_i}{\alpha_i + n_i} - \gamma (c_i) n_i +   q_{i-1,R} P_{i-1} n_{i-1} + \dots \notag \\ & \quad \quad \quad \quad \quad \quad \quad \quad \quad \quad \quad \quad \quad \quad \quad \quad \quad \quad \dots + q_{i+1,L}P_{i+1} n_{i+1}  -  P_i n_{i}, \label{eqlatticei2}  \\
\dfrac{\mathrm{d}u_\ell}{\mathrm{d}t} &= r_1 u_\ell \left( 1- \dfrac{u_\ell+n_\ell}{K_\ell}  \right) - \dfrac{\theta (c_\ell) n_{\ell}u_{\ell}}{\alpha_{\ell} + n_{\ell}}  +  q_{\ell-1,R} P_{\ell-1} u_{\ell-1}  -  q_{\ell,L} P_\ell u_{\ell}, \label{eqlatticel1}\\
\dfrac{\mathrm{d}n_\ell}{\mathrm{d}t} &= r_2 n_\ell \left( 1 - \dfrac{u_\ell+n_\ell}{K_\ell}  \right) + \dfrac{\theta (c_\ell) n_\ell u_\ell}{\alpha_\ell + n_\ell} - \gamma (c_\ell) n_\ell  + q_{\ell-1,R} P_{\ell-1} n_{\ell-1} - q_{\ell,L} P_{\ell} n_{\ell}, \label{eqlatticel2}  \\
\dfrac{\mathrm{d}c_k}{\mathrm{d}t} &= \phi_k - \beta c_k - q_1u_kc_k - q_2n_kc_k,  \label{eqlatticel3}
\end{align}
\endgroup
where $i = 1, 2, 3, \dots, \ell-1$, $k = 0, 1, 2, 3, \dots, \ell$, and $P_k := P_k (u_k + n_k)$. Note that $\ell$ is the number of lymph nodes in the network. We set $\phi_0 = \phi$ and $\phi_k = 0$ for $k = 1, 2, 3, \dots, \ell$ for the purpose of following an individual over the course of treatment. The amount of time which tumour cells spend in compartment $i$, in the case of a large density of tumour cells in the compartment, is $1/\eta_i$, in days. We may therefore consider $\eta_i$ to give a per capita \textit{spreading speed} of tumour cells \textit{away} from compartment $i$. We assume that the speed with which tumour cells travel to the next node is equal to the speed with which they travel to the previous node. In practice, $q_{i,R} \approx 1$ and $q_{i,L} \approx 0$, so this assumption is not typically needed. It will be necessary for the purpose of tractability of the subsequent mathematical analysis. In this work, we only consider the case of a linear lymphatic network.

With all these considerations, the spreading rate of tumour cells leaving the lymph node and spreading to adjacent nodes is given by
\begin{equation}
P_i(x) = \eta_i \left[1 -  e^{-\lambda_i x} \right], \quad \lambda_i > 0. \label{P_iprobequation}
\end{equation}
The rationale behind defining $P_i$ in such a way is based on experimental results relating the size of a primary lesion to the probability of the cancer reaching the sentinel lymph nodes. See, for example, \cite{lymphprob1}. For the purpose of simulations, we assume that once the carrying capacity is reached, the probability of spreading is 0.7. Hence, we take $(1/\eta_i) P_i(K_i) = 0.7$  and solve this equation to determine the value $\lambda_i$ to be $\lambda_i = -\ln(0.3)/(K_i)$. 

We can show that the regional model is also well-posed in the sense of existence, uniqueness, non-negativity of the solution of the corresponding initial value problem with non-negative initial conditions, and boundedness of solutions. We summarize this result in the following theorem.

\begin{theorem}
There exists a unique solution of the initial value problem (\ref{eqlattice1}) - (\ref{eqlatticel3}), with non-negative initial conditions, which remains non-negative and bounded for all $t \geq 0$. \label{regional_theorem}
\end{theorem}

\begin{proof}
Existence and unqiueness of solutions follow directly from the fundamental theory of ODEs. To address the non-negativity of solutions, we apply Theorem 2.1 in Chapter 5 of  \cite{zounotes}.

Let $(u_0 (t), n_0(t), u_1 (t), n_1 (t), \dots u_\ell (t), n_\ell (t), c_0 (t), c_1 (t), \dots, c_\ell(t)) \in \mathbb{R}^{3\ell + 3}$ be the solution of the initial value problem consisting of system (\ref{eqlattice1}) - (\ref{eqlatticel3}) with non-negative initial conditions. We begin by showing that none of the components $u_0, u_1, \dots, u_\ell$ become negative. Assume to the contrary that at some time, $t^*$, at least one of the components of the solution becomes negative. By continuity, these components must first cross 0. If $u_0$ is one of these components, then plugging in $u_0 = 0$ into equation (\ref{eqlattice1}) gives $$ \dfrac{\mathrm{d}u_0}{\mathrm{d}t} =  q_{1,L} u_1 P_1 (u_1 + n_1) \geq 0, $$ which implies that $u_0$ is non-decreasing at $t = t^*$. Therefore, $u_0$ cannot become negative, leading to a contradiction. The same argument can be used to show that none of the $u_i$ components may become negative. 
 
Similarly, this contradiction argument can be used to conclude the non-negativity of $n_i$ for $i = 0, 1, 2, \dots, \ell$.

Finally, it can be seen that for $k = 0, 1, 2, \dots, \ell$, equation (\ref{eqlatticel3}) gives
\begin{equation*}
c_k(t) = c_k (0) \exp \left[ -\int_0^t \beta + q_1 u_k(s) + q_2 n_k(s) \mathrm{d}s  \right]  +  \phi_k \int_0^t \exp \left[ -\int_s^t \left( \beta + q_1 u_k (\xi) + q_2 n_k (\xi) \right)  \mathrm{d} \xi \right] \mathrm{d}s,
\end{equation*}
from which non-negativity of $c_k (t)$ follows. 

Hence, for non-negative initial conditions, $u_i (0), n_i (0), c_i (0) \geq 0$, for $i = 0, 1, 2, \dots, \ell$, it follows that the solution of the initial value problem remains non-negative for all $t \geq 0$.

Next, we show that solutions of the regional model remain bounded. Define $U(t) := u_0 (t) + u_1 (t) + \dots + u_\ell (t)$. Then adding equations (\ref{eqlattice1}), (\ref{eqlatticei1}), and (\ref{eqlatticel1}) for $i = 1, 2, \dots, \ell-1$, gives
\begin{equation*}
    \dfrac{\mathrm{d}U}{\mathrm{d}t} \leq \sum_{i=0}^\ell r_1 u_i \left( 1 - \dfrac{u_i}{K} \right),
\end{equation*}
where $K := \max_{i=0,1, \dots, \ell} \{K_i\}$. Then,
\begin{align*}
    \dfrac{\mathrm{d}U}{\mathrm{d}t} &\leq r_1 \left[ (u_0 + u_1 + \dots + u_\ell) - \dfrac{u_0^2 + u_1^2 + \dots + u_\ell^2}{K} \right] \\[12pt]
    &\leq r_1 \left[ (u_0 + u_1 + \dots + u_\ell) - \dfrac{(u_0 + u_1 + \dots + u_\ell)^2}{(\ell+1)K} \right].
\end{align*}
The last inequality follows from the Cauchy-Schwarz inequality, namely, 
\begin{equation*}
    (\ell+1)(u_0^2 + u_1^2 + \dots + u_\ell^2) \geq (u_0 + u_1 + \dots + u_\ell)^2.
\end{equation*}
Hence, we have
\begin{equation*}
    \dfrac{\mathrm{d}U}{\mathrm{d}t} \leq r_1 U \left[ 1 - \dfrac{U}{(\ell+1)K} \right] \implies \limsup_{t \rightarrow \infty} U(t) \leq (\ell+1)K.
\end{equation*}
Therefore, the sum $U(t)$ is a bounded function. Since each component of the sum is non-negative, we conclude that each $u_i(t)$ is bounded for each $i = 0, 1, \dots, \ell$. 

We can similarly show that the infected tumour cells remain bounded at each node by defining $N(t) := n_0 (t) + n_1 (t) + \dots + n_\ell (t)$. Adding equations (\ref{eqlattice2}), (\ref{eqlatticei2}), and $(\ref{eqlatticel2})$ for $i = 1, 2, \dots, \ell-1$ gives
\begin{align*}
    \dfrac{\mathrm{d}N}{\mathrm{d}t} &\leq \sum_{i=0}^\ell \left[r_2 n_i \left( 1 - \dfrac{n_i}{K} \right) + \theta_\infty u_i \right] \\[12pt]
    &\leq r_2 N \left[ 1 - \dfrac{N}{(\ell+1)K} \right] + \theta_\infty \bar{U},
\end{align*}
where $\bar{U}$ is any upper bound for $U(t)$. The last inequality follows by applying the Cauchy-Schwarz inequality as in the previous case. Hence, 
\begin{equation*}
    \limsup_{t \rightarrow \infty} N(t) \leq \dfrac{r_2 (\ell+1)K + \sqrt{r_2^2 (\ell+1)^2K^2 + 4r_2 (\ell+1)K \theta_\infty \bar{U}}}{2r_2}.
\end{equation*}
Since the components of the sum $N(t)$ are non-negative, we conclude that $n_i(t)$ is bounded for each $i = 0, 1, \dots, \ell$.

Next, it is clear to see by a comparison argument that
\begin{equation*}
    \limsup_{t \rightarrow \infty} c_i (t) \leq \dfrac{\phi_k}{\beta}, \quad i = 0, 1, 2, \dots, \ell,
\end{equation*}
and so we may conclude that $c_i (t)$ are bounded.

We have successfully shown that solutions of the regional model with non-negative initial conditions are non-negative and bounded. 
\end{proof}

While an analytic investigation of system (\ref{eqlattice1}) - (\ref{eqlatticel3}) can be challenging to perform due to the potentially large number of equations, we may establish a result which is analogous to Proposition \ref{prop2} of the local model. In particular, we may establish the a sufficiently large oncolysis rate leads to stability of a tumour-dominant steady state. We begin by showing the existence of this steady state. 

Since we are also interested in obtaining results related to the spreading speed away from node $i$, i.e., $\eta_i$, we rewrite the function $P_i(x)$ from equation (\ref{P_iprobequation}) by defining the dimensionless quantity $p_i(x) := 1 - e^{-\lambda_i x}$, hence allowing us to formulate $P_i(x)$ in terms of the spreading speed. That is,
\begin{equation}
P_i(x) = \eta_i p_i (x).
\end{equation}
For the remainder of this subsection, we consider only the case where $\phi_k = 0$ for $k = 0, 1, 2, \dots, \ell$. Biologically, this condition corresponds to the case with no external oxygen input. Furthermore, we consider the case where a tumour cell may only travel forward through the network (i.e., in the right, $R$, direction). Hence, we set $q_{j,R} = 1$ and $q_{j,L}= 0$ for $j = 0, 1, 2, \dots, \ell$. This is biologically consistent with the unidirectional flow of tumour cells through the lymphatic system \cite{zimmerman_lymphfact}. 

A virus-free or \textit{tumour-dominant} steady state is of the form
\begin{equation}
   E_u := (u_0, n_0, u_1, n_1, \dots, u_\ell, n_\ell, c_0, c_1, \dots, c_\ell) = (u_0^*, 0, u_1^*, 0, \dots, u_\ell^*, 0, 0, 0, \dots, 0), \label{tumour_dominantsteadystate_regionalmodel}
\end{equation}
for $i \in \{ 0,1, \dots, \ell \}$, where $u_i^* > 0$.

It follows from system (\ref{eqlattice1}) - (\ref{eqlatticel3}) that the components of the tumour-dominant steady state satisfy the equations
\begingroup
\allowdisplaybreaks
\begin{equation}
    \begin{dcases}
        r_1 \left( 1 - \dfrac{u_0}{K_0} \right) &= \eta_0 \left(1 - e^{-\lambda_0 u_0}\right), \\[12pt]
        r_1 u_i\left( 1 - \dfrac{u_i}{K_i}  \right) &= -\eta_{i-1} u_{i-1} \left(1-e^{-\lambda_{i-1}u_{i-1}}\right) + \eta_i u_i \left(1-e^{-\lambda_{i}u_i}\right), \quad i \in \{ 1, 2, \dots, \ell - 1 \}, \\[12pt]
        r_1 u_\ell \left( 1 - \dfrac{u_\ell}{K_\ell} \right) &= -\eta_{\ell-1} u_{\ell - 1}\left( 1 - e^{-\lambda_{\ell-1}u_{\ell-1}} \right) .
    \end{dcases} \label{341system}
\end{equation}
\endgroup
The existence of the solution $u_0^* < K_0$ of the first equation of system (\ref{341system}) is clear. The solutions of the remaining equations of this system may subsequently be obtained by solving for $u_i^*$ recursively, given $u_{i-1}^*$. 

Based on a numerical exploration of the system (see Section  6), we also require that $u_i^* > K_i$ for $i \in \{1, 2, \dots, \ell \}$. It is trivial to see that this inequality holds for $i = \ell$. To ensure that this inequality is true for all other values of $i$, it is sufficient to consider the additional condition
\begin{equation}
    K_i < \dfrac{10 \eta_{i-1}}{7\eta_i}u_{i-1}^* \left( 1 - e^{-\lambda_{i-1}u_{i-1}^*} \right), \quad i \in \{ 1, 2, \dots, \ell - 1 \}.
\end{equation}
These upper bounds on $K_i$ come from system (\ref{341system}) and from $p_i (K_i) = 7/10$. They may be obtained recursively given $u_{i-1}^*$. 

Let $J = [J_{ij}] \in \mathbb{R}^{(3\ell + 3) \times (3\ell + 3)}$ be the Jacobian matrix of system (\ref{eqlattice1}) - (\ref{eqlatticel3}).

We are now in a position to establish the stability of $E_u$. To do so, we make use of the Gershgorin Disc Theorem \cite{bib_weisstein_chap3} which is stated as follows.

\begin{lemma}[Gershgorin Disc Theorem \cite{bib_weisstein_chap3}]
Consider an $n \times n$ matrix $A = [A_{ij}]$ in $\mathbb{C}^{n \times n}$. Define
\begin{equation*}
    R_i := \sum_{\substack{j=1\\[2pt]
                  j \not= i}}^n \left| A_{ij} \right|, \quad \left| z \right| \  \text{is the modulus of $z \in \mathbb{C}$}.
\end{equation*}
If $\lambda \in \mathbb{C}$ is an eigenvalue of $A$, then 
\begin{equation*}
    \lambda \in {\bigcup^n_{i=1}} \left\{ z \in \mathbb{C} : \left| z - A_{ii} \right| \leq R_i \right\}.
\end{equation*}\label{gershgorinlemma}
\end{lemma}
The circles of the form $\{ z \in \mathbb{C}: \left| z - A_{ii} \right| \leq R_i \} \subset \mathbb{C}$ in Lemma \ref{gershgorinlemma} are also called  \textit{Gershgorin discs}. Since all eigenvalues of $A$ are contained in these discs, we may bound the real part of these eigenvalues above by $0$ by ensuring that all of the Gershgorin discs lie in the left half of the complex plane.

We use the following approach in order to find sufficient conditions for the local asymptotic stability of $E_u$: 
\begin{enumerate}
    \item Linearize system (\ref{eqlattice1}) - (\ref{eqlatticel3}) at the steady state $E_u$. Let $J(E_u) = [J_{ij}(E_u)]$ denote this matrix.
    \item For all $i$, compute $R_i$ by adding the absolute value of all of the off-diagonal elements in row $i$ of $J(E_u)$, as in Lemma \ref{gershgorinlemma}.
    \item Find conditions (if any) such that $\forall i \in \{ 1, 2, \dots, 3\ell + 3 \} : J_{ii}(E_u) + R_i < 0$. If this is possible, then all of the eigenvalues of $J(E_u)$ have negative real part and hence, $E_u$ is locally asymptotically stable.
\end{enumerate}

For notational convenience, note that system (\ref{eqlattice1}) - (\ref{eqlatticel3}) may be written in the form
\begingroup
\allowdisplaybreaks
\begin{align*}
    \dfrac{\mathrm{d}u_i}{\mathrm{d}t} &= \mathcal{U}_i (u_0, n_0, u_1, n_1,\dots, u_\ell, n_\ell, c_0, c_1, \dots, c_{\ell}), \\[12pt]
    \dfrac{\mathrm{d}n_i}{\mathrm{d}t} &= \mathcal{N}_i (u_0, n_0, u_1, n_1,\dots, u_\ell, n_\ell, c_0, c_1, \dots, c_{\ell}), \\[12pt]
    \dfrac{\mathrm{d}c_i}{\mathrm{d}t} &= \mathcal{C}_i (u_0, n_0, u_1, n_1,\dots, u_\ell, n_\ell, c_0, c_1, \dots, c_{\ell}),
\end{align*}
\endgroup
where $i = 0, 1, 2, \dots, \ell$, for appropriately defined functions $\mathcal{U}_i, \mathcal{N}_i,$ and $\mathcal{C}_i$.

We begin by noting that the diagonal elements of $J(E_u)$ are
\begingroup
\allowdisplaybreaks
    \begin{align*}
    J_{jj}(E_u) &= \begin{dcases} \dfrac{\partial \mathcal{U}_i}{\partial u_i} \bigg|_{E_u} &= r_1 - \dfrac{2r_1}{K_i}u_i^* - \eta_i \left[ p_i(u_i^*) + \lambda_i u_i^* e^{-\lambda_i u_i^*} \right], \quad j = 2i+1, \quad i \in \{ 0,1,\dots,\ell-1 \}, \\[12pt] \dfrac{\partial \mathcal{N}_i}{\partial n_i} \bigg|_{E_u} &= r_2 - \dfrac{r_2}{K_i}u_i^* + \dfrac{\theta_0 u_i^*}{\alpha_i} - \gamma_0 - \eta_i p_i(u_i^*), \quad \quad \ \ j = 2i+2, \quad i \in \{ 0, 1,\dots,\ell-1 \}, \\[12pt] 
    \dfrac{\partial \mathcal{U}_\ell}{\partial u_\ell} \bigg|_{E_u} &= r_1 - \dfrac{2r_1}{K_\ell}u_{\ell}^*, \ \ \qquad \qquad \qquad \qquad  \qquad \qquad j = 2\ell + 1,\\[12pt]
    \dfrac{\partial \mathcal{N}_\ell}{\partial n_\ell} \bigg|_{E_u} &= r_2 - \dfrac{r_2}{K_\ell}u_{\ell}^* + \dfrac{\theta_0 u_{\ell}^*}{\alpha_\ell} - \gamma_0, \quad  \qquad \ \ \quad \ \qquad \ j = 2\ell + 2,\\[12pt]
    \dfrac{\partial \mathcal{C}_i}{\partial c_i} \bigg|_{E_u} &=  - \beta - q_1 u_i^*, \qquad \qquad \qquad \qquad \qquad \quad \quad \ \ j = 2\ell + 3 + i, \quad i \in \{ 0, 1, \dots \ell \}, \end{dcases}
    \end{align*}
\endgroup
To ensure that the eigenvalues lie in the left half of the complex plane, it is sufficient to find conditions such that $J_{ii}(E_u) + R_i < 0$ for all $i$. 

We begin by computing $J_{11}(E_u) + R_1$ which yields
\begingroup
\allowdisplaybreaks
\begin{align*}
    J_{11}(E_u) + R_1 &= \dfrac{\partial \mathcal{U}_0}{\partial u_0} \bigg|_{E_u}  + \left| \dfrac{\partial \mathcal{U}_0}{\partial n_0} \bigg|_{E_u}  \right| + \left| \dfrac{\partial \mathcal{U}_0}{\partial c_0} \bigg|_{E_u}  \right| + \sum_{j=1}^{\ell} \left( \left| \dfrac{\partial \mathcal{U}_0 }{\partial u_j} \bigg|_{E_u} \right| + \left| \dfrac{\partial \mathcal{U}_0 }{\partial n_j} \bigg|_{E_u} \right| + \left| \dfrac{\partial \mathcal{U}_0 }{\partial c_j} \bigg|_{E_u} \right| \right) \\[15pt]
    &= r_1 - \dfrac{2r_1}{K_0}u_0^* - \eta_0 \left[ p_0(u_0^*) + \lambda_0 u_0^* e^{-\lambda_j u_0^*} \right] + \dfrac{r_1}{K_0}u_0^* + \dfrac{\theta_0 u_0^*}{\alpha_0} + \eta_0 \lambda_0 u_0^* e^{-\lambda_0 u_0^*}  \\[15pt]
    &= r_1 \left(1- \dfrac{u_0^*}{K_0}\right) + \dfrac{\theta_0 u_0^*}{\alpha_0} - \eta_0 p_0 (u_0^*).
\end{align*}
\endgroup
This quantity is negative for a sufficiently large spreading speed away from the primary tumour site, $\eta_0$. In particular, this is true if 
\begin{equation}
    \eta_0 > \dfrac{1}{p_0 (u_0^*)}\left[ r_1 \left( 1 - \dfrac{u_0^*}{K_0} \right) + \dfrac{\theta_0 u_0^*}{\alpha_0} \right]. \label{condition_for_eta0}
\end{equation}
Next, computing $J_{22}(E_u) + R_2$ yields
\begingroup
\allowdisplaybreaks
\begin{align*}
    J_{22}(E_u) + R_2 &=  \dfrac{\partial \mathcal{N}_0}{\partial n_0} \bigg|_{E_u}  + \left| \dfrac{\partial \mathcal{N}_0}{\partial u_0} \bigg|_{E_u}  \right| + \left| \dfrac{\partial \mathcal{N}_0}{\partial c_0} \bigg|_{E_u}  \right| + \sum_{j=1}^{\ell} \left( \left| \dfrac{\partial \mathcal{N}_0 }{\partial u_j} \bigg|_{E_u} \right| +  \left| \dfrac{\partial \mathcal{N}_0 }{\partial n_j} \bigg|_{E_u} \right| + \left| \dfrac{\partial \mathcal{N}_0 }{\partial c_j} \bigg|_{E_u} \right| \right) \\[15pt]
    &= r_2 - \dfrac{r_2}{K_0}u_0^* + \dfrac{\theta_0 u_0^*}{\alpha_0} - \gamma_0 - \eta_0 p_0(u_0^*).
\end{align*}
\endgroup
Condition (\ref{condition_for_eta0}) is sufficient for the negativity of this quantity since $r_1 > r_2$ and hence, no additional conditions are necessary. 

We now consider the case $\ell > 1$, i.e., there are at least two lymph nodes in the network. Let $\mathcal{I} := \{1, 2, \dots, \ell-1\}$. For $i \in \mathcal{I}$, we define $k:= 2i+1$. Then we have
\begingroup
\allowdisplaybreaks
\begin{align*}
    J_{kk}(E_u) + R_k &= \dfrac{\partial \mathcal{U}_i}{\partial u_i} \bigg|_{E_u}  + \left| \dfrac{\partial \mathcal{U}_i}{\partial n_i} \bigg|_{E_u}  \right| + \left|  \dfrac{\partial \mathcal{U}_i }{\partial c_i} \bigg|_{E_u}  \right| +  \sum_{j \in \mathcal{I} \setminus \{i\}} \left( \left| \dfrac{\partial \mathcal{U}_i }{\partial u_j} \bigg|_{E_u} \right| + \left| \dfrac{\partial \mathcal{U}_i }{\partial n_j} \bigg|_{E_u} \right| + \left| \dfrac{\partial \mathcal{U}_i }{\partial c_j} \bigg|_{E_u} \right| \right) \\[15pt]
    &= r_1 - \dfrac{2r_1}{K_i}u_i^* - \eta_i \left[ p_i(u_i^*) + \lambda_i u_i^* e^{-\lambda_i u_i^*} \right] + \dfrac{r_1}{K_i}u_i^* + \dfrac{\theta_0 u_i^*}{\alpha_i} + \dots \\[15pt]
    & \quad \quad \quad \dots + \eta_i \lambda_i u_i^* e^{-\lambda_i u_i^*} + \eta_{i-1} \left[   p_{i-1}(u_{i-1}^*) + 2\lambda_{i-1}u_{i-1}^*e^{-\lambda_{i-1}u_{i-1}^*} \right] \\[15pt]
    &= r_1 \left( 1 - \dfrac{u_i^*}{K_i} \right) + \dfrac{\theta_0 u_i^*}{\alpha_i} - \eta_i p_{i}(u_i^*) + \eta_{i-1} \left[ p_{i-1}(u_{i-1}^*) + 2\lambda_{i-1}u_{i-1}^*e^{-\lambda_{i-1}u_{i-1}^*} \right].
\end{align*}
\endgroup
Since $u_i^* > K_i$, the negativity of the above quantity follows given the following condition:
\begin{align}
    r_1 \left( 1 - \dfrac{u_i^*}{K_i} \right) &< -\dfrac{\theta_0 u_i^*}{\alpha_i} - \eta_{i-1} \left[ p_{i-1}(u_{i-1}^*) + 2 \lambda_{i-1}u_{i-1}^* e^{-\lambda{i-1}u_{i-1}^*} \right] \\[12pt]
    \iff K_i &< \dfrac{r_1 u_{i}^*}{r_1 + \left[ \dfrac{\theta_0 u_i^*}{\alpha_i} + \eta_{i-1} \left[ p_{i-1}(u_{i-1}^*) + 2 \lambda_{i-1}u_{i-1}^* e^{-\lambda{i-1}u_{i-1}^*} \right] \right]}. \label{K_i_importantcondition}
\end{align}
Biologically, this condition corresponds a sufficiently small carrying capacity of the lymph nodes.

Next, for $i \in \mathcal{I}$, let $k = 2i+2$. We have
\begingroup
\allowdisplaybreaks
\begin{align*}
    J_{kk}(E_u) + R_k &= \dfrac{\partial \mathcal{N}_i}{\partial n_i} \bigg|_{E_u}  + \left| \dfrac{\partial \mathcal{N}_i}{\partial u_i} \bigg|_{E_u}  \right| + \left| \dfrac{\partial \mathcal{N}_i}{\partial c_i} \bigg|_{E_u}  \right| + \sum_{j \in \mathcal{I} \setminus \{i \}} \left( \left| \dfrac{\partial \mathcal{N}_i }{\partial u_j} \bigg|_{E_u} \right| + \left| \dfrac{\partial \mathcal{N}_i }{\partial n_j} \bigg|_{E_u} \right| + \left| \dfrac{\partial \mathcal{N}_i }{\partial c_j} \bigg|_{E_u} \right| \right) \\[15pt]
   &= r_2 - \dfrac{r_2}{K_i}u_i^* + \dfrac{\theta_0 u_i^*}{\alpha_i} - \gamma_0 - \eta_i p_i(u_i^*) + \eta_{i-1} p_{i-1}(u_{i-1}^*).
\end{align*}
\endgroup
In order to ensure negativity of this quantity, it suffices to impose the condition
\begin{equation}
    \gamma_0 > \dfrac{\theta_0 u_i^*}{\alpha_i} + \eta_{i-1}p_{i-1}(u_{i-1}^*), \label{gamma0condition_regional}
\end{equation}
which is a condition for the local asymptotic stability of the tumour-dominant steady state which is similar to that of the local model in Section 3.

Finally, we consider the final node in the network, lymph node $\ell$. We have
\begingroup
\allowdisplaybreaks
\begin{align*}
    J_{(2\ell + 1)(2 \ell + 1)}(E_u) + R_{2\ell + 1} &=  \dfrac{\partial \mathcal{U}_\ell}{\partial u_\ell} \bigg|_{E_u}  + \left| \dfrac{\partial \mathcal{U}_\ell}{\partial n_\ell} \bigg|_{E_u}  \right| + \left| \dfrac{\partial \mathcal{U}_\ell}{\partial c_\ell} \bigg|_{E_u}  \right| +  \sum_{j=0}^{\ell-1} \left( \left| \dfrac{\partial \mathcal{U}_\ell }{\partial u_j} \bigg|_{E_u} \right| + \left| \dfrac{\partial \mathcal{U}_\ell }{\partial n_j} \bigg|_{E_u}  \right| +  \left| \dfrac{\partial \mathcal{U}_\ell }{\partial c_j} \bigg|_{E_u} \right| \right) \\[15pt]
    &= r_1  - \dfrac{2r_1}{K_\ell}u_{\ell}^* + \dfrac{r_1}{K_\ell}u_{\ell}^* + \dfrac{\theta_0 u_{\ell}^*}{\alpha_\ell} + \eta_{\ell-1} \left[   p_{\ell-1}(u_{\ell-1}^*) + 2\lambda_{\ell-1}u_{\ell-1}^*e^{-\lambda_{\ell-1}u_{\ell-1}^*} \right]  \\[15pt]
    &= r_1 \left( 1 - \dfrac{u_{\ell}^*}{K_\ell} \right) + \dfrac{\theta_0 u_{\ell}^*}{\alpha_\ell} + \eta_{\ell-1} \left[   p_{\ell-1}(u_{\ell-1}^*) + 2\lambda_{\ell-1}u_{\ell-1}^*e^{-\lambda_{\ell-1}u_{\ell-1}^*} \right].
\end{align*} 
\endgroup
It is clear that this quantity is negative if condition (\ref{K_i_importantcondition}) is satisfied for $i = \ell$. Next,
\begingroup
\allowdisplaybreaks
\begin{align*}
    J_{(2\ell + 2)(2 \ell + 2)}(E_u) + R_{2\ell + 2} &= \dfrac{\partial \mathcal{N}_\ell}{\partial n_\ell} \bigg|_{E_u}  + \left| \dfrac{\partial \mathcal{N}_\ell}{\partial u_\ell} \bigg|_{E_u}  \right| + \left| \dfrac{\partial \mathcal{N}_\ell}{\partial c_\ell} \bigg|_{E_u}  \right| +  \sum_{j=0}^{\ell-1} \left( \left| \dfrac{\partial \mathcal{N}_\ell }{\partial u_j} \bigg|_{E_u} \right| + \left| \dfrac{\partial \mathcal{N}_\ell}{\partial n_j} \bigg|_{E_u} \right| + \left| \dfrac{\partial \mathcal{N}_\ell }{\partial c_j} \bigg|_{E_u} \right| \right) \\[15pt]
    &= r_2 - \dfrac{r_2}{K_\ell}u_{\ell}^* + \dfrac{\theta_0 u_{\ell}^*}{\alpha_{\ell}} - \gamma_0 + \eta_{\ell - 1} p_{\ell - 1}(u_{\ell-1}^*).
\end{align*} 
\endgroup
It is again clear that this quantity is negative if condition (\ref{gamma0condition_regional}) is satisfied for $i = \ell$.

Finally, define $\tilde{\mathcal{I}} := \{ 0, 1, 2, \dots, \ell\}$. For $i \in \tilde{\mathcal{I}}$, we define $k:= 2\ell + 3 + i$. It follows that
\begingroup
\allowdisplaybreaks
\begin{align*}
    J_{kk}(E_u) + R_K &= \dfrac{\partial \mathcal{C}_i}{\partial c_i} \bigg|_{E_u}  + \left| \dfrac{\partial \mathcal{C}_i}{\partial u_i} \bigg|_{E_u}  \right| + \left| \dfrac{\partial \mathcal{C}_i}{\partial n_i} \bigg|_{E_u}  \right| +  \sum_{j \in \tilde{\mathcal{I}}\setminus \{i\}} \left( \left| \dfrac{\partial \mathcal{C}_i }{\partial u_j} \bigg|_{E_u} \right| + \left| \dfrac{\partial \mathcal{C}_i }{\partial n_j} \bigg|_{E_u} \right| + \left| \dfrac{\partial \mathcal{C}_i}{\partial c_j} \bigg|_{E_u} \right| \right) \\[15pt]
    &= -\beta - q_1 u_i^*.
\end{align*}
\endgroup
Hence, $ J_{kk}(E_u) + R_K < 0$. By Lemma \ref{gershgorinlemma}, we conclude that if all of the above conditions are satisfied, then $E_u$ is locally asymptotically stable. We state this result in the following proposition.

\begin{proposition}
Consider system (\ref{eqlattice1}) - (\ref{eqlatticel3}) when $\phi_k = 0$, $q_{k,R} = 1$, and $q_{k,L} = 0$ for $k = 0, 1, 2, \dots, \ell$. The tumour-dominant steady state $E_u$ is locally asymptotically stable if the following conditions are satisfied:
\begin{enumerate}
    \item $$\eta_0 > \dfrac{1}{p_0 (u_0^*)} \left[ r_1 \left(1 - \dfrac{u_0^*}{K_0} \right) + \dfrac{\theta_0 u_o^*}{\alpha_0} \right]. $$
    \item For $i = 1, 2, \dots, \ell$, $$ K_i < \min \left\{ \dfrac{10 \eta_{i-1}}{7\eta_i}u_{i-1}^* \left( 1 - e^{-\lambda_{i-1}u_{i-1}^*} \right), \dfrac{r_1 u_{i}^*}{r_1 + \left[ \dfrac{\theta_0 u_i^*}{\alpha_i} + \eta_{i-1} \left[ p_{i-1}(u_{i-1}^*) + 2 \lambda_{i-1}u_{i-1}^* e^{-\lambda{i-1}u_{i-1}^*} \right] \right]} \right\}.$$
    \item For $i = 1, 2, \dots, \ell$, $$\gamma_0 > \dfrac{\theta_0 u_i^*}{\alpha_i} + \eta_{i-1}p_{i-1}(u_{i-1}^*).$$
\end{enumerate}
\label{regional_propchap3}
\end{proposition}
Proposition \ref{regional_propchap3} has some significant biological implications. Since this proposition gives conditions for the stability of the tumour-dominant steady state, the conditions being satisfied represents a clinically unfavourable outcome. The condition $\phi_k = 0$ represents no external oxygenation. Similarly to the local model, we see that hypoxic environments are beneficial to the tumour cells and reduce the efficacy of the adenovirus. Condition 1 of the proposition represents a sufficiently large rate of spreading of tumour cells away from the primary tumour. Condition 2 represents smaller carrying capacities of the lymph nodes -- this is not surprising, as tumour cells will more easily spread away from lymph nodes with lesser carrying capacities, i.e., due to less available resources. Condition 3 once again mirrors an important insight from the local model -- the oncolysis rate must not be too large in relation to the infection rate for an OV to be effective. However, this condition now comes with the additional consideration of incoming tumour cells from the previous lymph node in the network. In general, in a clinical setting, the model suggests effective treatment with an OV requires the engineering of a virus with a sufficiently large infection rate under hypoxic environments, which takes into account the spreading speed of the tumour cells as well as the carrying capacities of the lymph nodes. We further explore the implications of the regional model in the next section.

\section{Numerical simulations: regional model}

Due to the lack of analytic tractability of system (\ref{eqlattice1}) - (\ref{eqlatticel3}), we perform simulations to investigate the dynamics of this system. The primary tumour parameters (except for $\theta$ and $\gamma$) are pulled directly from Table 2. The parameters $K_i$ and $\alpha_i$ are estimated by taking into account the corresponding tumour parameters, $K_0$ and $\alpha_0$. In particular, for all $i$, we take $\eta_i =  0.0002$ days$^{-1}$, $K_i = K_0/10$ and $\alpha_i = \alpha_0/10$. Note that these parameters are the same for all lymph nodes. We make the biologically reasonable assumption that cells have a higher probability of migrating away from the primary tumour, i.e., in the direction of increasing node index. Hence, we set $q_{1,L} = q_{2,L} = q_{3,L} = 0.05$ and $q_{1,R} = q_{2,R} = 0.95$. 

\begin{figure}[h]
\centering
	\begin{subfigure}[t]{0.45\textwidth}\centering
\includegraphics[width=\textwidth]{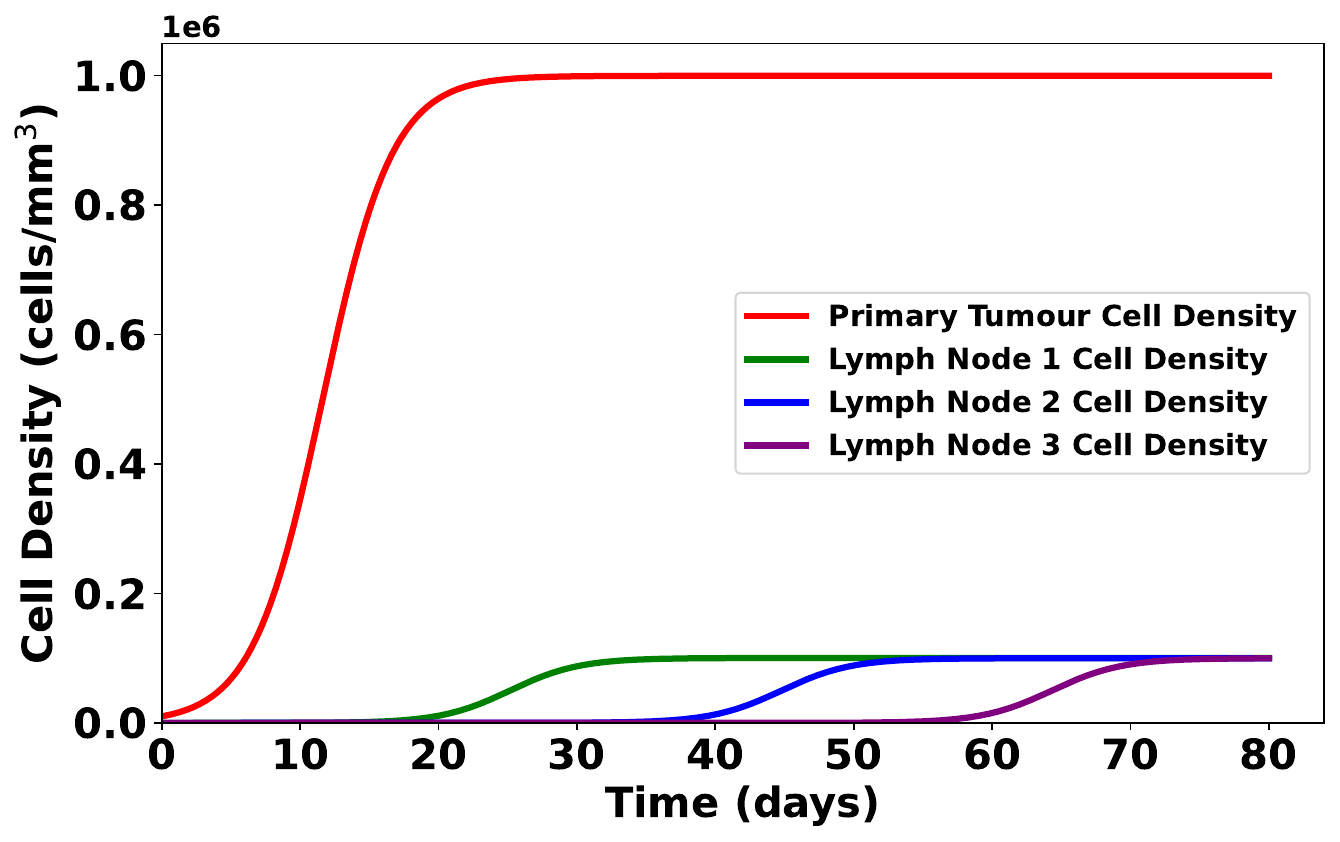}
         \caption{}\label{fig_lattice1a}
    \end{subfigure}
    \hspace{5mm}
\begin{subfigure}[t]{0.478\textwidth}\centering
\includegraphics[width=\textwidth]{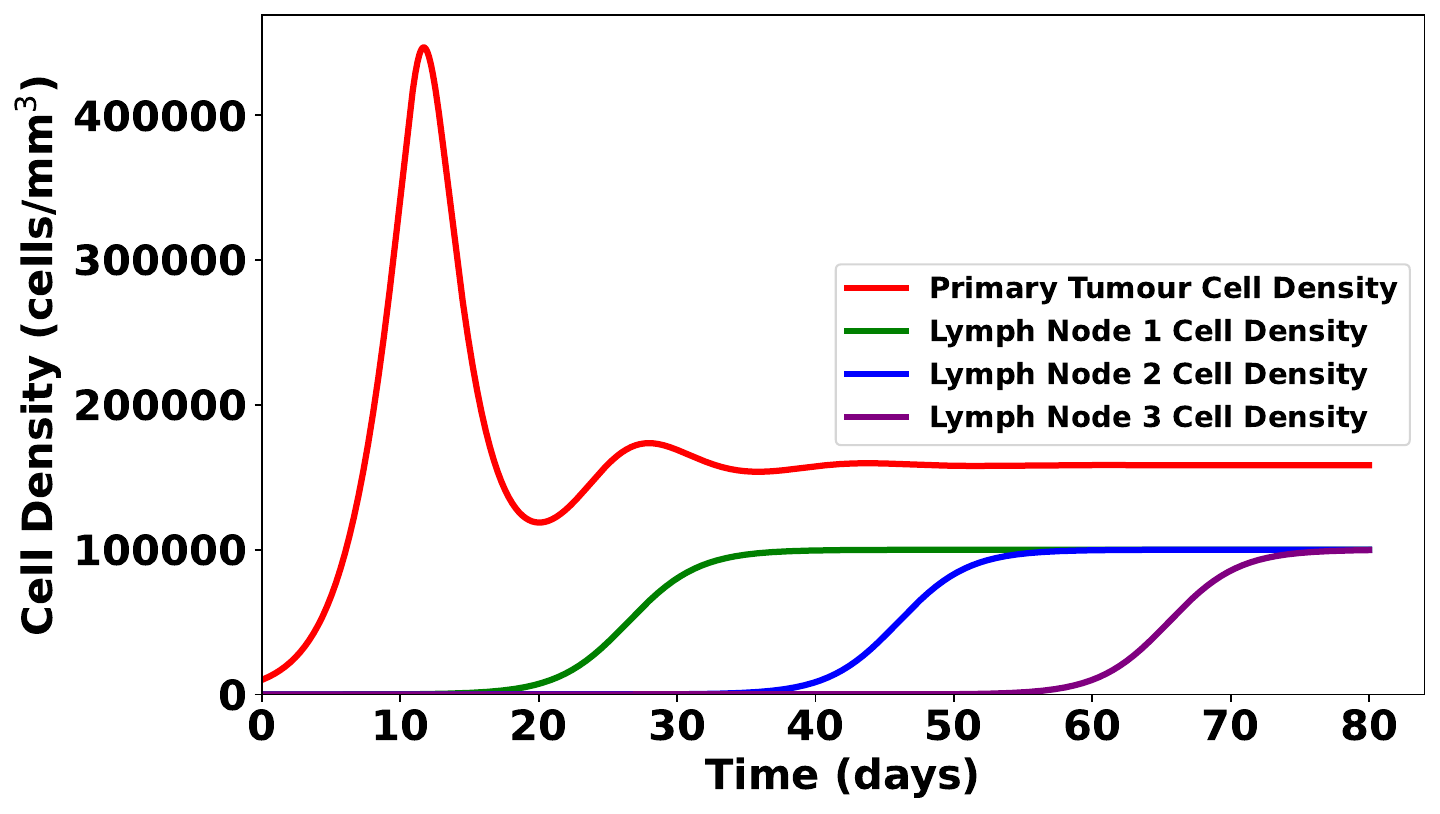}\\
 \caption{}\label{fig_lattice_1b}
\end{subfigure}
 \caption{The impact of oxygenation rate of the primary tumour, $\phi_0$, on the tumour cell density at the primary site and the first three lymph nodes in a network. (a)  The case of no external oxygenation, $\phi_0 = 0$. (b) The case where $\phi_0 = 10^4$ mM day$^{-1}$. There is a marked reduction in long-term tumour cell density when oxygenation is increased.}
\label{3graph10}
\end{figure}

In Figure \ref{3graph10}, we compare the case of no external oxygen input, $\phi_0 = 0$, to the case of some external oxygen input, $\phi_0 = 10^4$ mM day$^{-1}$. We graph the total number of tumour cells, $u_i (t) + n_i(t)$ over the course of 80 days. The functions $\theta(c)$ and $\gamma(c)$ are given by equations (\ref{theta_gamma_equations1}), where $\theta_0 = 0.005115, \theta_\infty =  1.0, k_\theta = 0.08, \gamma_0 = 0.1, \gamma_\infty = 0.9,$ and $k_\gamma = 0.08$. These parameter values are similar to the ones used in Figure \ref{graph04} -- they yield a favourable clinical outcome in the local model.  The model assumes that external oxygenation may only be performed on the primary tumour site -- not at the lymph nodes. From Figure \ref{3graph10} (a), we see that in the case where no external oxygen is provided, the tumour cells ultimately dominate at the primary tumour site and also approach a value near the carrying capacity at the lymph nodes. This unfavourable result is in stark contrast to the results of Figure \ref{graph04}, in which the tumour cells are either eradicated or kept under control. On the other hand, in the case of external oxygenation seen in Figure \ref{3graph10} (b), there is a sharp drop in the total tumour cell density. Namely, from a peak value approaching the carrying capacity at the primary site to approximately $4.47 \times 10^5$ cells/mm$^3$. This is a result of the benefit which the OV acquires as a result of an oxygen rich environment. This is consistent with the benefit consistently seen when treating cancer in oxygen-sufficient tumour microenvironments compared to hypoxic microenvironments. Even though the oxygenation occurs only at the primary tumour sites, the model allows for the proliferation of infected tumour cells through the lymphatic vessels into the lymph nodes, and hence the oxygenation also confers an increase in the efficacy of the OV treatment at the lymph nodes.

To this end, we turn our attention to the behaviour at the lymph nodes. In Figure \ref{3graph10} (a), the total tumour cell density across all three lymph nodes in the long-term is approximately given by the sum of their carrying capacities. In Figure \ref{3graph10} (b), as a result of external oxygenation, it takes a longer period of time for the tumour cell densities at the lymph nodes to reach their carrying capacities. This is because the benefit of the oxygenation here is less direct -- the oxygenation is only occurring at the primary site. There is still an indirect benefit, however, as a marked decrease of tumour cells at the primary site will result in slower spreading rates.

In summary, Figure \ref{3graph10} further illustrates the importance of the oxygen concentration in treatment with adenoviruses, a result which is consistent with the existing oncology literature \cite{onyx_adeno1}. It may also be worth noting that in contrast to Figure \ref{3graph10} (a), the tumour cell density at lymph node 3 eventually dominates the tumour cell density at lymph node 1 in Figure \ref{3graph10} (b). This may be explained by the fact that oxygenation occurs at the primary tumour site and, therefore, the infected cells are initially closer to the lymph nodes closer to the primary site rather than the subsequent lymph nodes in the network. Hence, lymph node 1 has a slightly greater benefit from the OV treatment than do lymph nodes 2 and 3. 

From the local model, we found that having a lower virus-induced death rate compared to the infection rate tends to yield more favourable clinical outcomes. To this end, we investigate the dynamics of the regional model in the case where $\theta(c) > (\alpha/K)\gamma(c)$ for all $c \geq 0$. We set $\theta_0 = 0.05115, \theta_\infty = 2.115, k_\theta = 0.016$ and $\gamma(c) = 0.005115$ for all $c\geq 0$. We once again plot the cases $\phi_0 = 0$ and $\phi_0 = 10^4$ mm day$^{-1}$.

\begin{figure}[h]
\centering
	\begin{subfigure}[t]{0.45\textwidth}\centering
\includegraphics[width=\textwidth]{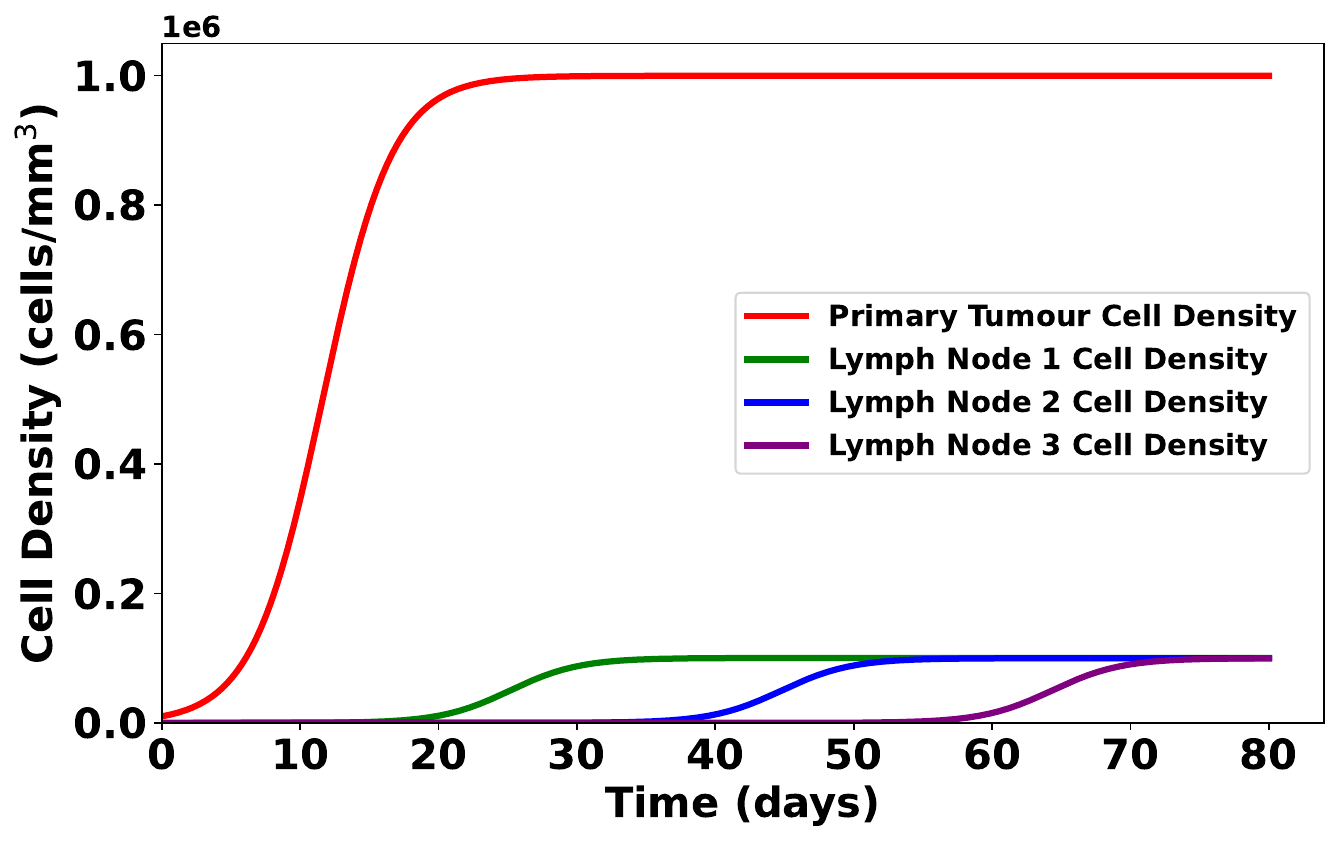}
         \caption{}\label{fig_lattice1a}
    \end{subfigure}
    \hspace{5mm}
\begin{subfigure}[t]{0.478\textwidth}\centering
\includegraphics[width=\textwidth]{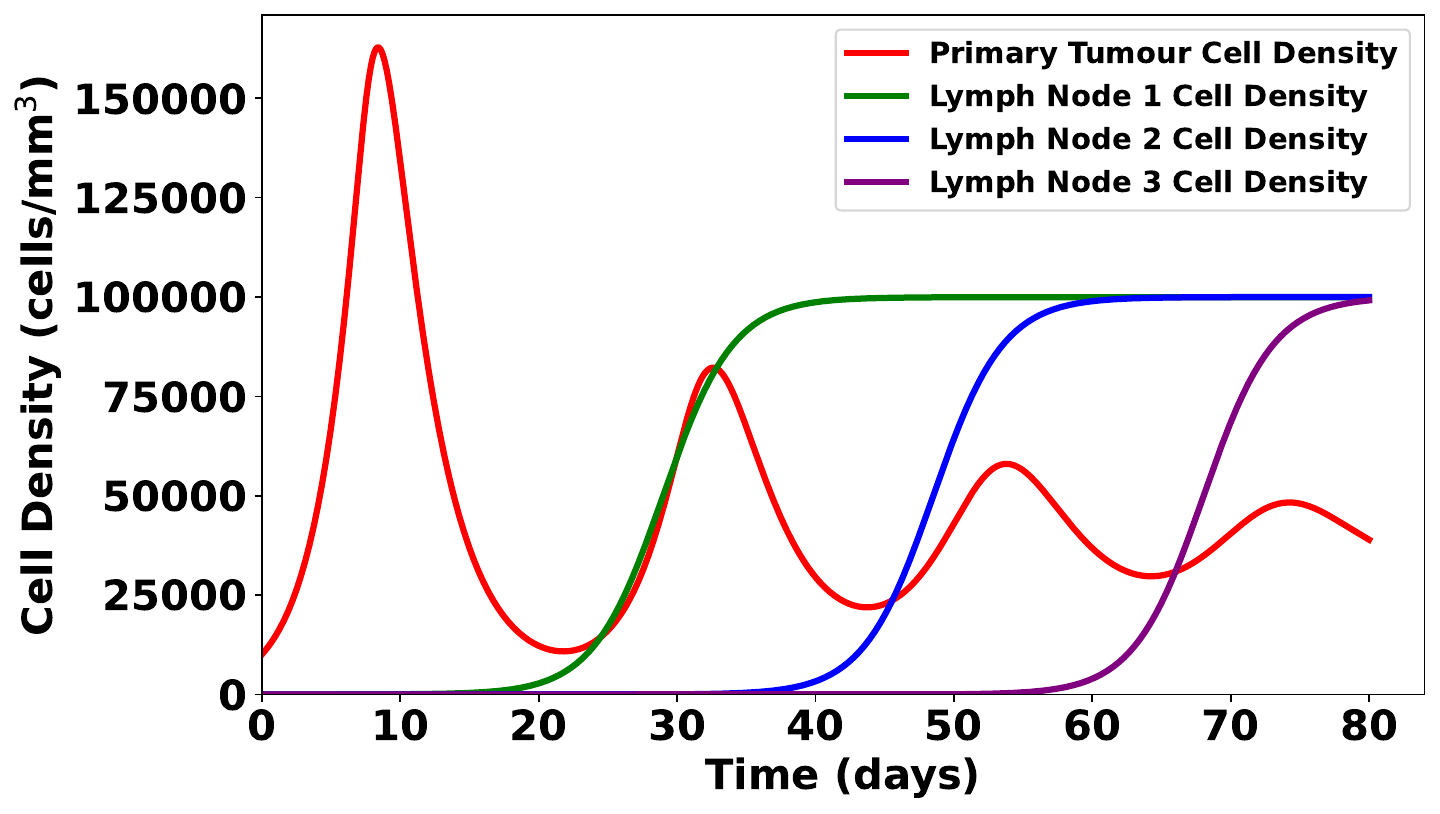}\\
 \caption{}\label{fig_lattice_1b}
\end{subfigure}
 \caption{Dynamics of the regional model in a case where $\theta(c) > (\alpha/K)\gamma(c)$. (a)  The case of no external oxygenation, $\phi_0 = 0$. (b) The case where $\phi_0 = 10^4$ mM day$^{-1}$. There is a very sharp reduction in long-term tumour cell density at the site of the primary tumour but the impact on the lymph nodes is much less pronounced.}
\label{3graph11}
\end{figure}

Figure \ref{3graph11} shows that the impact of having a sufficiently low virus oncolysis rate in the regional model is consistent with the local model. Once again, the effect of an increased external oxygenation rate is much more pronounced at the primary tumour compared to the lymph nodes.

\begin{figure}[H]
    \centering
    \includegraphics[scale=0.6]{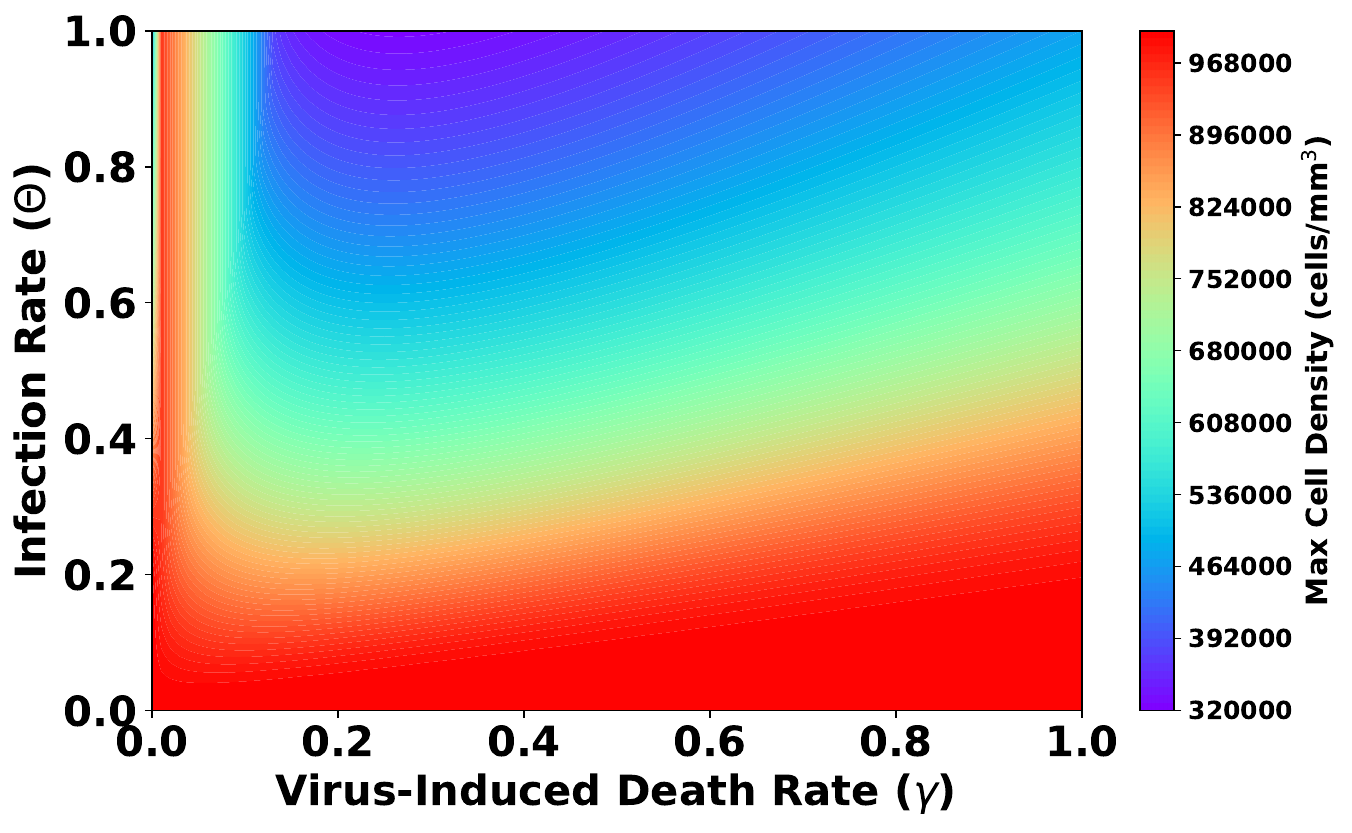}
    \caption{Maximum tumour cell density at the location of the primary tumor over the course of 80 days after treatment for various values of $\theta$ and $\gamma$.}
    \label{chap3heatmap}
\end{figure}

Motivated by Figure \ref{3graph11}, we now consider the impact of the infection rate, $\theta$, and the virus-induced death rate, $\gamma$, on the regional model. These parameters were considered extensively in the the numerical simulations of the local model in Section 4. In this case, we consider keeping $\theta$ and $\gamma$ constant rather than as  functions of oxygen concentration. The results of the simulations are plotted in Figure \ref{chap3heatmap}. 

In Figure \ref{chap3heatmap}, we plot the maximum value of the tumour cell density at the primary tumour site over the course of 80 days after OV treatment. That is, we plot $\max \{ u_0(t) + n_0(t) \}$ for different values of $\theta$ and $\gamma$. Consistent with our prior results, we can again visualize the relationship between infection and oncolysis. We see that increasing the virus-induced death rate to a much greater value relative to the infection rate leads to an unfavourable outcome (red region). This also occurs if the virus-induced death rate is too small, regardless of the value of the infection rate. Therefore, this provides further evidence of the importance of a high infection rate and a  oncolysis rate that is \textit{not too low} in oder to obtain favourable results (blue region).

Finally, we address the case where $\theta$ and $\gamma$ depend on the oxygen concentration. In particular, we assume that we have some mechanism through which to administer external oxygen to the lymph nodes and set $\phi_k = 10^4$ mM day$^{-1}$ for $k = 1, 2, \dots, \ell$. 

\begin{figure}[H]
\centering
	\begin{subfigure}[t]{0.47\textwidth}\centering
\includegraphics[width=\textwidth]{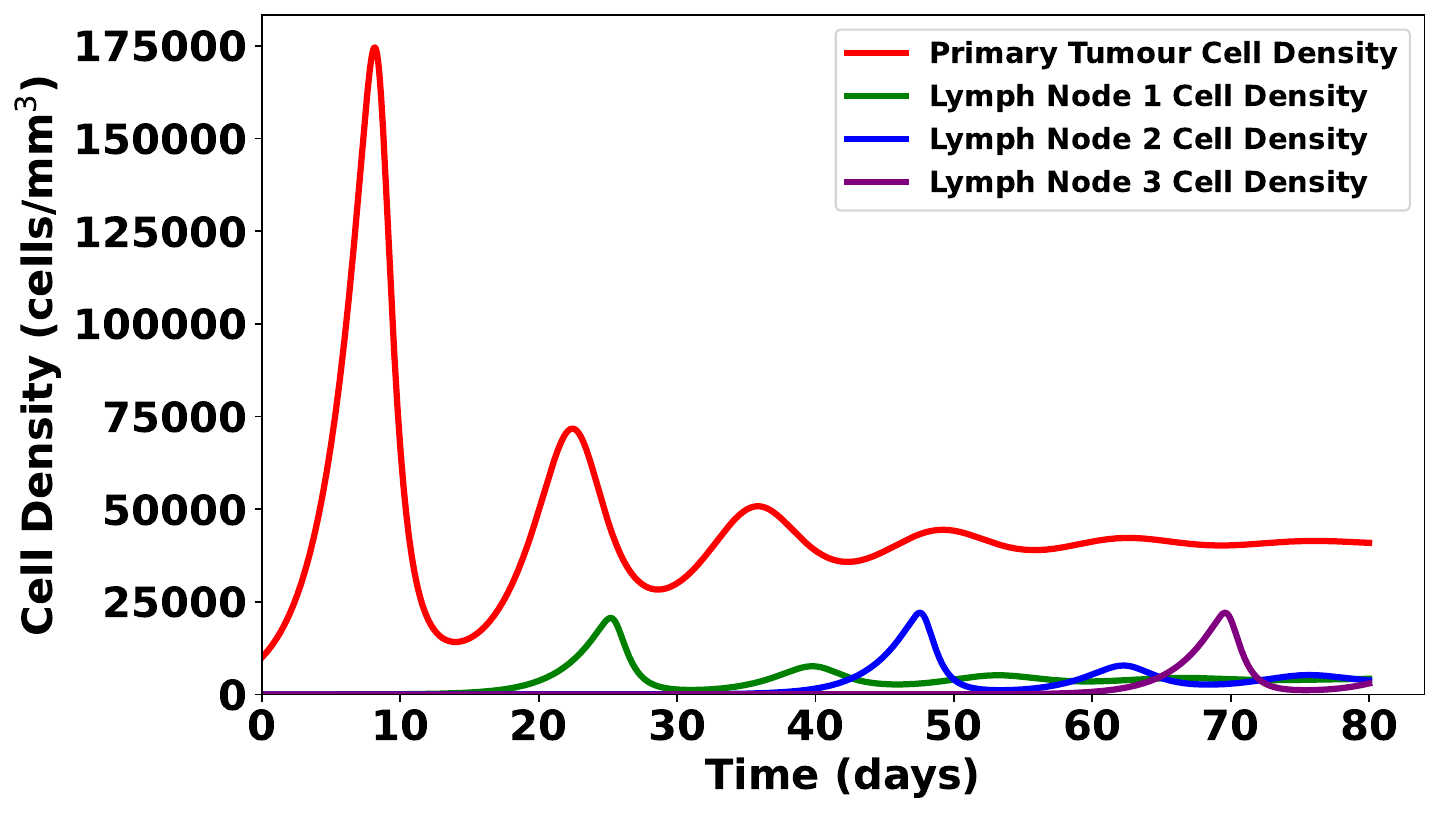}
         \caption{}\label{fig_lattice1a}
    \end{subfigure}
    \hspace{5mm}
\begin{subfigure}[t]{0.47\textwidth}\centering
\includegraphics[width=\textwidth]{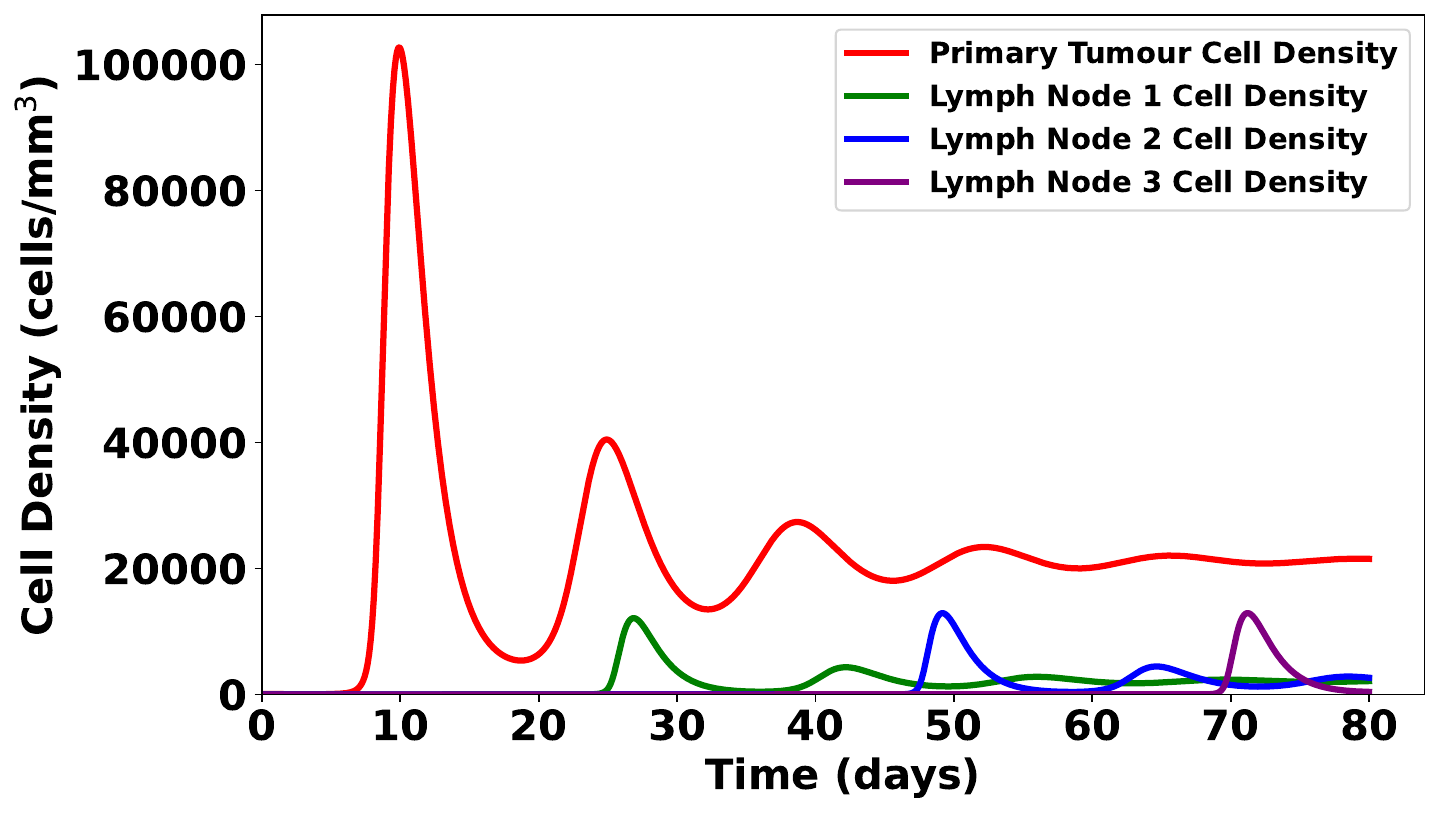}\\
 \caption{}\label{fig_lattice_1b}
\end{subfigure}
    \caption{The impact of oxygenation rate at the primary tumour site and the first three lymph nodes in the network. (a) Uninfected tumour cell density at the primary tumour and at each node. (b) Infected tumour cell density at the primary tumour and at each node}
\label{chap3_everynodeoxygenated}
\end{figure}

Figure \ref{chap3_everynodeoxygenated} shows the case of oxygen dependence of the infection rate and the oncolysis rate. In this case, there is a marked reduction in the total tumour cell density in all compartments. Figure \ref{chap3_everynodeoxygenated} (a) shows the uninfected tumour cell density, $u_i$ and Figure \ref{chap3_everynodeoxygenated} (b) shows the infected tumour cell density, $n_i$. In contrast to the case of no oxygen input at the lymph nodes, the infected tumour cell density at each lymph node asymptotically approaches a value below the carrying capacity of its corresponding node. This provides further evidence which supports the lack of efficacy of oncolytic adenoviruses in hypoxic environments and the increased efficacy of these OVs when external oxygenation is provided.

\section{Conclusion and discussion}

From the mathematical results of this work, as well as the simulations, the importance of the functions $\theta$ and $\gamma$ are emphasized. Biologically, this refers to the interplay between viral infection rate and the virus-induced death rate of the cancer cells. If the virus-induced cancer cell death rate is too large compared to the infection rate, the cancer cells end-up dominating in the long-run. This is a reflection of the virus not being able to infect cells faster than the infected cells are destroyed. On the other-hand, if the infection rate of the OV is significantly  large compared to the virus-induced death rate in \textit{all} oxygen environments, the infected tumour cells will dominate in the long-run and will reach some steady state. This steady state may represent the case where we avoid uncontrollable cancer cell growth as the number of cells will not approach the carrying capacity. This translates to a favourable clinical result. On the other hand, it may also represent a state in which the infected tumour cells dominate \textit{at} the carrying capacity if the OV tumour-destroying capabilities are \textit{too} low. For this reason, we suggest that when engineering OVs, it is important to ensure that these viruses have greater infection capabilities than they have oncolytic capabilities while ensuring that the virus-induced death rate is not too low. In particular, our results suggest that maintaining high viral infection rates tends to lead to clinically favourable results regardless of oxygen concentration of the tumour microenvironment. {Two approaches which have been identified for the production of efficient OVs under hypoxic conditions are direct genetic engineering and directed evolution \cite{onyx_adeno1}. While successes from taking a direct engineering approach have not yielded consistently favourable results, the approach of using directed evolution may offer a solution to creating potent OVs. While beyond the scope of this paper, more details on this approach can be found in \cite{directed_evolution}. } Our findings on the importance of the infection rate are consistent with \cite{Jenner, more_evidence} and offer further insight on a growing body of literature regarding efficacy of engineered viruses \cite{Jenner, more_evidence,  Kuro}. 

Another important component of this paper is the modelling of the impact of hypoxic conditions on OV treatment efficacy. As previously stated, the modelling suggests that significantly high infection rates are preferable under any oxygen conditions. However, another layer of complexity is added when considering the threat of toxicity which the OV poses toward healthy cells \cite{simpson}. Furthermore, having a virus-induced death rate which is \textit{too low} will lead to a decreased mortality of cancer cells. Hence, it is not sufficient to simply conclude that engineering extremely infectious viruses is the solution. Instead, we proposed taking into account the effect of different oxygen conditions and hypoxia when constructing the OVs. To address this, we considered the case in which which function dominates, $\theta$ or $\gamma$, depends on the oxygen concentration. A favourable result occurs when $\gamma$ dominates for low oxygen concentrations but $\theta$ dominates for high oxygen concentrations. Hence, we conjecture that another consideration of engineering OVs is whether or not the tumour microenvironment is hypoxic. The preferential virus characteristic would be to have greater oncolytic capabilities in hypoxic environments and greater infectious capabilities in more oxygen-rich environments. According to the modelling, this may lead to stability of a steady tumour load rather than uncontrollable growth. However, we also found that making the virus-induced death rate \textit{too} great under hypoxic conditions also leads to a reduction in the efficacy of the treatment, as the infected tumour cells die faster than they may infect the remaining susceptible tumour cells. 

We extended the model to a regional model which incorporated spatial structure through considering the axillary lymph nodes. This was done by considering ODEs on a one-dimensional lattice. This natural extension captures the invasive nature of melanoma (and many other invasive cancers). Once again, the importance of considering oxygen cannot be understated. When considering a system with three lymph nodes, we found that providing oxygenation at the site of the primary lesion (through an external oxygen source) yields an approximately 72\% decrease in tumour cell density at the site of the primary lesion. Lymph nodes closer to the site of oxygenation similarly obtained benefit from more  hyperoxic conditions. This benefit of external oxygenation in (various forms of) the treatment of cancer has also been observed clinically, such as in the use of hyperbaric-oxygen therapy \cite{tibbles}. Our simulations further support these experimental findings. We also found that the impact of the infection rate, $\theta$, is also present in the regional model and the findings were consistent with those of the local model. This leads us to further stress the importance of oxygen rich microenvironments being used in tandem with highly infectious OVs.

This model may be further enhanced by the addition of a variable which accounts for the free virus particles. Although this would increase the complexity of system in terms of mathematical analysis, it would lead to more interesting dynamics, biologically. In terms of the parameters, the growth rate of tumour cells also depends on the available oxygen of the tumour microenvironment \cite{hypoxia_last}. Hence an important next step is the use of growth rates which depend on the oxygen concentration, i.e., $r_1 (c), r_2(c)$. Future work also includes adding a \textit{continuous} spatial structure to the model, i.e., through the use of PDE modelling. This can take into account the spatial properties of the tumour as well as the efficacy of OV treatment in the context of metastatic disease by modelling cancer cell spreading at the site of the primary lesion. Extending the types of geometry of the lattice representing the lymphatic network is also an important next step. For example, this involves allowing certain lymph nodes in the network to have connections with multiple neighbouring lymph nodes. From a clinical perspective, incorporating the use of conventional chemotherapy along with the virotherapy is also likely to provide potentially useful insights. Finally, the toxic effects of an increased tumour cell infection rate may also be worth considering in order to model a more comprehensive treatment approach. 

\section*{Acknowledgements}

This research is partially supported by the Natural Sciences and Engineering Research Council of Canada (NSERC) and the province of Ontario via the Ontario Graduate Scholarship (OGS). The authors would also like to thank the two anonymous reviewers for their helpful comments, particularly for helping better shape this paper from a biological perspective. 

\section*{Conflict of interest}

The authors declare there is no conflict of interest.

\section{Appendix}

The code used to plot the solutions of both the local ((\ref{3eq1}) - (\ref{3eq3})) and regional (system (\ref{eqlattice1}) - (\ref{eqlatticel3})) models is given below. To obtain plots for the local model, we may set the spreading rate of tumour cells away from the primary tumour, $\eta_0$, to $0$. For non-negative $\eta_i$ values, the code produces plots which include lymph node involvement. 

\begin{lstlisting}[language=Python]
import numpy as np
import matplotlib.pyplot as plt
from scipy.integrate import odeint

#Parameters:

r1, r2 = 0.3954, 0.21;
K, alpha = 1e6, 1e5;
phi, beta, q1, q2 = 1e4, 5.0976, 5.47e-5, 0.5*(5.47e-5);
theta0, thetainf, k_theta = 0.005115, 2.115, 0.016
gamma0, gammainf, k_gamma = 0.1, 0.9, 0.08;

K1, K2, K3 = K/10, K/10, K/10;
alpha1, alpha2, alpha3 = alpha/10, alpha/10, alpha/10;
K_values = [K,K1,K2,K3];

eta = 0.0002; #Comment this out if eta not constant.
eta_values = [eta for i in range(4)] #The case with 3 lymph nodes.

#Functions:
def theta(x):
    
    return thetainf*theta0/(theta0 + (thetainf - theta0)*( np.exp((-1)*k_theta*x)))

def gamma(x):
    
    return gammainf*gamma0/(gamma0 + (gammainf - gamma0)*( np.exp((-1)*k_gamma*x) ))

def p(i,x):
    
    Lambda = ((-1)*np.log(0.3))/K_values[i] #start at i = 0
    
    return 1 - np.exp((-1)*Lambda*x)

def ODEs(x,t): 
    
    u, n, c, u1, n1, c1 = x[0], x[1], x[2], x[3], x[4], x[5]
    u2, n2, c2 = x[6], x[7], x[8]
    u3, n3, c3 = x[9], x[10], x[11]
    
    dudt = r1*u*(1 - (u+n)/K) - ((theta(c))*n*u)/(alpha + n) - eta_values[0]*u*p(0,u+n) + 0.05*eta_values[1]*u1*p(1,u1+n1)
    dndt = r2*n*(1 - (u+n)/K) + ((theta(c))*n*u)/(alpha + n) - (gamma(c))*n - eta_values[0]*n*p(0,u+n) + 0.05*eta_values[1]*n1*p(1,u1+n1)
    dcdt = phi - beta*c - q1*u*c - q2*n*c
    
    du1dt = r1*u1*(1 - (u1+n1)/K1) - ((theta(c1))*n1*u1)/(alpha1 + n1) - eta_values[1]*u1*p(1,u1+n1) +  eta_values[0]*u*p(0,u+n)  + 0.05*eta_values[2]*u2*p(2,u2+n2)
    dn1dt = r2*n1*(1 - (u1+n1)/K1) + ((theta(c1))*n1*u1)/(alpha1 + n1) - (gamma(c1))*n1 - eta_values[1]*n1*p(1,u1+n1) + eta_values[0]*n*p(0,u+n) + 0.05*eta_values[2]*n2*p(2,u2 + n2)
    dc1dt =  (-1)*beta*c1 - q1*u1*c1 - q2*n1*c1
    
    du2dt = r1*u2*(1 - (u2+n2)/K2) - ((theta(c2))*n2*u2)/(alpha2 + n2) - eta_values[2]*u2*p(2,u2+n2) +  0.95*eta_values[1]*u1*p(1,u1+n1) + 0.05*eta_values[3]*u3*p(3,u3+n3)
    dn2dt = r2*n2*(1 - (u2+n2)/K2) + ((theta(c2))*n2*u2)/(alpha2 + n2) - (gamma(c2))*n2 - eta_values[2]*n2*p(2,u2+n2) + 0.95*eta_values[1]*n1*p(1,u1+n1) + 0.05*eta_values[3]*n3*p(3,u3+n3) 
    dc2dt =  (-1)*beta*c2 - q1*u2*c2 - q2*n2*c2
    
    du3dt = r1*u3*(1 - (u3+n3)/K3) - ((theta(c3))*n3*u3)/(alpha3 + n3) - 0.05*eta_values[3]*u3*p(3,u3+n3) +  0.95*eta_values[2]*u2*p(2,u2+n2)
    dn3dt = r2*n3*(1 - (u3+n3)/K3) + ((theta(c3))*n3*u3)/(alpha3 + n3) - (gamma(c3))*n3 - 0.05*eta_values[3]*n3*p(3,u3+n3) + 0.95*eta_values[2]*n2*p(2,u2+n2) 
    dc3dt =  (-1)*beta*c3 - q1*u3*c3 - q2*n3*c3
    
    return [dudt, dndt, dcdt, du1dt, dn1dt,dc1dt, du2dt, dn2dt, dc2dt, du3dt, dn3dt, dc3dt]

#Initial conditions:
u0, n0, c0, u10, n10, c10 = 10000,100, 4.3751, 0, 0, 4.375;
u20, n20, c20 = 0, 0,  4.375;
u30, n30, c30 = 0, 0,  4.375;
init_0 = [u0, n0, c0, u10, n10, c10, u20, n20, c20, u30, n30, c30];

#Numerically solving and plotting the solution of the regional model:
t = np.linspace(0,80,10000);#domain

x = odeint(ODEs, init_0,t);  #integrating
        
u, n, c = x[:,0], x[:,1], x[:,2];
u1, n1, c1 = x[:,3], x[:,4], x[:,5];
u2, n2, c2 = x[:,6], x[:,7], x[:,8];
u3, n3, c3 = x[:,9], x[:,10],x[:,11];

plt.plot(t,u+n,'red',label='Primary Tumour',linewidth=3);
plt.plot(t,u1+n1,'green',label='Lymph Node 1',linewidth=3);
plt.plot(t,u2+n2,'blue',label='Lymph Node 2',linewidth=3);
plt.plot(t,u3+n3,'purple',label='Lymph Node 3',linewidth=3);
plt.xlim(0)
plt.ylim(0)
plt.legend(('Primary Tumour Cell Density', 'Lymph Node 1 Cell Density', 'Lymph Node 2 Cell Density', 'Lymph Node 3 Cell Density'),
              loc='upper right')
plt.ylabel("Cell Density (cells/mm$^3$)")
plt.xlabel("Time (days)")

    \end{lstlisting}
    
The following code produces the heatmap in Figure \ref{chap3heatmap}. This multi-parametric analysis shows the peak tumour density value for various values of constant $\theta$ and $\gamma$ over an interval of 100 days after initial OV treatment is administered. 

 \begin{lstlisting}[language=Python]   
t_val = 100; #Solve over this interval.
N = 100; #N+1 values of theta and gamma used.
theta_values = [0.01*i for i in range(0,N+1)];
gamma_values = [0.01*i for i in range(0,N+1)];
max_cancer_cells = [[] for i in range(0,N+1)];
    
i = 0;
for j in theta_values:
    theta0 = j;
    thetainf = j;
        
    for k in gamma_values:
        gamma0 = k;
        gammainf = k;
         
        x = odeint(ODEs, init_0,t);
    
        u = x[:,0];
        n = x[:,1];
        c = x[:,2];
            
        max_cancer_cells[i].append((max(u+n)));
            
    i = i + 1;
    
plt.xlabel("Virus-Induced Death Rate ($\gamma$)")
plt.ylabel("Infection Rate ($\Theta$)")
    
img = plt.contourf(theta_values,gamma_values,max_cancer_cells,100,cmap='rainbow')
plt.colorbar(img)
    
    \end{lstlisting}

\end{document}